\documentclass[final,leqno]{siamltex704}

\usepackage{amsmath,amssymb}% for \begin{bmatrix}, \begin{smallmatrix}, \begin{align} etc.
\usepackage{amsfonts}% for \mathbb{} etc. `blackboard' for uppercase only.
\usepackage{exscale}
\usepackage{graphicx}
\usepackage{comment}
\usepackage{subcaption}
\usepackage[ruled]{algorithm2e}
\usepackage{algorithmic}
\usepackage{url}
\usepackage{appendix}

\usepackage[pdftex]{thumbpdf}      %%% thumbnails for pdflatex. Run thumbpdf at the end
\usepackage[pdftex,                %%% hyper-references for pdflatex
pdfstartview=FitH,%
bookmarks=true,%                   %%% generate bookmarks ...
bookmarksnumbered=true,%           %%% ... with numbers
bookmarksopen=true,%
hypertexnames=false,%              %%% needed for correct links to figures !!!
breaklinks=true,%                  %%% break links if exceeding a single line
  colorlinks=true,%
  linkcolor=blue,anchorcolor=blue,%
  citecolor=blue,filecolor=blue,%
  menucolor=blue]%
{hyperref} 
\usepackage{etoolbox}

\makeatletter
\def\@cite#1#2{{\rm [}{{\rm#1}\if@tempswa , #2\fi}{\rm ]}}
\makeatother

\title{Preconditioned Locally Harmonic Residual Method for Computing Interior
Eigenpairs of Certain Classes of Hermitian Matrices
\thanks{%
Preliminary version posted at \url{http://arxiv.org/}.
The results presented in this 
work are partially based on the PhD thesis of the first coauthor~\cite{thesis}.
}
}

\author{Eugene Vecharynski\thanks{Computational Research Division, Lawrence Berkeley National Laboratory, Berkeley, CA 94720
(eugene.vecharynski[at]gmail.com)} \and Andrew Knyazev\thanks{
Mitsubishi Electric Research Laboratories; 201 Broadway
Cambridge, MA 02139 (knyazev[at]merl.com)}}

\begin{document}
%siam_id=72574
%CODEN=SJMAEL
%\slugger{sisc}{?}{?}{?}{?--?}
\maketitle

\setcounter{page}{1}

\begin{abstract}
We propose a Preconditioned Locally Harmonic Residual (PLHR) method for
computing several interior eigenpairs of a generalized Hermitian eigenvalue problem,
without traditional spectral transformations, matrix factorizations, or inversions.
PLHR is based on a short-term recurrence, easily extended to a block form, computing eigenpairs simultaneously.
PLHR 
%relies on the
can take advantage of Hermitian positive definite preconditioning, e.g.,
based on an approximate inverse of an absolute value of a shifted matrix,
introduced in [SISC, 35 (2013), pp. A696--A718].
Our numerical experiments demonstrate that
%PLHR is efficient and robust for large-scale interior eigenvalue computations,
PLHR is efficient and robust for certain classes of large-scale interior eigenvalue problems,
involving Laplacian and Hamiltonian operators, 
especially if memory requirements are tight.
\end{abstract}

\begin{keywords}
Eigenvalue, eigenvector, Hermitian, absolute value preconditioning, linear systems
%Preconditioned eigenvalue solvers, absolute value preconditioning, linear systems
% Preconditioning, linear system, preconditioned minimal residual method, 
% polar decomposition, matrix absolute value, geometric multigrid
\end{keywords}

%\begin{AMS}
%TBD %15A06, 65F08, 65F10, 65N22, 65N55 
%\end{AMS}
% 15A06 -- Linear equations 
% 65F08 -- Preconditioners for iterative methods
% 65F10 -- Iterative methods for linear systems
% 65N22 -- Solution of discretized equations
% 65N55 -- Multigrid methods; domain decomposition

%\begin{DOI}
% 
%\end{DOI}

\pagestyle{myheadings}
\thispagestyle{plain}
\markboth{EUGENE VECHARYNSKI AND ANDREW KNYAZEV}{COMPUTING INTERIOR EIGENPAIRS OF HERMITIAN MATRICES} %50 Characters Limit
%\markboth{EUGENE VECHARYNSKI AND ANDREW KNYAZEV}{COMPUTING INTERIOR EIGENPAIRS OF HERMITIAN PROBLEMS} %50 Characters Limit

\section{Introduction}\label{sec:intro}

We are interested in computing a subset of eigenpairs of the generalized 
Hermitian eigenvalue problem
\begin{equation}\label{eqn:eigp}
A v = \lambda B v, \quad A = A^* \in \mathbb{C}^{n \times n}, \quad B = B^*>0 \in \mathbb{C}^{n \times n}, 
%A v = \lambda B v, 
%\quad A = A^* \in \mathbb{R}^{n \times n}, \quad B = B^*>0 \in \mathbb{R}^{n \times n}, 
\end{equation} 
that correspond to eigenvalues closest to a given shift $\sigma$, which is real and points to the eigenvalues in the interior of the spectrum.  
We refer to~\eqref{eqn:eigp} as an \textit{interior} eigenproblem,
and call the targeted eigenpairs \textit{interior eigenpairs}. 

Interior eigenpairs are of fundamental interest in some physical models, e.g.,\ 
for first-principles electronic structure analysis 
of materials using a semi-empirical potential or a charge patching 
method~\cite{Wang.Li:04, Wang.Zunger:95, Wang.Zunger:99}.
There, the requested eigenvalues correspond to energy levels around a material-dependent
energy shift $\sigma$, and the eigenvectors represent the associated wave functions. 

We assume that the size $n$ of the matrices is so large that eigenvalue solvers based on explicit matrix transformations are impractical.  
Other conventional ways of solving the interior eigenvalue problems are based on
spectral transformations that allow reducing the interior eigenproblem to an easier task
of computing extreme (smallest or largest) eigenvalues of a transformed pencil.  
Such transformations include the so-called \textit{shift-and-invert} (SI) and \textit{folded spectrum} (FS) 
approaches; see, e.g.,~\cite{templates:00}. While both techniques are extensively used, 
they face serious issues as $n$ increases. 

In particular, SI relies on repeated solution of large shifted linear systems,
which generally becomes extremely inefficient or infeasible for large $n$.
Similar difficulties are experienced in methods relying on contour integration, e.g.,\ \cite{sakurai2007}.
 
The FS approach for \textit{standard}, i.e.,\ with identity $B=I$, eigenproblems is inverse-free. 
However, in exact arithmetic, it is equivalent to solving an eigenproblem for the matrix $(A - \sigma I)^2$, which  
%it 
squares the condition number, resulting in severely tightened clustering of the targeted
eigenvalues. For the \textit{generalized}, i.e.,\ with $B\neq I$, eigenproblems, FS requires inverting
the matrix $B$, which can be problematic. 
Additionally, available preconditioning, normally based on approximating
an inverse of $A - \sigma B$, may not be applicable for FS, which should be preconditioned by an approximate
inverse of the squared operator.     

%%%%%%%%%REFEREE
For certain classes of eigenproblems, multigrid (MG) algorithms have proved to be very efficient; see, 
e.g.,~\cite{Borzi.Borzi:05, Hetmaniuk:07, Kushnir.Galun.Brandt:10, Livshits:04} and references therein. 
However, the applicability of such eigensolvers faces the same kind of limitations 
as exist for the MG linear solvers in that either geometric (grid) information must be available at solution time
or the problem's structure should be appropriate for the \textit{algebraic} MG approaches.  

Another popular technique for interior eigenproblems is based on combining polynomial transformations, called \emph{filters},
with the standard Lanczos method, as in filtered Lanczos procedure~\cite{Fang.Saad:12}.
The approach cannot easily be extended to generalized
eigenproblems and does not take advantage of preconditioners that may be available.

Residual minimization scheme, direct inversion in the iterative subspace (RMM-DIIS), see~\cite{wood1985new}, 
provides an alternative, but is known to be unreliable if the initial approximations are not close enough to 
the targeted eigenpairs. RMM-DIIS has to be preceded by, e.g.,\  generalized Davidson~\cite{Morgan:91} or 
Jacobi--Davidson~\cite{Sleijpen.Vorst:96} methods, such as, e.g.,\ in Vienna \textit{ab initio} 
simulations (VASP) package~\cite{VASP}. 
Combining \textit{harmonic} Rayleigh--Ritz (RR) procedure ~\cite{templates:00} with the Davidson's
family of methods has recently aimed at supplementing RMM-DIIS in VASP; see~\cite{Jordan.Marsman.Kim.Kresse:12}. 

Success of the Davidson type methods is often determined by two interdependent algorithmic components:
the maximum allowed size of the search subspace and the quality of the preconditioner, typically given as a form 
of an approximate solve of the system with a shifted matrix $A - \sigma B$. 
The maximum size of the search subspace has to be increased if the preconditioner quality deteriorates.  
In practical computations, preconditioners that lack the desired quality are common, 
especially if the shift $\sigma$ targets the eigenvalues
that are deep in the interior of the spectrum. 

We develop new preconditioned iterative schemes for interior eigenproblems 
that exhibit a lower sensitivity to the preconditioner quality
and at the same time lead to fixed and relatively modest memory requirements.  
A connection between iterative methods for singular homogeneous Hermitian indefinite systems 
and eigenproblems in Section~\ref{sec:idealized} motivates our development.  
Our base method, called \textit{Preconditioned Locally Harmonic Residual} (PLHR) and
presented in Section~\ref{sec:plmr}, computes a single eigenpair. 
In Section~\ref{sec:bplhr}, we generalize PLHR to block iterations,
where targeted eigenpairs, as determined by $\sigma$, are computed simultaneously, 
similar to LOBPCG~\cite{Knyazev:01}.
 
PLHR requires that the preconditioner is  Hermitian positive definite (HPD). 
In Section~\ref{sec:prec}, we 
%assume 
suggest
that the preconditioner $T$ 
should
approximate 
the inverted matrix 
absolute value $|A - \sigma B|$. Preconditioners of this kind, called the absolute value (AV)
preconditioners, have been recently introduced in~\cite{Ve.Kn:13}. In particular,~\cite{Ve.Kn:13}
describes an MG scheme for approximating the inverse of the AV of the shifted Laplacian.
Thus, the approach can be directly applied as a PLHR preconditioner for computing interior eigenpairs
of the Laplacian matrix, which is demonstrated in our numerical experiments in Section~\ref{sec:numr}.      
Furthermore, the suitable HPD preconditioners are readily available in electronic structure calculations, 
where we combine PLHR with the existing state-of-the-art preconditioner~\cite{Teter.Payne.Allan:89}, 
also in Section~\ref{sec:numr}. 

While, in this work, we focus only on several model problems with already 
available HPD preconditioners, the proposed PLHR approach can be applicable to broader classes
of eigenproblems 
%provided that suitable preconditioners can be ?? 
that admit appropriate HPD preconditioning. 
Investigation of such problems, as well as the development of the 
corresponding preconditioners, is a matter of future research, and is 
% and, therefore, 
outside the scope of the current paper.

An important novelty of the proposed approach is a modification of the harmonic
RR procedure, obtained by a proper utilization of the HPD preconditioner 
$T$ in the Petrov-Galerkin condition for extracting the approximate eigenvectors,
including its real arithmetic implementation.
The new extraction technique, that we call \textit{the $T$-harmonic RR procedure} in  Section~\ref{sec:plmr}, is critical. 
Our tests demonstrate remarkable robustness if the $T$-harmonic procedure is used.
The gains are especially evident if the preconditioner quality deteriorates,
whereas the memory requirements remain fixed.

\section{Preconditioned null space computations}\label{sec:idealized}
We start with an easier problem of finding a null space component of a Hermitian 
matrix. Let $\lambda_q$ be a targeted eigenvalue of~\eqref{eqn:eigp} closest to 
the shift $\sigma$. We assume that $\lambda_q$ is known.
% and simple. 
Then eigenproblem~\eqref{eqn:eigp} turns into the singular homogeneous 
linear system    
\begin{equation}\label{eqn:eigp_ideal}
(A - \lambda_q B) v = 0.
\end{equation} 

%We also assume that the vector $v_q$ is normalized to have the unit $B$-norm, 
%and hence is unique up to a sign.
It is clear that (\ref{eqn:eigp_ideal}) has a non-trivial null space solution determining
an eigenvector $v_q$ associated with $\lambda_q$.
Thus, in the idealized setting, where $\lambda_q$ is available,
an eigensolver for~\eqref{eqn:eigp} 
can be given by an appropriate solution scheme for linear system~\eqref{eqn:eigp_ideal}. 
%
%For this reason, as a next step we consider possible approaches   
%that is capable of finding a null space solution.    
In this sense, linear solvers for~\eqref{eqn:eigp_ideal} can be viewed
as prototypical, or idealized, methods for computing the eigenpair $(\lambda_q, v_q)$.

The described connection between linear and eigenvalue solvers
has been emphasized in literature; see, e.g.,\ ~\cite{Knyazev:86, Knyazev:01}.
For example, in~\cite{Knyazev:01}, a special case is considered where $\lambda_q$ is the 
smallest eigenvalue of~\eqref{eqn:eigp} and, hence,
%, so that
the system~\eqref{eqn:eigp_ideal} is Hermitian positive semidefinite. 
A proper choice of the linear solver---a three-term recurrent form
of the preconditioned conjugate gradient (PCG) method---has 
led to derivation of the popular LOBPCG algorithm for finding extreme 
eigenpairs; see~\cite{Knyazev:01, Kn.Ar.La.Ov:07}.  
 
We follow a similar approach. Assuming that $\lambda_q$ is known,
we select efficient preconditioned linear solvers 
capable of computing a non-trivial solution of~\eqref{eqn:eigp_ideal}.  
Viewing these linear solvers as idealized eigensolvers, we then extend them
to the practical case where $\lambda_q$ is unknown. 
 
Since the targeted eigenvalue $\lambda_q$ can be located anywhere in the 
interior of the spectrum of the pencil $A - \lambda B$, the Hermitian coefficient matrix
of system~\eqref{eqn:eigp_ideal} is generally indefinite
in contrast to~\cite{Knyazev:01}, 
%the above mentioned case where $\lambda_q$ is the smallest eigenvalue, 
which makes PCG inapplicable. We look for suitable 
preconditioned short-term recurrent Krylov subspace type methods that can be applied to singular 
\textit{Hermitian indefinite} systems of the form~\eqref{eqn:eigp_ideal} and
guarantee linear convergence.   

An optimal technique in this class of methods is the preconditioned minimal
residual (PMINRES) algorithm~\cite{Paige.Saunders:75}. While most commonly applied to
non-singular systems, the method is also guaranteed to work for a class of 
\textit{singular consistent} Hermitian systems,
%~[{\color{blue} ???}], 
%which includes problem
such as~\eqref{eqn:eigp_ideal}. The preconditioner for PMINRES should normally be HPD, and 
the convergence of the method is governed by the spectrum distribution of the 
preconditioned matrix. 
However, the standard implementation of the algorithm is deeply rooted in the Lanczos procedure, 
where the availability of the matrix $A - \lambda_q B$ is crucial at every stage
of the process. Since this assumption is not realistic in the context of eigenvalue
computations, where only an approximation of $\lambda_q$ is at hand, PMINRES
fails to provide a proper insight into the structure of \textit{local subspaces} that 
determine a new approximate 
solution. Similar arguments apply to methods that are mathematically equivalent
to PMINRES, such as, e.g., orthodir($3$)~\cite{Greenbaum:97, Saad:03,Young.Jea:80}. 
%As we shall see in the next section, 
%capturing such local subspaces is an important ingredient in the derivation of an 
%locally optimal eigensolver.  
     
Returning back for a moment to the special case where $\lambda_q$ is the 
smallest eigenvalue of~\eqref{eqn:eigp}, a preconditioned steepest descent (PSD) 
linear solver can be viewed as a predecessor of the PCG linear solver, on the one hand, and as 
restarted preconditioned  Lanczos linear solver, on the other hand. 
We use this viewpoint to describe an analog of a PSD-like linear solver for the indefinite case.

To that end, we restart PMINRES in such a way that linear convergence is 
preserved (possibly with a lower rate), whereas the number of steps between the 
restarts is kept to a minimum possible. %This problem was extensively studied in~\cite{thesis}.  
It is shown in~\cite{thesis} that in order to maintain PMINRES convergence, 
the number of steps between the restarts should be no smaller than two.
Thus, the simplest convergent residual minimizing iteration for~\eqref{eqn:eigp_ideal} is
\begin{equation}\label{eqn:plmr_ideal_short}
%\begin{array}{lcl}
v^{(i+1)}  =  v^{(i)} + \alpha^{(i)} T r^{(i)} + 
%\beta^{(i)} T\left(A  - \lambda_q B \right) T\left(A  - \lambda_q B \right) v^{(i)},  \\
\beta^{(i)} T\left(A  - \lambda_q B \right) T r^{(i)},  \quad i = 0, 1, \ldots; \\
%w^{(i)} & = & T\left(Ax^{(i)} - \lambda_q B x^{(i)}\right)} = 0; \ i  = 0,1,\ldots,
%\end{array}
\end{equation}    
where $r^{(i)} = \left(\lambda_q B - A \right) v^{(i)}$ and 
the iteration parameters $\alpha^{(i)}$, $\beta^{(i)}$  
are chosen to minimize the $T$-norm of the residual $r^{(i+1)}$,
%over the affine space
i.e., are such that
\begin{equation}\label{eqn:coeff}
\|r^{(i+1)}\|_T = \min_{u \in  \mathcal{K}^{(i)}} 
\|r^{(i)} -  (A - \lambda_q B) u \|_T, 
%\mathcal{K}^{(i)} = \text{span}\left\{T r^{(i)}, T(A - \lambda_q B) T r^{(i)})\right\}.
%\arg \min_{\alpha \in \mathbb{R}} \|r^{(i)} - \alpha A l^{(i)}\|_T.
\end{equation} 
where 
\begin{equation}\label{eqn:subsp1}
\mathcal{K}^{(i)} = \text{span}\left\{T r^{(i)}, 
T(A - \lambda_q B) T r^{(i)} \right\}.
\end{equation}
%%%%%%%%%%%% REFEREE
Here and throughout, for a given HPD matrix $M$, the corresponding vector $M$-norm is defined 
as $\|\cdot\|_M = (\cdot, M \cdot)^{1/2}$.
%
%For simplicity, let us assume that $\lambda_q$ is a simple eigenvalue. 
\begin{theorem}[\cite{thesis}]\label{thm:ideal_cv}
Given an HPD preconditioner $T$, the iteration~\eqref{eqn:plmr_ideal_short}--\eqref{eqn:subsp1}
converges to a nontrivial solution of~(\ref{eqn:eigp_ideal}), 
provided that the initial guess $x^{(0)}$ has a non-zero projection onto the null space of $A - \lambda_q B$.
%provided that the 
%initial guess $x^{(0)}$ has a nonzero component from the null space of $A - \lambda_q B$. 
%in the expansion using the basis of the eigenvectors of the pencil $A - \lambda B$. 
Moreover, if $\lambda_q$ is a simple eigenvalue\footnote{This assumption is made only to simplify the statement of the theorem. The 
theorem can be similarly formulated for the case of multiple eigenvalue.} 
of $A - \lambda B$ and 
$\mu_1 \leq \mu_2 \leq \ldots \leq \mu_{q-1} < \mu_q = 0 < \mu_{q+1} < \ldots \mu_n$
are the eigenvalues of the preconditioned matrix $T(A - \lambda_q B)$, 
then the residual norm reduction is given~by 
\begin{equation}\label{eqn:cv_ideal_nonopt}
\frac{\| r^{(i+1)} \|_T}{\|  r^{(i)}\|_T} \leq \frac{\tilde{\kappa} - 1}{\tilde{\kappa} + 1} < 1,
\end{equation}
where 
%assuming that $\mu_1 < \mu_{q-1} < \mu_q = 0 < \mu_{q+1} < \mu_n$ 
%are the nonzero eigenvalues of the
%preconditioned operator $T(A - \lambda_q B)$, the expression for $\tilde \kappa$ is given by 
\begin{equation}\label{eqn:cv_ideal_kappa}
\tilde{\kappa} =
\left\{ \begin{array}{l} 
\displaystyle \left( \frac{\mu_n}{\mu_{q+1}} \right) \left( 1 + \frac{\mu_n - \mu_{q+1}}{\left|\mu_{q-1}\right|} \right), 
\ \mbox{if} \  \left|\mu_1\right| - \left|\mu_{q-1}\right| \leq \mu_{n} - \mu_{q+1} \\
 \\
\displaystyle \left(\frac{\mu_1}{\mu_{q-1}}\right)\left(1 + \frac{\left|\mu_1\right| - \left|\mu_{q-1}\right|}{\mu_{q+1}}\right), \ \mbox{if} \  \left|\mu_1\right| - \left|\mu_{q-1}\right| > \mu_n - \mu_{q+1}.   
\end{array} \right. 
\end{equation}
\end{theorem}

%%%%%%%REFEREE
The proof of Theorem~\ref{thm:ideal_cv} relies on the analysis of iterations of the 
form~\eqref{eqn:plmr_ideal_short}, where the parameters $\alpha^{(i)}$ and $\beta^{(i)}$ are fixed.
In this case,~\eqref{eqn:plmr_ideal_short} can be viewed as a stationary Richardson type scheme for a polynomially
preconditioned system, and a standard convergence analysis (see, e.g.,~\cite[Theorem 5.6]{Axelsson:94}) applies to determine
the values of the parameters that yield an optimal convergence, which is given by~\eqref{eqn:cv_ideal_nonopt}--\eqref{eqn:cv_ideal_kappa}.
Since the convergence of the stationary iteration cannot be better then that of~\eqref{eqn:plmr_ideal_short}
with the residual minimizing parameters~\eqref{eqn:coeff}, bound~\eqref{eqn:cv_ideal_nonopt}--\eqref{eqn:cv_ideal_kappa}
immediately applies to method~\eqref{eqn:plmr_ideal_short}--\eqref{eqn:subsp1}. For more details, we refer the reader to~\cite{thesis}. 
%%%%%%%

A natural way to  enhance~\eqref{eqn:plmr_ideal_short}--\eqref{eqn:subsp1}
is by introducing an additional term that holds information from the previous step.  
By analogy with a three-term recurrent form of PCG, 
%we define a vector $v^{(i)} - v^{(i-1)}$ ($v^{(-1)} = 0$) and 
consider the scheme
\begin{equation}\label{eqn:plmr_ideal}
%\begin{array}{lcl}
v^{(i+1)}  =  v^{(i)} + \alpha^{(i)} T r^{(i)} + 
\beta^{(i)} T\left(A  - \lambda_q B \right) T r^{(i)} 
+ \gamma^{(i)}(v^{(i)} - v^{(i-1)}), \ i  = 0,1,\ldots ,  \\
%w^{(i)} & = & T\left(Ax^{(i)} - \lambda_q B x^{(i)}\right)} = 0; \ i  = 0,1,\ldots,
%\end{array}
\end{equation}    
where $v^{(-1)} = 0$.
Here, the scalar parameters $\alpha^{(i)}$, $\beta^{(i)}$, and $\gamma^{(i)}$  
are chosen according to local minimality condition~\eqref{eqn:coeff}
with 
%\begin{equation}\label{eqn:coeff_ideal}
%\|r^{(i+1)}\|_T = \min_{u \in \mathcal{V}^{(i)}} 
%\|r^{(i)} + u \|_T,
%%\arg \min_{\alpha \in \mathbb{R}} \|r^{(i)} - \alpha A l^{(i)}\|_T.
%\end{equation} 
\begin{equation}\label{eqn:subsp2}
\mathcal{K}^{(i)} = \text{span}\left\{T r^{(i)}, 
T(A - \lambda_q B) T r^{(i)}, v^{(i)} - v^{(i-1)} \right\}.
\end{equation}
%$\mathcal{U}^{(i)} = \text{span}\left\{(A-\lambda_q B )T r^{(i)}, 
%(A-\lambda_q B )T(A - \lambda_q B) T r^{(i)}, (A - \lambda_q B)(v^{(i)} - v^{(i-1)})\right\}$.
%where $\mathcal{V}^{(i)} = \mathcal{U}^{(i)} + 
%\text{span}\left\{(A-\lambda_q B ) (v^{(i)} - v^{(i-1)}) \right\}$ and $v^{(-1)} = 0$. 
%
% It is clear that~the method is equivalent to PMINRES restarted every 
% two step with the~vector~$v^{(i)} - v^{(i-1)}$. ??????

Since $\mathcal{K}^{(i)}$ in~\eqref{eqn:subsp2} contains the subspace in~\eqref{eqn:subsp1}, 
the employed minimality condition~\eqref{eqn:coeff} implies that 
convergence of~\eqref{eqn:plmr_ideal} is not worse than that
of~\eqref{eqn:plmr_ideal_short}.
Hence, the results of Theorem~\ref{thm:ideal_cv} also apply to~\eqref{eqn:plmr_ideal}
with~\eqref{eqn:coeff} and~\eqref{eqn:subsp2}.   
In this case, however, bound~\eqref{eqn:cv_ideal_nonopt}--\eqref{eqn:cv_ideal_kappa}
is likely to be an overestimate, and in practice the presence of the additional
vector in the recurrence leads to a faster convergence.
In fact, the scheme~\eqref{eqn:plmr_ideal} has been observed
to exhibit behavior similar to PMINRES up to the occurrence of
superlinear convergence~\cite{thesis}, which is a consequence of the 
\textit{global} optimality of the latter. 
%(as opposed to only \text{local}
%optimality of~\eqref{eqn:plmr_ideal}--\eqref{eqn:coeff_ideal}).  
Additionally, iteration~\eqref{eqn:plmr_ideal} 
reveals the structure of local subspaces that are used to determine the improved
approximate solution, which we exploit in the next section for constructing the
trial subspaces in the context of the eigenvalue calculations.

For these reasons,
we choose~\eqref{eqn:plmr_ideal} with~\eqref{eqn:coeff} and~\eqref{eqn:subsp2} 
to be a \textit{``base''} linear solver for the null space problem~\eqref{eqn:eigp_ideal} 
and, in what follows, use it as a starting point for deriving preconditioned 
interior eigensolvers.

%Therefore,
%the choice of (\ref{eqn:plmr_ideal}) as a \textit{base} idealized method
%for finding the eigenpair $(\lambda_q,v_q)$ can be viewed as a reasonable
%compromise. Scheme (\ref{eqn:plmr_ideal}) is no longer \textit{globally}
%optimal. 
%However, it reveals the structure of the local subspaces, 
%used to determine the next iterate, 
%and is convergent for any initial guess. Moreover, the convergence
%rate and the amount of computational work can be expected to 
%mimic those (up to the possible occurrence of effects apparently attributed to the 
%superlinear convergence, see, e.g., Figure \ref{fig:4}) 
%of the 
%globally optimal preconditioned minimal residual method
%in one of its robust short-term recurrent formulations.   

\section{The Preconditioned Locally Harmonic Residual method}\label{sec:plmr}

We now present an approach for computing an eigenpair $(\lambda_q, v_q)$ of~\eqref{eqn:eigp} that corresponds to the eigenvalue
closest to a given shift $\sigma$. 
The method is motivated by the preconditioned null space finders discussed in the previous section.
%In particular, 
Our idea is to extend the ``base'' linear solver~\eqref{eqn:plmr_ideal} to the case
of eigenvalue computations by introducing a series of approximations into 
recurrence~\eqref{eqn:plmr_ideal} and optimality 
condition~\eqref{eqn:coeff}. 
%(with $\mathcal{K}^{(i)}$ in~\eqref{eqn:subsp2}).
%to the ``base'' linear solver~\eqref{eqn:plmr_ideal}
%Our derivation is motivated
%by the requirement that the resulting method should as much as possible 
%resemble the ``base'' linear solver~\eqref{eqn:plmr_ideal} 
%in terms of the convergence behavior and the associated computational cost.
%
%We achieve this by properly approximating the terms 
%in the recurrence~\eqref{eqn:plmr_ideal} as well as the optimality 
%condition~\eqref{eqn:coeff} with $\mathcal{K}^{(i)}$ in~\eqref{eqn:subsp2}. 
%As we demonstrate next, 
%the former approximation suggests the choice of the trial (search) subspaces, 
%whereas the latter gives an appropriate projection procedure for the eigenpair extraction. 
%%which we call the \textit{$T$-harmonic} projection.    

The null space finding scheme~\eqref{eqn:plmr_ideal}
suggests that, if $\lambda_q$ is known, the improved
eigenvector approximation $v^{(i+1)}$ belongs to the subspace
\begin{equation}\label{eqn:plmr_ideal_span}
%\[
\mbox{span}\left\{ v^{(i)}, T (A - \lambda_q B) v^{(i)}, T(A  - \lambda_q B ) T (A - \lambda_q B) v^{(i)}, v^{(i-1)}\right\}.
%\]
\end{equation}   
Clearly, in practice, the exact value of $\lambda_q$ is unavailable and, hence, the computation
of~\eqref{eqn:plmr_ideal_span} cannot be performed.
In order to obtain a computable subspace in the context of eigenvalue problem~\eqref{eqn:eigp}, 
we approximate $\lambda_q$ by the Rayleigh Quotient (RQ) 
\[
\lambda^{(i)} \equiv \lambda(v^{(i)}) = (v^{(i)},Av^{(i)})/(v^{(i)},Bv^{(i)}).
\] 
As a result, at each iteration $i$, the following \textit{trial subspace} is introduced:
\begin{equation}\label{eqn:plmr_span}
\mathcal{Z}^{(i)} = \mbox{span}\left\{ v^{(i)}, w^{(i)}, s^{(i)}, v^{(i-1)}\right\},
\end{equation}
where $w^{(i)} = T(A v^{(i)} - \lambda^{(i)} v^{(i)})$ is the \textit{preconditioned residual}
for problem~\eqref{eqn:eigp}, and $s^{(i)} = T(A w^{(i)} - \lambda^{(i)} w^{(i)})$. 
%Our expectation is that $\mathcal{Z}^{(i)}$ captures well the idealized trial 
%subspace~\eqref{eqn:plmr_ideal_span} suggested by the linear solver, allowing to resemble its
%convergence behavior and computational costs, provided a proper eigenvector extraction procedure is defined.
%~\eqref{eqn:plmr_ideal}. 

Next, we address the question of extracting an approximate 
eigenpair from the subspace $\mathcal{Z}^{(i)}$.
Let us recall that minimality principle~\eqref{eqn:coeff}, utilized by the base null space finder,
implies the orthogonality relation
\begin{equation}\label{eqn:orth_ideal}
r^{(i+1)} \equiv (A - \lambda_q B) v^{(i+1)} \perp_T (A - \lambda_q B) \mathcal{K}^{(i)},
\end{equation}
where ``$\perp_T$'' denotes the ``orthogonality in the $T$-based inner product'' 
and $\mathcal{K}^{(i)}$ is defined in~\eqref{eqn:subsp2}.
Following the analogy with the linear solver, in the context of the eigenvalue problem, 
we introduce a similar condition:   
\begin{equation}\label{eqn:orth}
(A  - \theta B) v^{(i+1)} \perp_T (A - \sigma B) \mathcal{Z}^{(i)}, \quad \|v^{(i+1)}\|_B = 1. 
%\quad \theta \in \mathbb{C}, \ v^{(i+1)} \neq 0. 
%\quad \|v^{(i+1)}\|_B = 1. 
\end{equation}
Here, we find a $B$-unit vector $v^{(i+1)}$ and a scalar $\theta \in \mathbb{C}$ that
ensure the orthogonality of the vector $(A  - \theta B) v^{(i+1)}$ to the 
%low-dimensional 
subspace $(A - \sigma B) \mathcal{Z}^{(i)}$. 

Condition~\eqref{eqn:orth} has been obtained from~\eqref{eqn:orth_ideal} as a result of two 
modifications. First, the known eigenvalue $\lambda_q$ in the residual $r^{(i+1)}$ of linear 
system~\eqref{eqn:eigp_ideal} has been replaced by an unknown $\theta$, giving rise to the 
residual like vector $(A  - \theta B) v^{(i+1)}$ for the eigenvalue problem. Note that
this vector should be distinguished from the standard eigenresidual $(A  - \lambda^{(i+1)} B) v^{(i+1)}$,
where $\lambda^{(i+1)}$ is the RQ. Since the orthogonality of $(A  - \theta B) v^{(i+1)}$  
in~\eqref{eqn:orth} is invariant with respect to the norm of $v^{(i+1)}$, for definiteness, we request that
the vector has a unit $B$-norm. 

The second modification has been the replacement of the subspace 
$(A - \lambda_q B) \mathcal{K}^{(i)}$ in the right-hand side of~\eqref{eqn:orth_ideal}  
by $(A - \sigma B) \mathcal{Z}^{(i)}$, where $\mathcal{Z}^{(i)}$ is the trial subspace defined 
in~\eqref{eqn:plmr_span}. We expect that $(A - \sigma B) \mathcal{Z}^{(i)}$ captures well
the ideal subspace $(A - \lambda_q B) \mathcal{K}^{(i)}$. 
Note that in the case where $\lambda_q$ is known, $\mathcal{K}^{(i)} \subset \mathcal{Z}^{(i)}$.
Hence, if $\sigma$ is close to the targeted eigenvalue, %at least near the solution $(\lambda_q, v_q)$, 
both subspaces are close to each other.  
%It is easy to see that~\eqref{eqn:orth} has been obtained from~\eqref{eqn:orth_ideal} 
%by formally replacing the residual $r^{(i+1)}$ of linear system~\eqref{eqn:eigp_ideal} by 
%the residual of the eigenvalue problem, 
%in the left-hand side, and approximating 
%$(A - \lambda_q B) \mathcal{K}^{(i)}$ by $(A - \sigma B) \mathcal{Z}^{(i)}$, in the right-hand side.
%We also impose the condition that the new approximate eigenevector has a unit $B$-norm. This is done
%only to ensure that the constructed sequence of 
%
%Note that we replace $\mathcal{K}^{(i)}$ by the trial subspace $\mathcal{Z}^{(i)}$ defined in~\eqref{eqn:plmr_span}, so that 
%$(A - \sigma B) \mathcal{Z}^{(i)}$ is expected to capture $(A - \lambda_q B) \mathcal{K}^{(i)}$.
%in the absence of the exact value~$\lambda_q$. 

%In contrast to $r^{(i+1)} = (A - \lambda_q B)v^{(i+1)}$ in~\eqref{eqn:orth_ideal}, the residual of the eigenvalue problem 
%contains \textit{a pair of unknowns} ($\theta, v^{(i+1)}$).   
%%Thus, given a trial subspace $\mathcal{Z}^{(i)}$ in~\eqref{eqn:plmr_span},
%%at every iteration $i$, we would like to define a new approximate eigenvector
%%%$v^{(i+1)} \in \mathcal{Z}^{(i)}$ and a corresponding eigenvalue approximation 
%%$\theta$ that satisfy~\eqref{eqn:orth}. 
%%and give a better approximation to the desired eigenpair $(\lambda_q, v_q)$. 
We next discuss a procedure for computing the pair $(\theta, v^{(i+1)})$ in~~\eqref{eqn:orth}. 
%%the desired pair 
%%$(\tilde \lambda, v^{(i+1)})$.   
 
\subsection{The $T$-harmonic Rayleigh--Ritz procedure}\label{subsec:THarm} 
Motivated by~\eqref{eqn:orth}, let us consider a general problem, 
where we are interested in finding a number of pairs $(\theta, v)$ that approximate 
eigenpairs of~\eqref{eqn:eigp} corresponding to the eigenvalues closest a given shift $\sigma$, 
such that each $v$ belongs to a given $m$-dimensional subspace $\mathcal{Z}$, and  
\begin{equation}\label{eqn:orth_gen}
(A  - \theta B) v \perp_T (A - \sigma B) \mathcal{Z}, \quad \|v\|_B = 1.
\end{equation}
%where $(\theta,  v)$ approximates eigenpairs
%associated with the eigenvalues closest to the given shift $\sigma$.

One can immediately recognize that~\eqref{eqn:orth_gen} represents the 
Petrov-Galerkin condition~\cite{templates:00, Saad-book3} formulated with respect to the $T$-based inner product.
In the standard case with $T = I$, this condition leads to the well known 
\textit{harmonic} Rayleigh-Ritz (RR) procedure, where approximate eigenpairs
are delivered by the \textit{harmonic} 
Ritz pairs~\cite{Morgan:91, Morgan.Zeng:98, Paige.Parlett.Vorst:95}.  
The corresponding vectors $A v - \theta B v$ are sometimes called the \textit{harmonic}
residuals; see, e.g.,~\cite{Vomel:10}.

Let a basis of $\mathcal{Z}$ be given by columns of an $n$-by-$m$ 
matrix $Z$, so that any element of the subspace can be expressed in the 
form $v = Z y$, where $y$ is a vector of coefficients. 
Then the orthogonality constraint~\eqref{eqn:orth_gen} translates into 
%the equation   
\[
Z^* (A - \sigma B) T (A - \theta B) Z y = \textbf{0}, \quad \|Z y\|_B = 1.
\]  
This is equivalent to the eigenvalue problem 
\begin{equation}\label{eqn:THarm_evp}
Z^* (A - \sigma B) T (A - \sigma B ) Z y  = \xi Z^* (A - \sigma B) T B Z y, 
%\quad \xi = \theta - \sigma, 
\quad \|Z y\|_B = 1,
\end{equation}  
where $\xi \equiv \theta - \sigma$.

It is clear that the smallest in the absolute value eigenvalues $\xi$ correspond
to the values of $\theta$ closest to the shift $\sigma$. Thus, if $(\xi, y)$ are
the eigenpairs of~\eqref{eqn:THarm_evp} associated with the eigenvalues of
the smallest absolute value, then the candidate approximations to the desired
eigenpairs of the original problem~\eqref{eqn:eigp} can be chosen as $(\theta, Zy)$.    
\begin{definition}
Let $(\xi, y)$ be an eigenpair of~\eqref{eqn:THarm_evp}. 
We call~$(\theta, Zy)$ a~\emph{$T$-harmonic Ritz pair}, where $\theta = \xi+\sigma$ 
and $v = Zy$ are a \emph{$T$-harmonic Ritz value and vector}, respectively.
The corresponding vector $Av - \theta B v$ is then referred to as the $T$-harmonic
residual with respect to the subspace $\mathcal{Z}$ (or, simply, the $T$-harmonic
residual).  
\end{definition}

In order to see why $T$-harmonic Ritz pairs can be expected to deliver 
reasonable approximations to the interior eigenpairs, let us introduce  
the substitution $Q = T^{\frac{1}{2}} (A - \sigma B) Z$ and rewrite~\eqref{eqn:THarm_evp} as     
\begin{equation}\label{eqn:THarm_evp_sinv}
Q^* (T^{\frac{1}{2}} B) (A - \sigma B)^{-1} B (T^{\frac{1}{2}} B)^{-1} Q y = \tau Q^* Q y, 
\quad \tau  = \frac{1}{\xi} \equiv \frac{1}{\theta - \sigma}. 
\end{equation}  
We note that~\eqref{eqn:THarm_evp_sinv} is the projected problem in the RR 
procedure for the matrix 
$H = (T^{\frac{1}{2}} B) (A - \sigma B)^{-1} B (T^{\frac{1}{2}} B)^{-1}$ with respect to the 
subspace spanned by the columns of $Q$. The solution of this problem 
%follows from the Galerkin condition and 
yields the Ritz pairs $(\tau, Qy)$ that tend to approximate well the 
\textit{exterior} eigenpairs of $H$; see, e.g.,~\cite{Parlett:98, Saad-book3, Stewart:01}.  

Since the matrix $H$ is similar through $T^{\frac{1}{2}} B$ to the ``shift-and-invert'' 
operator $(A - \sigma B)^{-1} B$, each of its eigenpairs is of the form 
$(1/ (\lambda - \sigma), T^{\frac{1}{2}} B v)$, where $(\lambda, v)$ is an 
eigenpair of~\eqref{eqn:eigp}. 
Thus, the desired interior eigenvalues of~\eqref{eqn:eigp} near the target $\sigma$ 
correspond to the \textit{exterior} (largest in the absolute value) eigenvalues of $H$,
which can be well approximated by the Ritz values of~\eqref{eqn:THarm_evp_sinv}. 
%or,
%equivalently, by the $T$-harmonic Rayleigh--Ritz problem~\eqref{eqn:THarm_evp}.
The associated Ritz vectors $Qy$ then represent approximations to the eigenvectors 
$T^{\frac{1}{2}} B v$ of $H$, i.e., 
$T^{\frac{1}{2}} (A - \sigma B) Z y \equiv Qy \approx T^{\frac{1}{2}} B v$.    
This implies that $Zy \approx  \alpha v$ for some scalar $\alpha$.

\subsection{The PLHR algorithm}\label{subsec:plmr_alg}

We now summarize the developments of the preceding sections and 
introduce an iterative scheme that at each step constructs a 
low-dimensional trial subspace~\eqref{eqn:plmr_span} and uses 
the $T$-harmonic RR procedure
to extract an approximate eigenvector from this subspace. The extraction 
is based on solving a 4-by-4 (3-by-3, at the initial step) eigenvalue 
problem~\eqref{eqn:THarm_evp}, which defines an appropriate $T$-harmonic Ritz pair.
The $T$-harmonic Ritz vector is then used as a new eigenvector approximation.
We call this scheme the Preconditioned Locally Harmonic Residual (PLHR) algorithm; 
see Algorithm~\ref{alg:plmr} below. 

\begin{algorithm}[htbp]
\begin{small}
\begin{center}
  \begin{minipage}{5in}
\begin{tabular}{p{0.5in}p{4.5in}}
{\bf Input}:  &  \begin{minipage}[t]{4.0in}
The matrices $A = A^*$ and $B = B^* > 0$, a preconditioner $T = T^*>0$, the shift $\sigma$, and
                 the initial guess $v^{(0)}$ for the eigenvector;
                  \end{minipage} \\
{\bf Output}:  &  \begin{minipage}[t]{4.0in}
                 An eigenvalue $\lambda$ closest to $\sigma$
                 and the associated eigenvector $v$, $\|v\|_B = 1$;
                  \end{minipage}
\end{tabular}
\begin{algorithmic}[1]
\STATE $v \gets v^{(0)}$; $v \gets v / \|v\|_B$; 
       $\lambda \gets (v, A v)$; $p \gets [ \ ]$;
\WHILE {convergence not reached}
  \STATE Compute the preconditioned residual $w \gets T ( A v - \lambda B v )$;
  \STATE Compute $s \gets T (A w - \lambda B w)$;
%  \STATE Set $Z \gets \left[ v, w, s, p \right]$;
%  \STATE $B$-orthonormalize columns of $Z$;
  \STATE Find the eigenpair $(\xi, y)$ of~\eqref{eqn:THarm_evp},
         where $Z = \left[ v, w, s, p \right]$, such that $|\xi|$ is the smallest; 
         $y = (\alpha, \beta, \gamma, \delta)^T$  (at the initial step, 
         $y = (\alpha, \beta, \gamma)^T$ and $\delta = 0$); 
  \STATE $p \gets \beta w + \gamma s + \delta p$;
  \STATE $v \gets \alpha v + p$;
%  \STATE $p \gets p/\|v\|_B$, $v \gets v/\|v\|_B$; $\lambda \gets (v, A v)$; 
  \STATE $v \gets v/\|v\|_B$; $\lambda \gets (v, A v)$; 
\ENDWHILE
\STATE Return $(\lambda, v)$.
\end{algorithmic}
\end{minipage}
\end{center}
\end{small}
  \caption{Preconditioned Locally Harmonic Residual (PLHR) algorithm}
  \label{alg:plmr}
\end{algorithm}

%%%%%%%%%%REFEREE
The initial guess $v^{(0)}$ in Algorithm~\ref{alg:plmr} can be chosen as a random vector.
In practical applications, however, certain information about the solution is available, and 
$v^{(0)}$ can give a reasonable approximation to the eigenvector of interest. 
Utilizing such initial guesses typically leads to a substantial decrease in the iteration count
and is therefore advisable.     

Note that each PLHR iteration constructs a residual of the eigenvalue problem that can be
used to determine the convergence in step 2 of Algorithm~\ref{alg:plmr}. In some applications, 
appropriate stopping criteria rely on tracking the difference between eigenvalue approximations, evaluated in step 8,
at the subsequent iterations.  
%%%%%%%%%%

%%%%%%%%%%%%%

An important property of Algorithm~\ref{alg:plmr} is that,
in contrast to the original problem~\eqref{eqn:eigp},
the $T$-harmonic problem~\eqref{eqn:THarm_evp} is no longer Hermitian. 
As such,~\eqref{eqn:THarm_evp} can have complex eigenpairs, i.e.,
the iteration parameters $y = (\alpha, \beta, \gamma, \delta)^T$ can be complex.
Thus, the PLHR algorithm should be implemented using complex arithmetic, even if
$A$ and $B$ are real. This feature is undesirable, e.g., because of the need to double the required storage
to accommodate complex numbers.
A real arithmetic version of the method for the real case is presented in Section~\ref{subsec:plmr_real}.
%, where a real arithmetic version 
%of the method is presented.     
%or, in the real case, symmetric. 
%{\color{blue} FINISH, refer to sec. 4.1!!}
%As a result, even though the eigenpairs
%of~\eqref{eqn:eigp} are real, the $T$-harmonic Ritz values $\theta = \xi + \sigma$ 
%%can be complex. 

In Algorithm~\ref{alg:plmr}, we follow the common practice of discarding the harmonic Ritz values
and replacing them with the RQs associated with the newly computed Ritz vectors; 
see, e.g.,~\cite{Morgan:91, Morgan.Zeng:98} for the standard harmonic case. 
In particular, we skip the $\theta (\equiv \xi + \sigma)$ values 
(determined in step 5 of Algorithm~\ref{alg:plmr}) and replace them with $\lambda(v) = (v^*Av)/(v^*Bv)$,
where $v = Z y$ is the $T$-harmonic Ritz vector.  
%The latter gives a better eigenvalue approximations.

Note that, at each iteration $i$, Algorithm~\ref{alg:plmr} (step 6) constructs a conjugate 
direction that can be expressed, in the superscripted notation, as 
$$p^{(i+1)} = v^{(i+1)} - \alpha^{(i)} v^{(i)} 
\equiv \beta^{(i)} w^{(i)} + \gamma^{(i)} s^{(i)} + \delta^{(i)} p^{(i)},$$
%Then, at every step, 
and the method extracts an approximate eigenvector from 
the subspace spanned by $v^{(i)}$, $w^{(i)}$, $s^{(i)}$, and $p^{(i)}$.
In exact arithmetic, this subspace is the same as $\mathcal{Z}^{(i)}$ defined 
in~\eqref{eqn:plmr_span}. 
The above formula is known to yield an improved numerical stability in 
calculating the trial subspaces 
%and is similar to the one used  
%for defining the conjugate direction in step 6 of the algorithm 
%%is known to be stable and 
%is similar to the one used 
%for constructing conjugate directions 
in the LOBPCG algorithm~\cite{Knyazev:01}. We expect that a similar property
will hold for PLHR, and therefore follow the same style for computing the conjugate
directions in Algorithm~\ref{alg:plmr}.  

If it is possible to store additional vectors that contain results
of the matrix-vector multiplications with $A$, $B$, and $T(A - \sigma B)$, 
then Algorithm~\ref{alg:plmr} can be implemented with two matrix-vector multiplications
involving $A$ and $B$, and four applications of the preconditioner $T$.
Two of the preconditioning operations result from the construction of the trial subspace 
%construction
in steps 3 and 4, whereas the other two are the consequence of the $T$-harmonic
extraction.
%~\eqref{eqn:THarm_evp}. 
%which requires evaluating $T(A - \sigma B) Z$. 

%In practical implementations, it is desireable 
%to $B$-orthonormalize the basis of the trial subspace $\mathcal{Z}^{(i)}$ at each PLHR iteration.
%The orthogonalization step typically improves the numerical stability
%of the scheme and allows approximating the solution to a higher accuracy. Thus, the matrix $Z$
%in step 5 of Algorithm~\ref{alg:plmr} should be replaced by a matrix 
%$\hat Z = [\hat v, \hat w, \hat s, \hat p]$, whose columns represent a $B$-orthonormal
%basis of the trial subspace defined by $Z$. Then the new approximate eigenvector 
%is computed as a linear combination of $\hat Z$'s columns, similar to steps 6 and 7 of Algorithm~\ref{alg:plmr}:
%$p \gets \beta \hat w + \gamma \hat s + \delta \hat p$ and $v \gets \alpha \hat v + p$.
% 

%%%%%%%%%%%%%%%%%%%%%%%%% BEGIN PLMR

\subsection{Possible variations}\label{subsec:plmr_thesis}

The choice of the trial subspaces and of the extraction procedure,
motivated by the preconditioned null space finding framework in Section~\ref{sec:idealized}, 
is not unique. For example, one can define the ``s-vector'' as $s^{(i)} = T(A w^{(i)} - \theta^{(i)} B w^{(i)})$,
where instead of the RQ $\lambda^{(i)}$ the current $T$-harmonic value $\theta^{(i)}$ is used.
Alternatively, it is possible to set $s^{(i)} = T(A w^{(i)} - \sigma B w^{(i)})$, i.e., use
the fixed $\sigma$ instead of $\lambda^{(i)}$. However, our numerical experience suggests that these options do
not lead to any better results compared to the original choice of $s^{(i)}$ in~\eqref{eqn:plmr_span}.
  
The minimal residual condition~\eqref{eqn:coeff} can also motivate an extraction procedure that is
different from the $T$-harmonic RR described in Section~\ref{subsec:THarm}. In particular,
let us assume that $\tilde \lambda \in \mathbb{R}$ is some approximation to the targeted 
eigenvalue $\lambda_q$ and $v^{(i)}$ is the current approximate eigenvector.
%, i.e., $\tilde \lambda \approx \lambda_q$. 
In this case, one can extend~\eqref{eqn:coeff} to minimize the $T$-norm of the residual-like vector 
$\tilde r  = A v - \tilde \lambda B v$, such that 
%attempt to extract the corresponding eigenvector approximation $v \in \mathcal{V}$ 
%by satisfying the following optimality condition:
\begin{equation}\label{eqn:rRR}
%v = \mbox{arg}\min_{ z \in \mathcal{V}, \| z \|_{B} = 1} \|A z - \tilde \lambda B z \|, 
v^{(i+1)} = \underset{v \in \mathcal{Z}^{(i)}, \| v \|_{B} = 1} {\operatorname{argmin}} \|A v - \tilde \lambda B v \|_T.
\end{equation}  
%where $\|z\|^2_B = (z,B z)$, and 
%$\| z \| = \| z \|_I$ is the $2$-norm. 
This minimization principle represents
%in (\ref{eqn:rRR}),
%in fact, 
%defines the \textit{refined procedure}, 
%also called 
the \textit{refinement procedure}~\cite{Jia:97}, 
%where the residual refinement is 
performed with respect to the $T$-norm.
%which is
%straightforwardly extended to the case of the generalized eigenproblem
%(the original condition in \cite{Jia:97} was formulated for $B = I$). 
%The minimizer $v$ in (\ref{eqn:rRR}) is called the \textit{refined approximate eigenvector}.
%
%Given an SPD preconditioner $T$, we modify 
%condition (\ref{eqn:rRR}) to perform the minimization in the preconditioner-based 
%$T$-norm, rather than in the standard $2$-norm, i.e., 
%\begin{equation}\label{eqn:rRRT}
%v = \underset { z \in \mathcal{V}, \| z \|_{B} = 1} {\operatorname{argmin}} \|A z - \tilde \lambda B z \|_T, 
%\end{equation}  
%where $\|z\|^2_T = (z,T z)$.  

Assuming that $Z$ is a basis of the current trial subspace $\mathcal{Z}^{(i)}$,
one can show that~\eqref{eqn:rRR} yields the new approximate eigenvector $v^{(i+1)} = Z y_{\min}$,
where $y_{\min}$ minimizes the bilinear form 
\[
\theta^2(y) =  \frac{( y, Z^* (A  - \tilde \lambda B )T(A  - \tilde \lambda B) Z y)}{(y, Z^* B Z y)}     
\]
over all vectors $y$ of dimension $4$.
%the matrix $V \in \mathbb{R}^{n \times k}$ is such that $\mbox{col}(V) = \mathcal{V}$,
%where $\mbox{col}(V)$ denotes the column space of $V$, and, hence, 
%any $z \in \mathcal{V}$ is of the form $z = V y$, for some $y \in \mathbb{R}^k$, we get 
%\begin{eqnarray}
%\nonumber \|A z - \tilde \lambda B z \|^2_T 
%& = & (A z - \tilde \lambda B z, T(A z - \tilde \lambda B z)) = (z, (A  - \tilde \lambda B )T(A  - \tilde \lambda B) z) \\
%\nonumber
%& = & (V y, (A  - \tilde \lambda B )T(A  - \tilde \lambda B) V y) \\
%& = &  ( y, V^* (A  - \tilde \lambda B )T(A  - \tilde \lambda B) V y).
%\end{eqnarray} 
%Thus, (\ref{eqn:rRRT}) can be replaced by the problem of finding the minimizer $y_{min} \in \mathbb{R}^k$, such that
%\begin{eqnarray}
%\nonumber
%y_{min} & = & \underset{ y \in \mathbb{R}^k, \| Vy \|_{B} = 1}{\operatorname{argmin}} ( y, V^* (A  - \tilde \lambda B )T(A  - \tilde \lambda B) V y) \\
%\nonumber
%  & = & \underset{ y \in \mathbb{R}^k}{\operatorname{argmin}} \frac{( y, V^* (A  - \tilde \lambda B )T(A  - \tilde \lambda B) V y)}{(V y,B V y)} \\   
%\nonumber
%  & = & \underset{ y \in \mathbb{R}^k}{\operatorname{argmin}} \frac{( y, V^* (A  - \tilde \lambda B )T(A  - \tilde \lambda B) V y)}{(y, V^* B V y)},      
%\end{eqnarray}  
This is equivalent to solving 
%the 
%problem of finding eigenvector of 
the eigenvalue problem
%corresponding to the smallest eigenvalue $\theta^2_{min}$, of the
%$k$-by-$k$ generalized symmetric eigenvalue problem
%\begin{equation}\label{eqn:rRRevp} 
$
Z^* (A  - \tilde \lambda B )T(A  - \tilde \lambda B) Z y = \theta^2 Z^* B Z y,
$
%\end{equation}
where the desired $y_{\min}$ is the eigenvector associated with the smallest eigenvalue.
%The square root of the smallest eigenvalue in (\ref{eqn:rRRevp}), i.e.,
%$\theta_{min}$, gives the minimal value of norm (\ref{eqn:rRRT}),
%while the eigenvector $y_{min}$ determines the corresponding minimizer  
%\begin{equation}\label{eqn:v_refined}
%v = V y_{min}, \ \|v\|_B = 1,
%\end{equation} 
%which we set as the new eigenvector approximation.
%The value of $\theta_{min}$ is typically discarded.
%where $\theta_{\min}$ is the minimal value of the norm in (\ref{eqn:rRRT}).
%% over theunit (in the $B$-norm) sphere in $\mathcal{V}$. 
%We further discard $\theta_{\min}$,
%however, use $y_{\min}$ to define the new eigenvector approximation, i.e.,

The described approach can be of interest, in particular, if $v^{(i)}$ is already 
reasonably close to the targeted eigenvector. In this case, $\tilde \lambda$ 
in ~\eqref{eqn:rRR} can be set to the current RQ $\lambda^{(i)}$.
%, so that~\eqref{eqn:rRR} represents the residual norm minimization. 
After the minimizer $v^{(i+1)}$ is constructed and the corresponding RQ and the new trial subspace
$\mathcal{Z}^{(i+1)}$ are computed, the procedure is repeated once again. Such an approach,
called the Preconditioned Locally Minimal Residual (PLMR) method,
which combines the trial subspaces~\eqref{eqn:plmr_span} with the eigenvector extraction in~\eqref{eqn:rRR},
has been described in~\cite{thesis}.

\section{The block PLHR algorithm}\label{sec:bplhr}
%It is common in practical applications that \textit{several} eigenpairs are of interest. 
%In this case, a possible approach is to apply Algorithm~\ref{alg:plmr} for computing the desired solutions one after
%another by restricting iterations to the orthogonal complement of the already converged eigenvectors. 
%%This technique is called \textit{deflation by restriction} [{\color{blue} ???}].
%This gives a deflated variant of PLHR.
%
%Another standard option is to use a \textit{block iteration} which 
%attempts to approximate all the targeted eigenpairs simultaneously. 
%In modern extreme-scale eigenvalue computations, block methods are typically preferred over the 
%deflation based techniques due to their ability to reveal additional levels parallelism and to improve performance
%through an intensive use of BLAS3 operations. Block iterations also represent a traditional tool for  
%hanling clusters of eigenvalues. 
%%
In this Section, we extend PLHR to the block case, where several eigenpairs
closest to $\sigma$ are computed simultaneously. 
%We call this version the block PLHR (BPLHR).
%The  
%called the Block PLHR (BPLHR),
%and discuss a few implementational details of the algorithm.

We start by introducing the block notation. 
Let $\Lambda = \text{diag}\{ \lambda_1, \lambda_2, \ \dots, \ \lambda_k \}$ denote a $k$-by-$k$ 
diagonal matrix of the targeted eigenvalues, ordered according to their distances from the shift $\sigma$,
so that $| \lambda_l - \sigma | \leq | \lambda_j - \sigma|$ for $l < j$, 
where $j = 1, 2, \ldots, k$.   
Let $V$ be an $n$-by-$k$ matrix of the associated eigenvectors. 
We assume that 
\[
V^{(i)} = [v^{(i)}_1, v^{(i)}_2, \ \dots, \ v^{(i)}_k], \quad
\Lambda^{(i)} = \text{diag}\{ \lambda^{(i)}_1, \lambda^{(i)}_2, \ \dots, \ \lambda^{(i)}_k \} 
\]
are the matrices of the approximate eigenvectors and eigenvalues at iteration $i$,
respectively. The diagonal entries of $\Lambda^{(i)}$ are the RQs 
%and
%of approximate eigenvectors~$v_j^{(i)}$, such that $\|v_j^{(i)}\|_B = 1$, and a matrix  
%\[
%\Lambda^{(i)} = \text{diag}\{ \lambda^{(i)}_1, \lambda^{(i)}_2, \ \dots, \ \lambda^{(i)}_k \} 
%%\quad \lambda^{(i)}_j = (v_j^{(i)},A v_j^{(i)})
%\]
%is a diagonal matrix of the corresponding RQs 
\[
\lambda^{(i)}_j \equiv \lambda(v_j^{(i)}) = (v_j^{(i)},A v_j^{(i)})/(v_j^{(i)},B v_j^{(i)}),
%, \quad j = 1, 2, \ldots, k.
\]
evaluated at the corresponding vectors $v_j^{(i)}$, such that 
%ordered according to their distance to the targeted shift $\sigma$,
%so that 
$| \lambda^{(i)}_l - \sigma | \leq | \lambda^{(i)}_j - \sigma|$ for $l < j$. 
%where $j = 1, 2, \ldots, k$.   
%We assume that each $v_j^{(i)}$ has a unit $B$-norm. Furthermore, 
%Note that we do not require $B$-orthogonality of the approximate eigenvectors and only
%assume their linear independence.  
%Associated with $V^{(i)}$ is the diagonal matrix of the RQs 
%which contains approximations of the targeted eigenvalues.  

Given the approximate eigenpairs in $V^{(i)}$ and $\Lambda^{(i)}$, let us define the block 
\[
W^{(i)} \equiv [ w^{(i)}_1, w^{(i)}_2, \ \dots, \ w^{(i)}_k]  = T(AV^{(i)} - BV^{(i)}\Lambda^{(i)}) 
\] 
of the preconditioned residuals $w^{(i)}_j = T(A v_j^{(i)} - \lambda_j^{(i)} B v_j^{(i)})$, and
the block  
%is a preconditioned residual
%for the pair $(\lambda_j^{(i)}, v_j^{(i)})$. 
%for the approximate eigenvector $v^{(i)}_j$. 
%By analogy with Algorithm~\ref{alg:plmr}, 
%We also introduce a block of ``$s$-vectors'', 
\[
S^{(i)} \equiv [ s^{(i)}_1, s^{(i)}_2, \ \dots, \ s^{(i)}_k]  = T(AW^{(i)} - BW^{(i)}\Lambda^{(i)}), 
\] 
of ``s-vectors'' $s^{(i)}_j = T(A w_j^{(i)} - \lambda_j^{(i)} B w_j^{(i)})$.
%Analogously to the single-vector PLHR iteration, 
We can then introduce a trial subspace~$\mathcal{Z}^{(i)}$ spanned by the columns of $V^{(i)}$,
$W^{(i)}$, $S^{(i)}$, and $V^{(i-1)}$ ($V^{(-1)} = \textbf{0}$), i.e., 
\begin{equation}\label{eqn:plmr_span_block}
\mathcal{Z}^{(i)} = 
\mbox{span}\left\{ v_1^{(i)},\ldots, v_k^{(i)}, 
w_1^{(i)} \ldots w_k^{(i)}, 
s_1^{(i)}, \ldots, s_k^{(i)}, 
v_1^{(i-1)}, \ldots v_k^{(i-1)} \right\}.
\end{equation}
Clearly,~\eqref{eqn:plmr_span_block} is a block generalization of the trial subspace~\eqref{eqn:plmr_span}
constructed at each single-vector PLHR iteration.

Let $Z$ be a matrix whose columns represent a basis of~\eqref{eqn:plmr_span_block}. 
Analogously to Algorithm~\ref{alg:plmr}, we require that
new eigenvector approximations $V^{(i+1)}$ are the $T$-harmonic Ritz vectors 
extracted from~\eqref{eqn:plmr_span_block}. 
Thus, the vectors $V^{(i+1)}$ are defined by solving the $4k$-by-$4k$ (or, at the initial step, $3k$-by-$3k$) 
eigenvalue problem~\eqref{eqn:THarm_evp}, and correspond to the $k$ $T$-harmonic Ritz values that are closest to the shift $\sigma$.
The entire approach, referred to as the \textit{block PLHR (BPLHR)}, is summarized in Algorithm~\ref{alg:bplmr}. 

\begin{algorithm}[htbp]
\begin{small}
\begin{center}
  \begin{minipage}{5in}
\begin{tabular}{p{0.5in}p{4.5in}}
{\bf Input}:  &  \begin{minipage}[t]{4.0in}
The matrices $A$ and $B$, a preconditioner $T = T^*>0$, the shift $\sigma$, and
                 the initial guess $V^{(0)}$ for $k$ eigenvectors;
                  \end{minipage} \\
{\bf Output}:  &  \begin{minipage}[t]{4.0in}
                 Diagonal matrix $\Lambda$ of eigenvalues closest to the target $\sigma$
                 and the matrix $V$ of the associated eigenvectors;
                  \end{minipage}
\end{tabular}
\begin{algorithmic}[1]

\STATE $V \gets V^{(0)}$; $P \gets [ \ ]$; 
\STATE Normalize columns of $V$ to have a unit $B$-norm; $\Lambda \gets \mbox{diag}\left( V^* A V \right)$;  
%\STATE Compute and store $AV$, $BV$, and $T(A - \sigma B) V$;
%\STATE $P \gets [ \ ]$, $AP \gets [ \ ]$, $BP \gets [ \ ]$, and $T(A - \sigma B) P \gets [ \ ]$; 
%\STATE $P \gets [ \ ]$; 
%, and $T(AP - \sigma BP) \gets [ \ ]$;
\WHILE {convergence not reached}
  \STATE Compute the preconditioned residuals $W \gets T(AV - B V \Lambda)$;
%  \STATE Compute and store $AW$, $BW$, and $T(A - \sigma B)W$;  
  \STATE Compute $S \gets  T (A W - B W \Lambda)$;
%  \STATE Compute and store $AS$, $BS$, and $T(A - \sigma B)S$;  
%  \STATE Compute $[V, R] \gets \texttt{QR\_B}(V)$\footnote{$[Q,R] \gets \texttt{QR\_B}(Y)$ denotes a truncated QR decomposition of 
%the matrix $Y$, such that the columns of $Q$ represent a $B$-orthonormal basis of $Y$'s column space and $R$ is upper triangular.   
%};
%  \STATE $AV \gets ( AV ) R^{-1}$, $BV \gets ( BV ) R^{-1}$,~and~$T(A - \sigma B)V \gets \left( T(A - \sigma B) V \right) R^{-1}$;  
%  \STATE $B$-orthogonalize $W$ against $V$. Update $AW$, $BW$, and $T(A - \sigma B)W$ accordingly. 
%%  \STATE Update $AW \gets AW - AV (V^* BW)$, $BW \gets BW - BV (V^* BW)$,~and~$T(A - \sigma B)W \gets T(A - \sigma B) W - T(A - \sigma B)V (V^* BW)$;  
%%   \STATE $B$-orthogonalize $W$ against $V$: $W \gets W - V (V^* B W)$;     
%  \STATE Compute $[W, R] \gets \texttt{QR\_B}(W)$;
%  \STATE $AW \gets ( AV ) R^{-1}$, $BW \gets ( BV ) R^{-1}$,~and~$T(A - \sigma B)W \gets \left( T(A - \sigma B) W \right) R^{-1}$;  
%  \STATE Orthogonalize $S$ against $V$ and $W$; Update $AW$, $BW$, and $T(A - \sigma B)W$ accordingly.
%  \STATE $B$-orthogonalize columns of $S$;
%  \STATE $B$-orthogonalize $P$ against $V$, $W$, and $S$. Update 
%$AP$, $BP$, and $T(A - \sigma B) P$ accordingly. 
  \STATE Set $Z \gets \left[ V,\ W, \ S, \ P \right]$. 
$B$-orthonormalize the columns of $Z$. 
Let $\hat Z = [ \hat V,\ \hat W, \ \hat S, \ \hat P ]$ be the matrix of the resulting
$B$-orthonormal columns. 
%$BZ \gets \left[ BV,\ BW, \ BS, \ BP \right]$, 
%and $T(A - \sigma B)Z \gets \left[ T(A - \sigma B)V,  \ T(A - \sigma B)W, \ T(A - \sigma B)S, \ T(A - \sigma B)P \right]$. 
  \STATE Find eigenpairs of the projected problem~\eqref{eqn:THarm_evp} with $Z \equiv \hat Z$. 
Define $Y \equiv [Y_V^T, \ Y_W^T, \ Y_S^T, \ Y_P^T]^T$
to be the matrix of $k$ eigenvectors of \eqref{eqn:THarm_evp} corresponding to the smallest in the absolute value eigenvalues. 
%$\Omega$ is the diagonal matrix
%         of eigenvalues and $Y \equiv [Y_V^T, \ Y_W^T, \ Y_S^T, \ Y_P^T]^T$ contains the 
%         corresponding eigenvectors of a unit $B$-norm (at the initial step $Y_P = [ \ ]$);
  \STATE  Compute $P \gets \hat W Y_W + \hat S Y_S + \hat P Y_P$;
%  \STATE  Update $AP \gets (AW) Y_W + (AS) Y_S + (AP) Y_P$; 
%  \STATE  Update $BP \gets (BW) Y_W + (BS) Y_S + (BP) Y_P$; 
%  \STATE  Update $T(A-\sigma B)P \gets \left(T(A - \sigma B)P \right) Y_W + \left(T(A - \sigma B)P\right) Y_S + \left( T(A - \sigma B)P \right) Y_P$;     

  \STATE  Compute new approximate eigenvectors $V \gets \hat V Y_V + P$; 

%  \STATE  Update $AV \gets (AV) Y_V + AP$; 
%  \STATE  Update $BV \gets (BV) Y_V + BP$; 
%  \STATE  Update $T(A-\sigma B)P \gets \left(T(A - \sigma B)V \right) Y_V + T(A - \sigma B)P$;
  \STATE  Normalize columns of $V$ to have a unit $B$-norm; $\Lambda \gets \mbox{diag}\left( V^* A V \right)$;
\ENDWHILE
\STATE Perform the standard RR procedure for~\eqref{eqn:eigp} with respect to $V$. 
Update $V$ to contain the Ritz vectors and $\Lambda$ the corresponding Ritz values.  
\STATE Return $(\Lambda, V)$.
\end{algorithmic}
\end{minipage}
\end{center}
\end{small}
  \caption{The block PLHR (BPLHR) Algorithm}
  \label{alg:bplmr}
\end{algorithm}

Note that the set of approximate eigenvectors constructed by Algorithm~\ref{alg:bplmr} is 
generally not $B$-orthogonal. This is a consequence of the $T$-harmonic projection that is 
based on solving the non-Hermitian reduced problem~\eqref{eqn:THarm_evp}. The eigenvectors $Y$
of this problem are non-orthogonal. Therefore the corresponding $T$-harmonic Ritz vectors $V = Z Y$
are also non-orthogonal.

%Theoretically, since~\eqref{eqn:THarm_evp} is a non-Hermitian problem, it is possible
%that deficient eigenvalues can occur, i.e., the basis of eigenvectors does not exist.
%While attention should be paid to this case in practical BPLHR implementations, such
%situations are extremely rare. In fact, in our tests, we have never observed deficient eigenvalues.

As iterations proceed and the approximate eigenpairs get closer to the solution,
the (near) $B$-orthogonality of the columns of $V$ starts to show up naturally,
due to the intrinsic property of the symmetric eigenvalue problem. However, since in
many applications the required accuracy of the solution may be low, in order to ensure 
that the returned eigenvector approximations are $B$-orthonormal, 
%of the resulting eigenvector approximations, 
the algorithm performs the post-processing (step 12 of Algorithm~\ref{alg:bplmr}),
where the final block $V$
%the final eigenvector block should be post-processed by performing the standard RR procedure after step 10
%of the algorithm, so that the $T$-harmonic Ritz vectors are 
is ``rotated'' to the $B$-orthogonal set of the Ritz vectors.

For the purpose of numerical stability, at each iteration, Algorithm~\ref{alg:bplmr} $B$-orthogona-lizes
the set of vectors $Z = [ V,\ W, \ S, \ P ]$ that span the trial subspace  
and performs the $T$-harmonic RR procedure with respect to the $B$-orthonormal basis
$\hat Z = [ \hat V,\ \hat W, \ \hat S, \ \hat P ]$. This is done by first $B$-orthonormalizing the 
columns of $V$, so that the column space of $\hat V$ and $V$ is the same. Then the block $W$ is
$B$-orthogonalized against $\hat V$ and the resulting set of vectors is $B$-orthonormalized to
obtain the block $\hat W$. The remaining blocks $\hat S$ and $\hat P$ are obtained in the same manner,
by orthogonalizing against the currently available $B$-orthogonal blocks.

The BPLHR algorithm requires storage for $16k$ vectors of $Z$, $AZ$, $BZ$, and $T(A - \sigma B) Z$.
If the storage is available, it can be implemented using $2$ matrix-block multiplications involving $A$ and $B$, 
and $4$ applications of the preconditioner $T$ to a block per iteration. 
%The products with $A$ and $B$, as well as the two block preconditioning operations, 
%result from construction of the trial subspace. The two additional
%applications of $T$ are needed to accomplish the $T$-harmonic projection,
%which requires computing $T(A - \sigma B) W$ and $T(A - \sigma B) S$
%(matrices $T(A - \sigma B) V$ and $T(A - \sigma B) P$ can be updated without additional preconditionings). 
Note that in the case of standard eigenvalue problem ($B = I$), the memory requirement reduces to 
storing $12k$ vectors, and the two matrix-block multiplications with $B$ are no longer needed.     
Additional cost reductions can result from locking the converged eigenpairs. In our BPLHR
implementation, this is done by the standard soft locking procedure, similar to~\cite{Kn.Ar.La.Ov:07}.

\subsection{The BPLHR algorithm in real arithmetic}\label{subsec:plmr_real}

Each BPLHR iteration relies on the solution of the projected
problem~\eqref{eqn:THarm_evp}, which is generally non-Hermitian. 
As a consequence,~\eqref{eqn:THarm_evp} can have complex eigenpairs and, therefore, 
complex arithmetic should be assumed for Algorithm~\ref{alg:bplmr}. 

However, in the case where $A$ and $B$ are real symmetric, the use of the complex arithmetic is unnatural,
since the solution of~\eqref{eqn:eigp} is real. Moreover, the presence of complex-valued iterations 
translates into the undesirable requirement to double the memory allocation for each vector, or block of vectors, 
utilized in the computation. The goal of the present subsection is to discuss ways to overcome this issue and
introduce a real arithmetic version of the BPLHR algorithm. 
%that handles the real symmetric case of problem~\eqref{eqn:eigp}  
%in real arithemetic. 

%As has been already mentioned, the $T$-harmonic Rayleigh-Ritz procedure,
%described in section~\ref{subsec:THarm},
%leads to the loss of symmetry in the projected problem. As a consequence,
%the reduced eigenvalue problem~\eqref{eqn:THarm_evp} can have complex eigenpairs.
%If $A$ and $B$ are complex Hermitian then this departure
%from symmetry causes no trouble. In this case Algorithm~\ref{alg:bplmr} 
%operates on complex-valued blocks to deliver the desired set of complex 
%eigenvectors and the associated eigenvalues. 
%%of~\eqref{eqn:eigp}.   

%If $A$ and $B$ (and the preconditioner $T$) are real symmetric, 
%then the presence of complex
%eigenpairs $(\xi, y)$ in~\eqref{eqn:THarm_evp} represents a more notable issue.
%In particular, a complex eigenvector $y$ yields a complex $T$-harmonic
%Ritz vector $v = Z y$, and hence the block of updated approximate eigenvectors 
%$V^{(i)}$ can have complex columns. In practice, this 
%%property of the 
%%$T$-harmonic eigenpair extraction 
%requires that Algorithm~\ref{alg:bplmr} 
%is either capable of switching to complex-valued calculations at runtime or is 
%entirely implemented in complex arithmetic. Since for real $A$ and $B$ 
%the eigenvectors of~\eqref{eqn:eigp} can also be chosen real, this feature 
%of the  method may not be desirable, e.g., due to doubling the size of processed data.   
%We next suggest a possible treatment for this problem, which allows us to maintain 
%iterations of Algorithm~\ref{alg:bplmr} in real arithmetic.     

Let $A$ and $B$ be real, where $A$ is symmetric and $B$ symmetric positive definite (SPD).
Let $T$ be an SPD preconditioner, and assume that 
%Let $A$ and $B$ be real symmetric, and assume that 
at the current iteration of Algorithm~\ref{alg:bplmr} 
the trial subspace is given by a real matrix $Z$. 
In this case the left- and right-hand side matrices in the $T$-harmonic problem~\eqref{eqn:THarm_evp} are real,
implying that the eigenvalues $\xi$ are either real or appear in the complex conjugate pairs. 
Specifically, if $(\xi, y)$ is a complex eigenpair of~\eqref{eqn:THarm_evp} 
then $(\bar \xi, \bar y)$ is also an eigenpair of~\eqref{eqn:THarm_evp}.
%, where the bar denotes the complex conjugation. 

Clearly, if $y$ is a real eigenvector of~\eqref{eqn:THarm_evp} then the corresponding $T$-harmonic
Ritz vector $v = Z y$ is also real. If $y$ is complex, then $v = Z y$ is complex, and the conjugate 
eigenvector $\bar y$ gives $\bar v = Z \bar y$, i.e., the presence of complex solutions in~\eqref{eqn:THarm_evp}
yields complex conjugate $T$-harmonic Ritz vectors $v$ and $\bar v$.

Complex solutions of~\eqref{eqn:THarm_evp} and, subsequently,  the 
complex conjugate $T$-harmonic Ritz vectors, are not rare in practice. Their presence
indicates that the extraction procedure is attempting to approximate eigenpairs associated with 
a (real) eigenvalue 
%of~\eqref{eqn:eigp} 
of multiplicity greater than one.  
In particular, if $v = v_{\text{R}} + i v_{\text{I}} \in \mathbb{C}^{n}$ and 
$\bar v = v_{\text{R}} - i v_{\text{I}} \in \mathbb{C}^{n}$ 
is a pair of the complex conjugate $T$-harmonic Ritz vectors, where $v_{\text{R}} \in \mathbb{R}^{n}$ 
and $v_{\text{I}} \in \mathbb{R}^{n}$
are the real and imaginary parts of $v$, then the eigenvalue approximations
given by the corresponding RQs coincide, i.e., 
\begin{equation}\label{eqn:lambda_complex}
\lambda(v) = \lambda(\bar v) = \frac{v_{\text{R}}^* A v_{\text{R}} + v_{\text{I}}^* A v_{\text{I}} }
{v_{\text{R}}^* B  v_{\text{R}} + v_{\text{I}}^* B v_{\text{I}}}. % \in \mathbb{R}.
\end{equation}
Thus, $v$ and $\bar v$ indeed approximate eigenvectors corresponding to the same eigenvalue. 
%For example, the $T$-harmonic Ritz values $\theta = \xi + \sigma$ and  
%$\bar \theta = \bar \xi + \sigma$ can be viewed as two equidistant
%approximations to such a multiple eigenvalue from the upper and lower 
%halves of the complex plain. 
%This can be 
%Additionally, note that the RQs 
%evaluated at the corresponding $T$-harmonic Ritz vectors $v = Z y$ and $\bar v = Z \bar y$ 
%are equal, yielding the same eigenvalue approximation for both $v$ and $\bar v$.

%In the former case, the eigenvectors
%of~\eqref{eqn:THarm_evp}, and hence the corresponding $T$-harmonic
%Ritz vectors, can also be chosen real. Therefore, maintaining
%real arithmetic for this common, according to our numerical experiments, 
%situation is straightforward. 
%Now let us consider the other possibility.
%Suppose that $(\xi, y)$ is an eigenpair of~\eqref{eqn:THarm_evp} corresponding to a 
%complex eigenvalue $\xi = \xi_{\text{R}} + i \xi_{\text{I}}$,  
%where $\xi_{\text{R}}$ and $\xi_{\text{I}} > 0$ are the real and imaginary
%parts of $\xi$, respectively. The eigenvector $y = y_{\text{R}} + i y_{\text{I}}$ 
%is complex, with real and imaginary parts given by the real vectors 
%$y_{\text{R}}$ and $y_{\text{I}}$. 
%%
%Then there exists an eigenpair $(\bar \xi, \bar y)$
%associated with the conjugate eigenvalue $\bar \xi = \xi_{\text{R}} - i \xi_{\text{I}}$,
%where the eigenvector $\bar y = y_{\text{R}} - i y_{\text{I}}$ is obtained from
%$y$ by the componentwise conjugation. 

Let $Y$ denote the matrix of eigenvectors of~\eqref{eqn:THarm_evp}
associated with $k$ smallest magnitude eigenvalues. For simplicity, we assume that if $y$ is a 
complex eigenvector in~$Y$, then the eigenvector $\bar y$ is also included into $Y$, i.e.,
$Y$ contains complex conjugate columns $y$ and $\bar y$ (the case where $Y$ contains only the column $y$ will be discussed below). 
Then, for a given $Y$, let us define a real matrix $Y^{\prime} = [Y_0 \ Y_{\text{R}} \ Y_{\text{I}}]$, 
where the subblock $Y_0$ consists of the real columns of $Y$, 
whereas $Y_{\text{R}}$ and $Y_{\text{I}}$ contain the real and 
imaginary parts of the complex columns of $Y$. More precisely, corresponding to each complex 
conjugate pair of columns $y = y_{\text{R}} + i y_{\text{I}}$ and 
$\bar y = y_{\text{R}} - i y_{\text{I}}$ in $Y$,
%, where $y_{\text{R}}$ and $y_{\text{I}}$ are the real and imaginary parts of $y$, 
we define two real columns $y_{\text{R}}$ and $y_{\text{I}}$, where the former is placed
to $Y_{\text{R}}$ and the latter to~$Y_{\text{I}}$. It is assumed that all columns 
$y_{\text{R}}$ and $y_{\text{I}}$ appear in $Y_{\text{R}}$ and $Y_{\text{I}}$
in the same order. That is, if $y_{\text{R}}$ is the $j$th column in the subblock $Y_{\text{R}}$,
then $y_{\text{I}}$ is the $j$th column of $Y_{\text{I}}$. Clearly, $Y_{\text{R}}$ and $Y_{\text{I}}$ 
contain the same number of columns, which is equal to the number of complex eigenvectors in $Y$ divided by two.
Thus, the sizes of $Y$ and $Y^{\prime}$ are identical.               

Let $V = Z Y$ be the block of approximate eigenvectors extracted by Algorithm~\ref{alg:bplmr} from
the trial subspace defined by $Z$ at the current BPLHR step, and let $[ V,\ W, \ S, \ P ]$ be the basis of 
the trial subspace computed from $V$ at the next iteration. It is clear that each block in this basis
will have complex columns, provided that $Y$ contains complex eigenvectors of~\eqref{eqn:THarm_evp}. 
Additionally, according to the assumption on $Y$, for each complex column $v$ of the block $V$ ,
there exists a conjugate column $\bar v$ in $V$. 
Below we show that it is possible to replace $[ V,\ W, \ S, \ P ]$ by a real basis 
$[  V^{\prime},\  W^{\prime}, \ S^{\prime}, \  P^{\prime} ]$ that spans exactly the same trial subspace.
This idea will lead to a real arithmetic version of the BPLHR algorithm.

%In order to introduce real arithmetic into the BPLHR algorithm, we show that it is always possible
%to define a real basis $Z^{\prime}$ that spans exactly the same trial subspace as the one given by 
%the matrix $Z$ constructed at each step of Algorithm~\ref{alg:bplmr}. As a consequence, all operations 
%can be performed in real arithmetic and no complex data should be stored.

Our  approach is based on the observation that
$\text{span}\{ v_{\text{R}}, v_{\text{I}} \} = \text{span}\{v, \bar v\}$ over the field 
of complex numbers. Therefore, if we define  
\begin{equation}\label{eqb:vcmplx}
V^{\prime} \equiv [V_0 \ V_{\text{R}} \ V_{\text{I}}] =  [Z Y_0 \  Z Y_{\text{R}} \ Z Y_{\text{I}}],
\end{equation}
then its column space will be exactly the same as that of $V$. 
Here the subblock $V_0$ of $V^{\prime}$ contains the real columns of the matrix $V$, 
whereas $V_{\text{R}}$ and $V_{\text{I}}$ correspond to the real and imaginary 
parts of $V$'s complex columns, respectively. 
Thus, one can replace the block $V$ of $T$-harmonic Ritz vectors 
by the corresponding real-valued block $V^{\prime}$ without changing the column span. 

%%%%%%%%%%%REFEREE
Note that both $V$ and $V^{\prime}$ have $k$ columns. However, all columns of $V^{\prime}$
are real. Therefore, placing them into the new trial subspace does not lead to any increase in storage, 
in contrast to using $V$, whose complex columns would require extra memory.   
%%%%%%%%%%%

Similarly, since the presence of the complex conjugate pair of columns $v$ and $\bar v$
in $V$ yields the conjugate columns $w$ and $\bar w$ in the block $W$, the latter  
can be replaced by the real block 
\begin{equation}\label{eqn:wcmplx}
W^{\prime} \equiv [W_0 \ W_{\text{R}} \ W_{\text{I}}] = [T(A V_0 - B V_0 \Lambda_0) \ T(A V_{\text{R}} - B V_{\text{R}}\Lambda_1) \ 
T(A V_{\text{I}} - B V_{\text{I}}\Lambda_1)],
\end{equation}
where $\Lambda_0$ is the diagonal matrix of the RQs $\lambda(v) = v^* A v / v^* B v$ evaluated at the 
columns $v$ of $V_0$; and $\Lambda_1$ is diagonal matrix that, for each column $v_{\text{R}}$ of $V_{\text{R}}$ and 
the corresponding $v_{\text{I}}$ of $V_{\text{I}}$, contains the RQ defined in~\eqref{eqn:lambda_complex}.
Note that the matrix $W^{\prime}$ can be expressed as $W^{\prime} = T(AV^{\prime} - BV^{\prime}  \Lambda^{\prime})$, 
where $\Lambda^{\prime} = \text{diag}(\Lambda_0, \Lambda_1, \Lambda_1)$. 

%This is a consequence of the fact that $v$ and $\bar v$ yield a complex conjugate pair 
%of preconditioned residuals $w$ and $\bar w$.

%In~\eqref{eqn:wcmplx}, $\Lambda_0$ is the diagonal matrix of the RQs $\lambda(v) = v^* A v / v^* B v$ evaluated at the 
%columns of $V_0$. The matrix $\Lambda_1$ is diagonal. 
%For each column $v_{\text{R}}$ of $V_{\text{R}}$ and 
%$v_{\text{I}}$ of $V_{\text{I}}$, such that $v = v_{\text{R}} + i v_{\text{I}}$ is the $T$-harmonic
%Ritz vector, it contains the RQ defined in~\eqref{eqn:lambda_complex}.

Finally, using the same argument, the block $S$ constructed by Algorithm~\ref{alg:bplmr} can be replaced by the real block
\begin{equation}\label{eqn:scmplx}
S^{\prime} \equiv [S_0 \ S_{\text{R}} \ S_{\text{I}}] = [T(A W_0 - B W_0 \Lambda_0) \ T(A W_{\text{R}} - B W_{\text{R}}\Lambda_1) \ 
T(A W_{\text{I}} - B W_{\text{I}}\Lambda_1)],
\end{equation}
%Clearly, given a real $Y^{\prime}$, the block $P^{\prime}$ of the new conjugate directions
%can also be chosen real.
%
i.e., $S^{\prime} = T(AW^{\prime} - BW^{\prime}  \Lambda^{\prime})$.
Thus, instead of using the possibly complex basis $\left[ V,\ W, \ S, \ P \right]$, constructed 
by Algorithm~\ref{alg:bplmr} for the next BPLHR iteration,
one can set up the real basis $\left[ V^{\prime},\ W^{\prime}, \ S^{\prime}, \ P^{\prime} \right]$ 
that defines the same trial subspace. The block $P^{\prime}$ is constructed exactly in the same way 
as $P$ in Algorithm~\ref{alg:bplmr}, with the difference that the coefficients given by $Y$ are replaced
by those in $Y^{\prime}$. As a result, $P^{\prime}$ represents a combination of the new $V^{\prime}$ and 
of the subblock of $Z$ that corresponds to the current eigenvector approximations. Clearly, this definition
ensures that $P^{\prime}$ is real, and hence the same trial subspace is spanned by $\left[ V,\ W, \ S, \ P \right]$
and $\left[ V^{\prime},\ W^{\prime}, \ S^{\prime}, \ P^{\prime} \right]$.      

%The details of the resulting scheme are summarized in Algorithm~\ref{alg:bplhr_real}.   
The above discussion is based on the assumptions that the matrix $Y$ of eigenvectors 
of~\eqref{eqn:THarm_evp} contains the conjugate $\bar y$ for each complex column $y$.
In practice, however, this may not necessarily be the case. Since $k$ is a fixed parameter,
one can encounter the situation where all the complex columns $y$ have the corresponding conjugate
columns $\bar y$ in $Y$, but there exists a column $\tilde y = \tilde y_{\text{R}} + i \tilde y_{\text{I}}$
whose conjugate is not among the $k$ columns of $Y$. This happens if $\tilde y$ 
corresponds to the $k$th eigenvalue of~\eqref{eqn:THarm_evp}, whereas the conjugation of $\tilde y$ corresponds
to the $(k+1)$st eigenvalue, where the eigenvalues are numbered in the ascending order of their magnitudes.   

In this case, we ignore the imaginary part of $\tilde y$ and construct the matrix $Y^{\prime}$ to be 
of the form $Y^{\prime} = [Y_0 \ Y_{\text{R}}\ \tilde y_{\text{R}} \ Y_{\text{I}}]$, where 
$Y_{\text{R}}$ and $Y_{\text{I}}$ contain the real and imaginary parts of those vectors $y$ that
appear in $Y$ in conjugate pairs.
%, i.e., such that both $y$ and $\bar y$ are in $Y$. 
Given $Y^{\prime}$,
we proceed in the same fashion as has been described above. We introduce the block 
$V^{\prime} \equiv [V_0 \ V_{\text{R}}\ \tilde v_{\text{R}} \ V_{\text{I}}]$ of approximate eigenvectors
and the diagonal matrix $\Lambda^{\prime} = \text{diag}(\Lambda_0, \Lambda_1, \tilde \lambda, \Lambda_1)$
of the corresponding eigenvalue approximations, where  
$\tilde v_{\text{R}} = Z \tilde y_{\text{R}}$ and $\tilde \lambda$ is the RQ evaluated at $\tilde v_{\text{R}}$,
i.e., $\tilde \lambda \equiv \lambda(\tilde v_{\text{R}}) = (\tilde v_{\text{R}}, A \tilde v_{\text{R}})/
(\tilde v_{\text{R}}, B \tilde v_{\text{R}})$. The blocks $W^{\prime}$ and $S^{\prime}$ then take the form
$W^{\prime} = T(AV^{\prime} - B V^{\prime} \Lambda^{\prime}) = [W_0 \ W_{\text{R}}\ \tilde w_{\text{R}} \ W_{\text{I}}]$
and $S^{\prime} = T(AW^{\prime} - B W^{\prime} \Lambda^{\prime}) = [S_0 \ S_{\text{R}}\ \tilde s_{\text{R}} \ S_{\text{I}}]$,
where $\tilde w_{\text{R}} = A \tilde v_{\text{R}} - \tilde \lambda B v_{\text{R}}$ and $\tilde s_{\text{R}} = 
A \tilde w_{\text{R}} - \tilde \lambda B \tilde w_{\text{R}}$.  

Clearly, if the imaginary part of $\tilde y$ is ignored, i.e., the vector 
$\tilde v_{\text{I}} = Z \tilde y_{\text{I}}$ is not in~$V^{\prime}$, 
then the columns of  $[V^{\prime}, \ W^{\prime}, \  S^{\prime}, \ P^{\prime}]$ span a slightly smaller subspace 
than those of $[V, \ W, \  S, \ P]$ in Algorithm~\ref{alg:bplmr}. This may lead to a slight convergence deterioration
of the real arithmetic scheme compared to the original BPLHR version. Therefore, prior to the run, 
we suggest increasing the 
block size (at least) by 1 and running the algorithm for the extended block size. This removes the 
possible effects of ``cutting'' in between the conjugate pair when eigenvectors of~\eqref{eqn:THarm_evp} are selected into $Y$.

\begin{algorithm}[htbp]
\begin{small}
\begin{center}
  \begin{minipage}{5in}
\begin{tabular}{p{0.5in}p{4.5in}}
{\bf Input}:  &  \begin{minipage}[t]{4.0in}
A symmetric matrix $A$, an SPD matrix $B$, an SPD preconditioner $T$,
the shift $\sigma$, and the real initial guess $V^{(0)}$ for $k$ eigenvectors;
                  \end{minipage} \\
{\bf Output}:  &  \begin{minipage}[t]{4.0in}
                 Diagonal matrix $\Lambda$ of eigenvalues closest to the shift $\sigma$
                 and the matrix $V$ of the associated eigenvectors;
                  \end{minipage}
\end{tabular}
\begin{algorithmic}[1]

\STATE $V \gets V^{(0)}$; $P \gets [ \ ]$; 
\STATE Normalize columns of $V$ to have a unit $B$-norm; $\Lambda \gets \mbox{diag}\left( V^* A V \right)$;  
\WHILE {convergence not reached}
  \STATE Compute the preconditioned residuals $W \gets T(AV - B V \Lambda)$;
  \STATE Compute $S \gets  T (A W - B W \Lambda)$;
  \STATE Set $Z \gets \left[ V,\ W, \ S, \ P \right]$. 
$B$-orthonormalize the columns of $Z$. 
Let $\hat Z = [ \hat V,\ \hat W, \ \hat S, \ \hat P ]$ be the matrix of the resulting
$B$-orthonormal columns. 
  \STATE Find eigenpairs of the projected problem~\eqref{eqn:THarm_evp} with $Z \equiv \hat Z$. 
%Define $Y \equiv [Y_V^T, \ Y_W^T, \ Y_S^T, \ Y_P^T]^T$
Sort the eigenvalues $\xi$ in the ascending order of their absolute values, ensuring that in the 
sorted set every complex eigenvalue $\xi$ is immediately followed by its conjugate $\bar \xi$.  
\STATE Select $k$ eigenvectors of \eqref{eqn:THarm_evp} associated with the smallest magnitude eigenvalues into 
the matrix $Y$.
  \STATE Move real columns of $Y$ into $Y_0$. Define $Y_{\text{R}}$ and $Y_{\text{I}}$ to be the matrices 
  containing the real and imaginary parts of the columns of $Y$ that appear in complex pairs.
  \STATE If $Y$ contains a complex column $\tilde y$, such that the conjugate of $\tilde y$ is not in $Y$,
  then set $\tilde y_{\text{R}}$ to be the real part of $\tilde y$. Otherwise, $\tilde y_{\text{R}} \gets [ \ ]$.  
%such that if $y$ is a column of $Y$ then $\bar y$ is also a column of $Y$.
%  \STATE Use $Y$ to form the real matrix 
  Define $Y^{\prime} = [Y_0 \ Y_{\text{R}} \ \tilde y_{\text{R}} \ Y_{\text{I}}]
 \equiv [Y_V^T, \ Y_W^T, \ Y_S^T, \ Y_P^T]^T$.
%         where
%         $Y_0$ contains the real columns of $Y$; and $Y_{\text{R}}$ and $Y_{\text{I}}$ correspond to the real
%         and imaginary parts of the complex columns of $Y$, respectively, in the corresponding order.
%     \STATE  
%Define the rowwise partitioning $Y^{\prime} \equiv [Y_V^T, \ Y_W^T, \ Y_S^T, \ Y_P^T]^T$.
  \STATE  Compute $P \gets \hat W Y_W + \hat S Y_S + \hat P Y_P$;
  \STATE  Compute new approximate eigenvectors $V \gets \hat V Y_V + P$.
    Define the column partitioning $V = [V_0 \ V_{\text{R}} \ \tilde v_{\text{R}}  \ V_{\text{I}}] (\equiv \hat Z Y^{\prime})$
          according to $Y^{\prime}$ ($\tilde v_{\text{R}} \gets [ \ ]$ if $\tilde y_{\text{R}} = [ \ ]$). 
\STATE For each column $v$ of $V_0$ compute the RQ $\lambda(v) \gets v^* A v/v^* B v$ and place it to the diagonal
of $\Lambda_0$. For each column $v_{\text{R}}$ of $V_{\text{R}}$ and the corresponding $v_{\text{I}}$ of $V_{\text{I}}$, 
compute the diagonal entry of $\Lambda_1$ by~\eqref{eqn:lambda_complex}. 
Set $\tilde \lambda \gets \tilde v_{\text{R}}^* A \tilde v_{\text{R}}/\tilde v_{\text{R}}^* B \tilde v_{\text{R}}$
($\tilde \lambda \gets [ \ ]$ if $\tilde v_{\text{R}} = [ \ ]$). 
%  $\Lambda_0 \gets \mbox{diag}\left( V_0^* A V_0 \right)$ 
  \STATE  Normalize columns of $V$ to have a unit $B$-norm; 
$\Lambda \gets \mbox{diag}( \Lambda_0, \ \Lambda_1, \ \tilde \lambda, \  \Lambda_1 )$. 
\ENDWHILE
\STATE Perform the standard RR procedure for~\eqref{eqn:eigp} with respect to $V$. 
Update $V$ to contain the Ritz vectors and $\Lambda$ the corresponding Ritz values.  
\STATE Return $(\Lambda, V)$.
\end{algorithmic}
\end{minipage}
\end{center}
\end{small}
  \caption{The real arithmetic version of the BPLHR algorithm}
  \label{alg:bplhr_real}
\end{algorithm}

The details of the real arithmetic version of BPLHR algorithm are summarized in Algorithm~\ref{alg:bplhr_real}. 
Note that if at each step of the algorithm problem~\eqref{eqn:THarm_evp} gives only real eigenvectors $Y$, then 
Algorithm~\ref{alg:bplhr_real} becomes equivalent to the original version of BPLHR in Algorithm~\ref{alg:bplmr}.

\section{Preconditioning}\label{sec:prec}

In order to motivate preconditioning strategy,
let us again consider the idealized trial subspaces~\eqref{eqn:plmr_ideal_span}. 
Assuming that the targeted eigenvalue $\lambda_q$ is known, we address 
the question of defining an optimal preconditioner $T$ which ensures 
that the corresponding eigenvector $v_q$ is exactly in 
the trial subspace~\eqref{eqn:plmr_ideal_span}.   

%A possible choice of such a preconditioner, pointed out in [?], is given by 
A possible choice of such a preconditioner is given by 
$T = (A - \lambda_q B)^\dagger$~\cite{Knyazev:98}. 
In this case,
\[
v_q \equiv v^{(i)} - T (A  - \lambda_q B) v^{(i)} = (I - K)v^{(i)}, 
\]
where 
$K = (A - \lambda_q B)^\dagger (A  - \lambda_q B )$.
Clearly, $v_q$ is in~\eqref{eqn:plmr_ideal_span}. The fact that 
$v_q$ is an eigenvector follows from the observation that 
$K$ is an orthogonal projector onto the range of $A - \lambda_q B$. 
Hence, $I - K$ projects $v^{(i)}$ onto the null space of $A - \lambda_q B$,
which gives the desired eigenvector.
%, which is clearly 
%in~\eqref{eqn:plmr_ideal_span}. 
Note that $v_q$ is nonzero provided that $v^{(i)}$ has a nontrivial 
component in the direction of the targeted eigenvector. 
%Clearly, $v_q$ is in~\eqref{eqn:plmr_ideal_span}.
%  \text{sign}(A - \lambda_q B)

In general, neither the idealized subspaces~\eqref{eqn:plmr_ideal_span} nor
the optimal preconditioner $T = (A - \lambda_q B)^\dagger$ are available
at the PLHR iterations. However, the above analysis suggests that practical 
preconditioners can be defined as approximations of $(A - \lambda_q B)^\dagger$.
For example, one can aim at constructing $T$, such that $T \approx (A - \sigma B)^{-1}$.
This gives rise to the inexact shift-and-invert type preconditioning,
which represents a traditional approach 
for preconditioning eigenvalue problems; see, e.g.,~\cite{Golub.Ye:02, Knyazev:01}.    

Since preconditioners for PLHR should be HPD, the definition
$T \approx (A - \sigma B)^{-1}$ may not always be suitable. 
While one can expect to construct the \textit{HPD} inexact shift-and-invert 
preconditioners for computing eigenpairs associated with the eigenvalues 
not far away from the ends of the spectrum,
%$T \approx (A - \sigma B)^{-1}$ 
the approach will result in the indefinite preconditioning  
as the eigenvalues deeper in the interior of the spectrum are sought.

%A traditional way to define a preconditioner for Algorithms~\ref{alg:plmr}
%and~\ref{alg:bplmr} would be to choose $T$ as a form of an approximate inverse 
%of the shifted matrix $A - \sigma B$; see~\cite{Golub.Ye:02, Knyazev:01}. 
%Since $A - \sigma B$ is
%generally symmetric indefinite, it is natural to expect that a preconditioner 
%$T \approx (A - \sigma B)^{-1}$ will also be indefinite. 
%This, however, violates our assumption on positive definiteness of $T$
%and, as a consequence, does not allow to define the $T$-harmonic 
%projection procedure. 
% used by the proposed techniques.   

In order to define a preconditioner that is HPD 
%regardless of the targeted eigenpair 
and, at the same time, that preserves the desirable effects of 
``shift-and-invert'', let us consider the operator $T = |A - \lambda_q B|^\dagger$,
where $|A - \lambda_q B|$ is defined as a matrix function~\cite{Golub.VanLoan:96} of 
$A - \lambda_q B$. 
Similar to $T = (A - \lambda_q B)^\dagger$, 
this choice of preconditioner is also optimal with respect to the idealized trial subspaces~\eqref{eqn:plmr_ideal_span}. Indeed, if $T = |A - \lambda_q B|^\dagger$, then
\[
v_q \equiv v^{(i)} - T(A - \lambda_q B)T(A - \lambda_q B)v^{(i)} = (I - K) v^{(i)},
\] 
where 
$K = \left( |A - \lambda_q B|^\dagger (A - \lambda_q B) \right)^2 = 
%\left( \text{sign}(A - \lambda_q B) (A - \lambda_q B)^\dagger (A - \lambda_q B) \right)^2 = $ $\text{sign}^2(A - \lambda_q B) = 
\left( (A - \lambda_q B)^\dagger (A - \lambda_q B) \right)^2 = 
(A - \lambda_q B)^\dagger (A - \lambda_q B)$. Thus, $v_q$ 
belongs to~\eqref{eqn:plmr_ideal_span} and represents the wanted eigenvector, since
$I - K$ projects $v^{(i)}$ onto the null space of $A - \lambda_q B$.

%This matrix has the same set of eigenvectors as $A - \sigma B$, whereas its 
%eigenvalues are equal to the absolute values of eigenvalues of $(A - \sigma B)^{-1}$,
%i.e., $\lambda_j(|A - \sigma B|^{-1}) = | \lambda_j((A - \sigma B)^{-1}) |$. 
%Clearly, $|A - \sigma B|^{-1}$ is SPD and, hence, constructing  
%$T \approx |A - \sigma B|^{-1}$ can result in SPD preconditioners, as required
%by the introduced algorithms.  

As discussed above, the computation of idealized trial subspaces and 
an optimal preconditioner may be infeasible in practice. Instead, however,
one can construct $T$ that \textit{approximates} the action of the (pseudo-) inverted
absolute value operator. In particular, we can define $T \approx |A - \sigma B|^{-1}$.
%%%%%%%%%%REFEREE
In the next section, we demonstrate 
that such preconditioners exist for certain classes of problems
% that available preconditioners
%representing such approximations can lead to 
and indeed lead to a rapid and robust convergence if used to 
construct the PLHR trial subspaces~\eqref{eqn:plmr_span} 
(or~\eqref{eqn:plmr_span_block}, in the block case) and to perform the $T$-harmonic 
projection introduced in Section~\ref{subsec:THarm}.

\section{Numerical experiments}\label{sec:numr}

We organize our numerical experiments into three sets. The first set concerns 
a model generalized eigenvalue problem of a small size. This example is mainly
of theoretical nature and is intended to demonstrate several features of the 
convergence behavior of the proposed method.

In the second test set, we consider a larger model problem, where a practical
AV preconditioner is used. Specifically, we address a problem of computing
a subset of interior eigenpairs of a discrete Laplacian. As an SPD preconditioner
we use the multigrid (MG) scheme proposed in~\cite{Ve.Kn:13} in the context
of solving symmetric indefinite Helmholtz type linear systems. 
We show that exactly the same
AV preconditioner can be employed for the interior eigenvalue computations in 
%the framework of 
the BPLHR algorithm.   

The third series of experiments aims at a particular application.
We consider several Hamiltonian matrices that arise in electronic
structure calculations, and apply BPLHR to reveal eigenvectors (wave functions) 
that correspond to the eigenvalues (energy levels) around a given reference 
energy. The preconditioners available for this type of problems are diagonal and 
SPD~\cite{Teter.Payne.Allan:89}.   
%All matrices in these examples are generated using the KSSOLV {\sc matlab} 
%toolbox~\cite{kssolv:09} for solving Kohn--Sham equations.
   
In our experiments, we compare BPLHR with the a block version of the 
Generalized Davidson (BGD) method that is based on the harmonic 
projection~\cite{Jordan.Marsman.Kim.Kresse:12,Morgan:91}.
This choice has been motivated by several factors. 
First, we want to restrict our comparisons to the class of block methods
which rely on the same type of computational kernels, 
such as multiplication of a \textit{block} of vectors by a matrix or a preconditioner,
dense matrix-matrix multiplication, etc.
%
%As has been discussed in introduction, our interest in block iterations is motivated
%by their capability to leverage the additional degree of parallelism introduced by 
%the block operations when implemented on parallel computing platforms.
%
Secondly, the BGD algorithm represents a state-of-the-art approach for interior
eigenvalue computations in a number of critical applications; 
see, e.g.,~\cite{Jordan.Marsman.Kim.Kresse:12}.     
%Therefore, our comparisons may result in recommendations on when the new
%method can be preferable than the well-established GD scheme. In particular, as 
%shown 
In particular, 
%subsequent sections, 
we demonstrate that BPLHR can give a more robust solution option compared to BGD 
in cases where memory is limited and the available 
preconditioner is of a ``moderate'' quality. 
%Note that all our tests are performed in {\sc matlab}.  

\subsection{A small model problem: the finite element (FE) Laplacian}
\label{subsec:small}

Let us consider the eigenvalue problem 
\begin{equation}\label{eqn:evp_cont}
-\Delta u (\mathrm x, \mathrm y)  =  \lambda u(\mathrm x,\mathrm y), \  (\mathrm x,\mathrm y) \in \Omega  = (0,1)\times(0,1), \ u|_\Gamma = 0,    
\end{equation} 
where $\displaystyle \Delta = \partial^2/\partial \mathrm x ^2 + \partial^2 /\partial \mathrm y^2$ is the Laplace operator and $\Gamma$ denotes the boundary of the
unit square $\Omega$. Assume that we are interested in approximating
a number of eigenpairs $(\lambda, u(\mathrm x,\mathrm y))$ of~\eqref{eqn:evp_cont}
that are closest to a given shift $\sigma$.

Discretizing~\eqref{eqn:evp_cont} with standard bilinear finite elements
results in the algebraic generalized eigenvalue problem $A v = \lambda B v$,
where $A$ and $B$ are the SPD stiffness and mass matrices, respectively.
In this example, we assume a relatively small number of finite
elements, $50$ along each side of the unit square, which results in the 
total of $2,401$ degrees of freedom, i.e., the size of the matrix problem is 
$n = 2,401$.

Due to the small problem size, in this example, we can model high-quality AV preconditioners
by perturbations of the operator $|A - \sigma B|^{-1}$. 
In particular, we define $T = |A - \sigma B|^{-1} + E$,
where $E$ is a random SPD matrix, such that 
$\|E\| \leq \epsilon \|(A - \sigma B)^{-1}\|$ and $\epsilon$
is a parameter that sets up the preconditioning quality; $\|\cdot\|$ denotes the spectral norm.  
%%%%%%%%%%%%%%REFEREE
As $\epsilon$ increases, so does the norm of the perturbation $E$, which means that
$T$ becomes more distant from $|A - \sigma B|^{-1}$ and hence the quality of the preconditioner deteriorates.
%%%%%%%%%%%%%%
The matrix $E$ is fixed during the iterations, but is updated for every new run.
We present typical results. 

%In our first test, 
%is to demonstrate the following. First, 
We first compare the convergence rate of PLHR to that of the 
``base'' preconditioned null space finder~\eqref{eqn:plmr_ideal}, with~\eqref{eqn:coeff} 
and~\eqref{eqn:subsp2}, applied to \eqref{eqn:eigp_ideal}.  
As discussed in Section~\ref{sec:idealized}, this scheme, further referred to as BASE-NULL,
represents an idealized interior eigenvalue solver with a proven convergence 
bound~\eqref{eqn:cv_ideal_nonopt}--\eqref{eqn:cv_ideal_kappa}, and was used as a prototype of the PLHR
algorithm. In particular, PLHR was derived in Section~\ref{sec:plmr} 
%Recall that BASE-NULL has a   
%%this linear solver, with a 
%proven convergence with bound~\eqref{eqn:cv_ideal_nonopt}--\eqref{eqn:cv_ideal_kappa}, and 
%%considered in section~\ref{sec:idealized}, 
%was used as a prototype of the PLHR algorithm.
%In particular, PLHR was obtained in section~\ref{sec:plmr} 
by introducing a number of approximations into BASE-NULL. 
Therefore, it is of interest to see to what extent  
%how 
these approximations affect the behavior of the resulting eigenvalue algorithm
if the understood BASE-NULL convergence is taken as a reference.   
\begin{figure}[ptbh]
\begin{center}%
\begin{tabular}
[c]{cc}%
	\includegraphics[width=6.5cm]{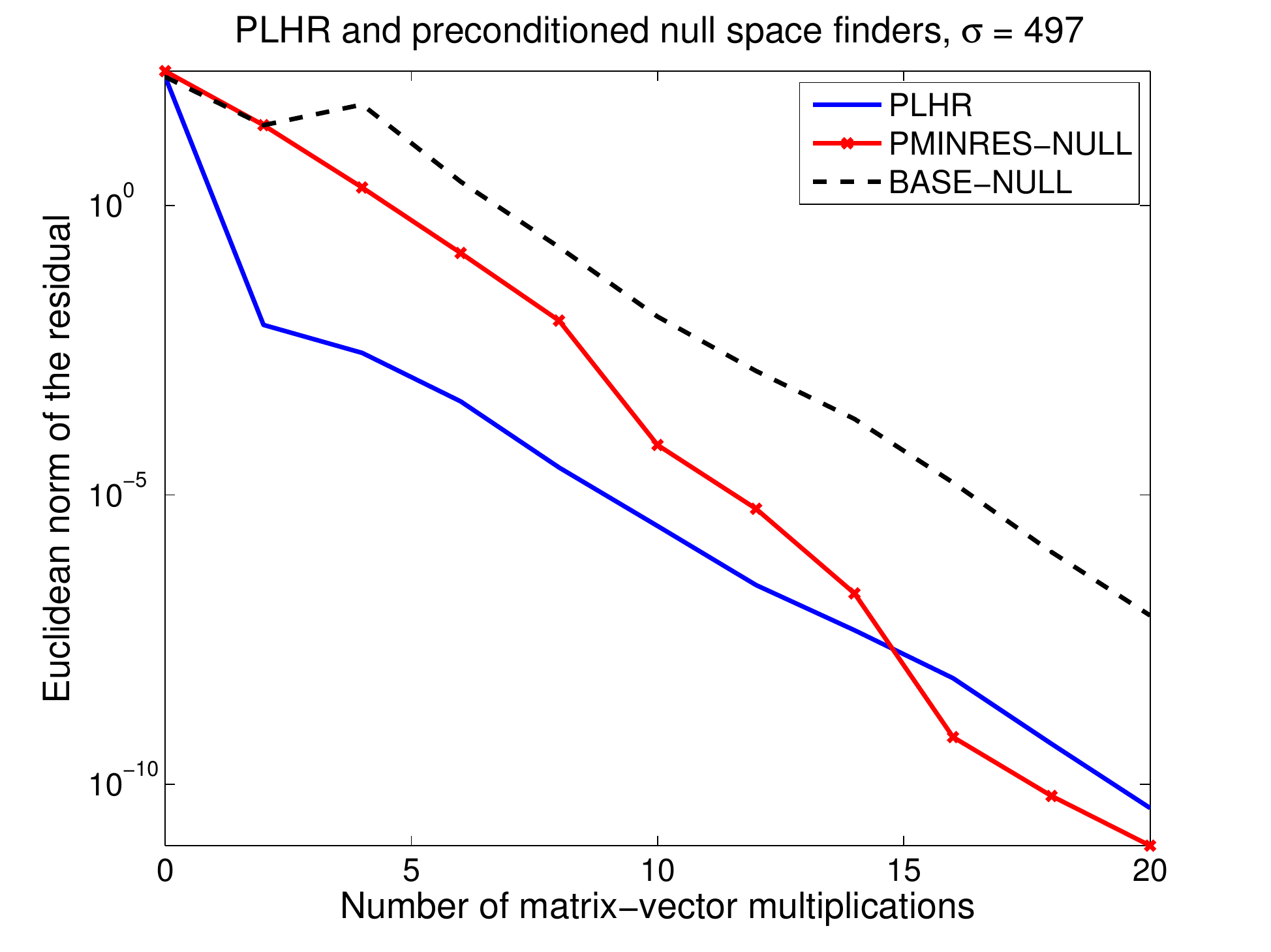}
	\includegraphics[width=6.5cm]{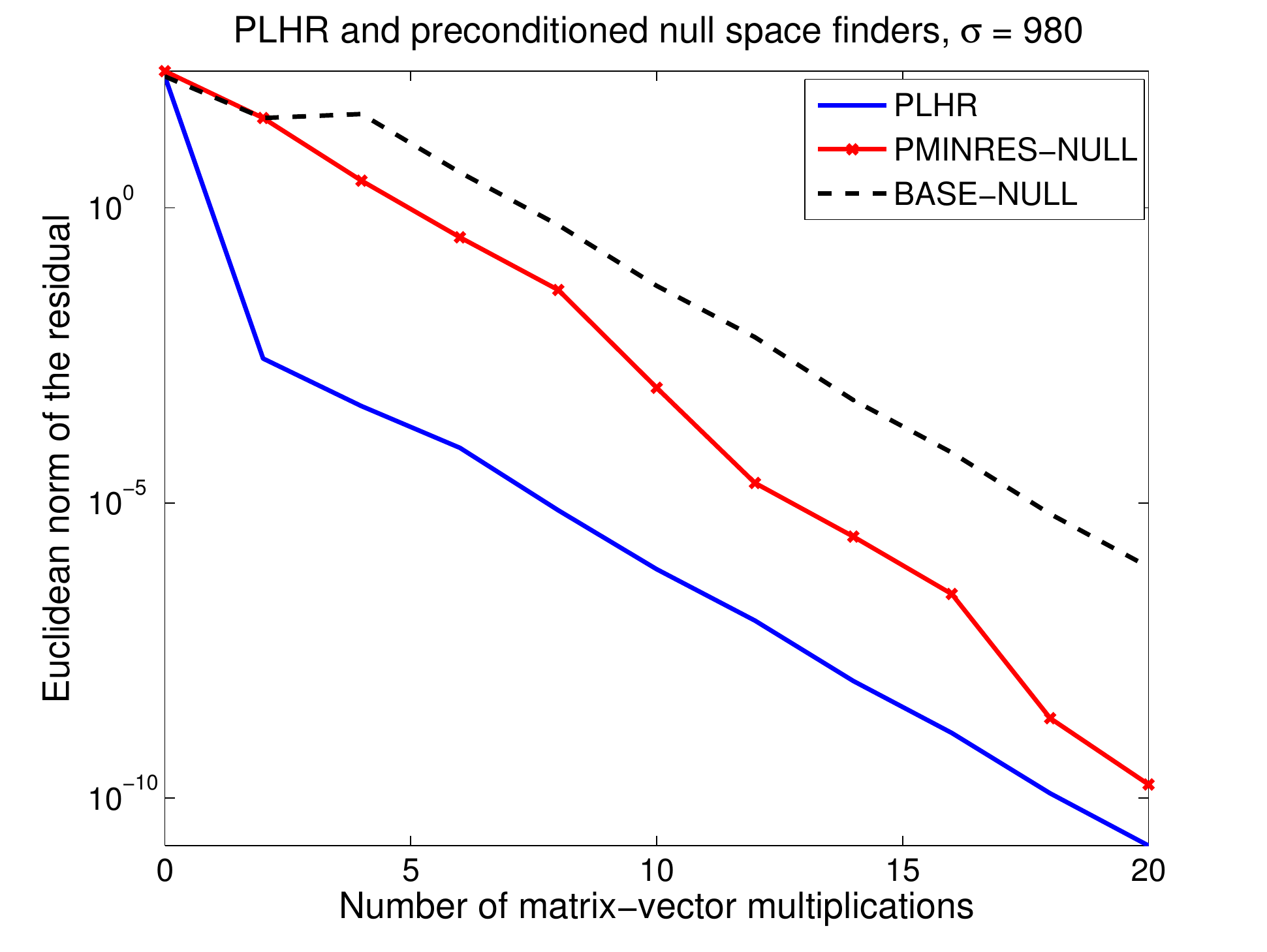}
\end{tabular}
\end{center}
\caption{A comparison of PLHR with 
%to the convergence of
preconditioned linear solvers for 
%the singular linear system 
$(A - \lambda_{31} B)x=0$ (left) and $(A - \lambda_{66} B)x=0$ (right), where 
$\lambda_{31} \approx 497.5521$ and $\lambda_{66} \approx 979.7072$.}
\label{fig:null_find}%
\end{figure}

In Figure~\ref{fig:null_find}, we apply PLHR (Algorithm~\ref{alg:plmr}) and BASE-NULL
to compute an eigenpair closest to the shift
%$\sigma$. 
%
%In particular, we choose 
$\sigma = 497$ and $\sigma = 980$. These shift values target
eigenpairs corresponding to the eigenvalues $\lambda_{31} \approx 497.5521$ and 
$\lambda_{66} \approx 979.7072$, respectively. 
%
%The convergence of PLHR is compared to the convergence of the above mentioned 
%idealized ``base'' locally optimal null space finder (denoted ``LO-BASE-NULL'')
%applied to systems with matrices 
Thus, BASE-NULL solves the homogeneous systems with matrices
$A - \lambda_{31} B$ and $A - \lambda_{66} B$. 
%Note that, in this test, we use the PLHR version given by Algorithm~\ref{alg:plmr},
%which can 
We also plot convergence curves that correspond to the runs of PMINRES (denoted PMINRES-NULL) 
applied to the same singular systems.
This gives us an opportunity to compare the convergence of PLHR to that of an optimal
Krylov subspace method used as an idealized eigenvalue solver; see the discussion in Section~\ref{sec:idealized}.
%To assess the convergence, for all schemes, we measure
%the norms of the residuals, $\|Av^{(i)} - \lambda^{(i)}B v^{(i)}\|$,
%of the eigenvalue problem.   

To assess the convergence, for all schemes in Figure~\ref{fig:null_find}, we measure
the norms of the residuals $\|Av^{(i)} - \lambda^{(i)}B v^{(i)}\|$ of the eigenvalue problem.   
In the definition of the random perturbation based preconditioner $T$,  
we set $\epsilon = 10^{-5}$.

Figure~\ref{fig:null_find} shows that PLHR and BASE-NULL exhibit essentially the same convergence behavior. 
Additionally, at a number of initial steps, their convergence is similar to that of
PMINRES-NULL. 
Thus, at least if the preconditioner is sufficiently strong, the approximations introduced
into the idealized BASE-NULL to obtain PLHR do not significantly alter its convergence behavior.
Moreover, the convergence is comparable to that of the optimal PMINRES-NULL.       
%, at a number of its initial steps, 
%also has a similar convergence behavior, however, accelerates at some point 
%due to the method's global optimality.  
%, possibly, with the occurrence
%of the superlinear convergence, which is frequently noticed for preconditioned globally optimal
%Krylov subspace methods, e.g., \cite{Beckermann.Kuijlaars:01, Simoncini.Szyld.05a}. 
%In fact, it is generally hard to expect the superlinear convergence for 
%the PLHR algorithm, which is likely to be the price, paid for the departure from 
%the \textit{global} optimality.

\begin{figure}[ptbh]
\begin{center}%
\begin{tabular}
[c]{cc}%
	\includegraphics[width=6.5cm]{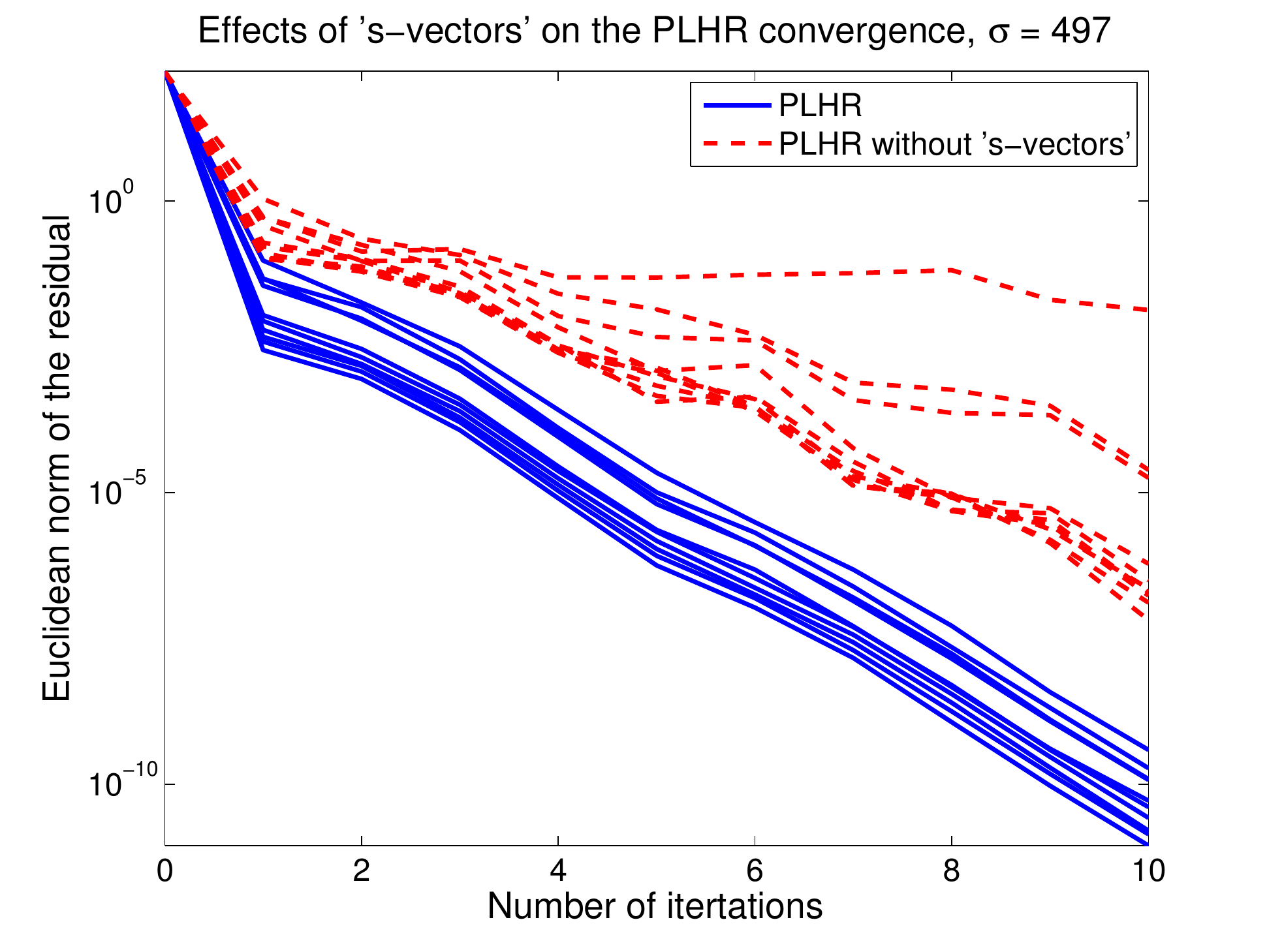}
	\includegraphics[width=6.5cm]{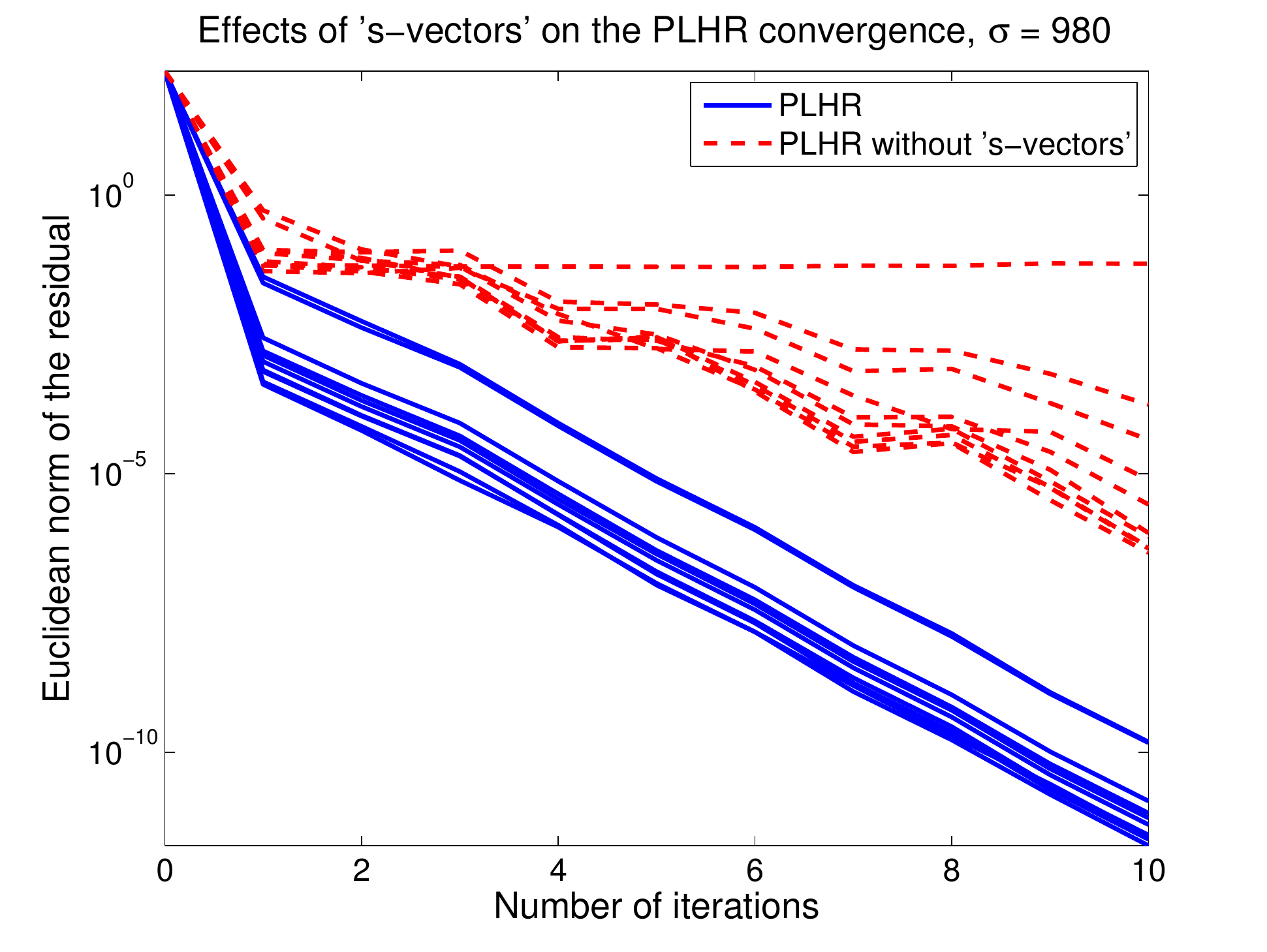}
\end{tabular}
\end{center}
\caption{The significance of ``s-vectors'' for the PLHR convergence.
The PLHR algorithm with $\sigma = 497$ (left) and $\sigma = 980$ (right) is
compared to its version without ``s-vectors''.
Each curve corresponds to a separate run with random initial guess and a preconditioner generated by a random 
SPD perturbation of $|A - \sigma B|^{-1}$ with~$\epsilon = 10^{-5}$.}  
%Convergence of PLHR compared to preconditioned null space 
%finders; 
%$\lambda_{31} \approx 497.5521$ (left) and $\lambda_{66} \approx 979.7072$ 
%(right) . 
%Prec. rel. pert = $10^{-5}$; 10 runs}%
\label{fig:plmr_noS}%
\end{figure}

Next, we would like demonstrate the effects of the vectors $s^{(i)}$ 
introduced into the PLHR trial subspaces~\eqref{eqn:plmr_span}. These 
subspaces can be viewed as the LOBPCG-like subspaces, spanned by $v^{(i)}$, $w^{(i)}$, and $p^{(i)}$,
extended by the additional vector. Therefore, a natural question is whether the occurrence of 
the new ``s-vectors'' has any impact on the convergence of the proposed scheme. 

%An important novelty of the proposed technique is the introduction of 
%``s-vectors'' $s^{(i)}$ into the trial subspaces. 
%In Figure~\ref{fig:plmr_noS}, we demonstrate that the presence of such vectors
%is indeed crucial for the method's stability and robustness. 
%In particular, we compare the PLHR algorithm with its variant that has no
%``s-vector'' in the trial subspaces, i.e., performs the $T$-harmonic
%extraction from the (LOBPCG-like) subspaces spanned only by $v^{(i)}$, $w^{(i)}$, 
%and $p^{(i)}$.     

Figure~\ref{fig:plmr_noS} compares the PLHR iteration in Algorithm~\ref{alg:plmr} 
to its variant where the ``s-vectors'' are not included into the trial subspaces, i.e.,~\eqref{eqn:plmr_span}
are spanned only by $v^{(i)}$, $w^{(i)}$, and $p^{(i)}$. We perform $10$ runs of both versions of the algorithm,
so that each curve in Figure~\ref{fig:plmr_noS} represents a separate execution with a random initial guess.
As in the previous example, we consider shifts $\sigma = 497$ and $\sigma = 980$. 
For each run, the preconditioners are given by random SPD perturbations of 
$|A - \sigma B|^{-1}$ with $\epsilon = 10^{-5}$. 

One can observe from Figure~\ref{fig:plmr_noS} that PLHR demonstrates a stable linear convergence at all runs
regardless of the initial guess and a particular instance of the preconditioner.
At the same time, despite the high preconditioning quality, the absence of ``s-vectors'' makes the method highly 
unstable, with a slower or stagnant convergence pattern. Therefore, the presence of $s^{(i)}$ 
in~\eqref{eqn:plmr_span} is important. 
%%%%% REFEREE
Such a behavior is consistent with the fact that linear solver~\eqref{eqn:plmr_ideal_short} 
generally does not converge if $s$-vectors, defined as $T\left(A  - \lambda_q B \right) T r^{(i)}$, 
are removed from the iterative scheme. 
%where they are
%defined, using the exact eigenvalue, as $T\left(A  - \lambda_q B \right) T r^{(i)}$.

%{\small
%\begin{table}
%\centering
%\begin{tabular}{|c|c|c|}
%\hline 
%The PLHR variant & Preconditioner & RR \tabularnewline
%\hline 
%\hline 
%(B)PLHR & SPD  & $T$-harmonic\tabularnewline
%\hline 
%(B)PLHR-HARM-AV & SPD & Harmonic\tabularnewline
%\hline 
%(B)PLHR-HARM-INV & Indefinite & Harmonic\tabularnewline
%\hline 
%\end{tabular}
%\caption{The summary of PLHR variants based on different choices of preconditioner and 
%extraction procedure.}\label{tab:plmrs}
%\end{table}
%}

Another new feature incorporated into PLHR is the $T$-harmonic 
RR procedure presented in Section~\ref{subsec:THarm}. Similar to the above, 
we would like to address the question of whether any advantage is gained 
by the $T$-harmonic approach compared, e.g., to the \textit{standard} harmonic RR~\cite{Morgan:91}. 
%orporating the preconditioner into the Petrov-Galerkin
%condition~\eqref{eqn:orth_gen}.
%% for the eigenvector extraction. 
%In order to verify this, 
To answer this question, let us compare the PLHR algorithm to
its variant where the $T$-harmonic projection is replaced by the standard harmonic RR,
whereas the same SPD (AV) preconditioner is used to generate the trial subspaces.
%We refer to this modification as PLHR-HARM-AV.  

Since PLHR with a standard harmonic RR no longer requires the preconditioner to be SPD, 
%(due to the absence of the $T$-harmonic RR),
we are also interested in the case where $T$ is \textit{indefinite}, i.e., the preconditioner is
an approximation of the ``shift-and-invert'' operator $(A - \sigma B)^{-1}$. For this reason, 
let us consider a variant of PLHR with an indefinite $T \approx (A - \sigma B)^{-1}$, combined with the standard harmonic 
RR. 
%We call this version PLHR-HARM-INV. The differences between different variants of PLHR are outlined
%in Table~\ref{tab:plmrs}. 

\begin{figure}[ptbh]
\begin{center}%
\begin{tabular}
[c]{cc}%
	\includegraphics[width=6.5cm]{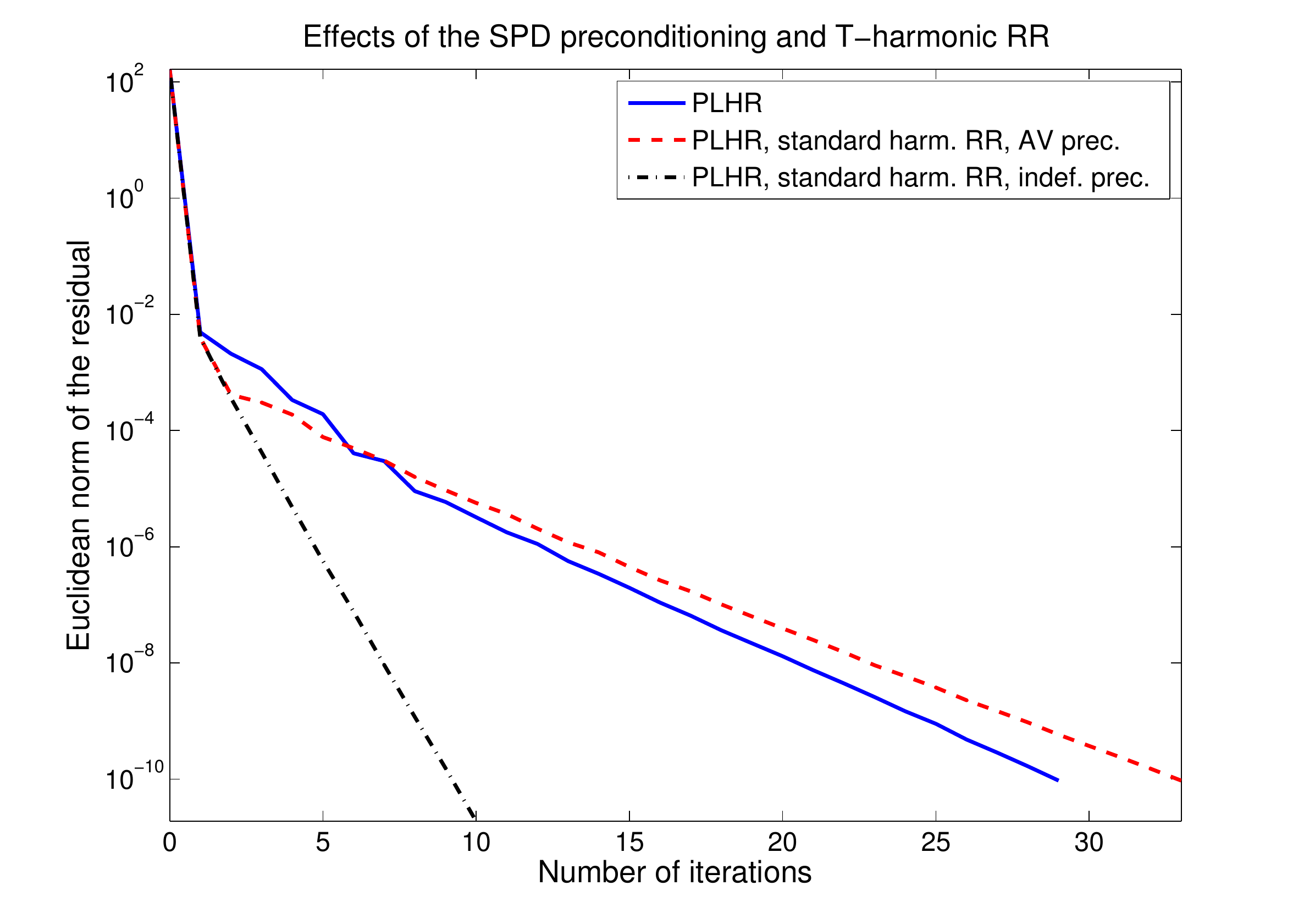}
	\includegraphics[width=6.5cm]{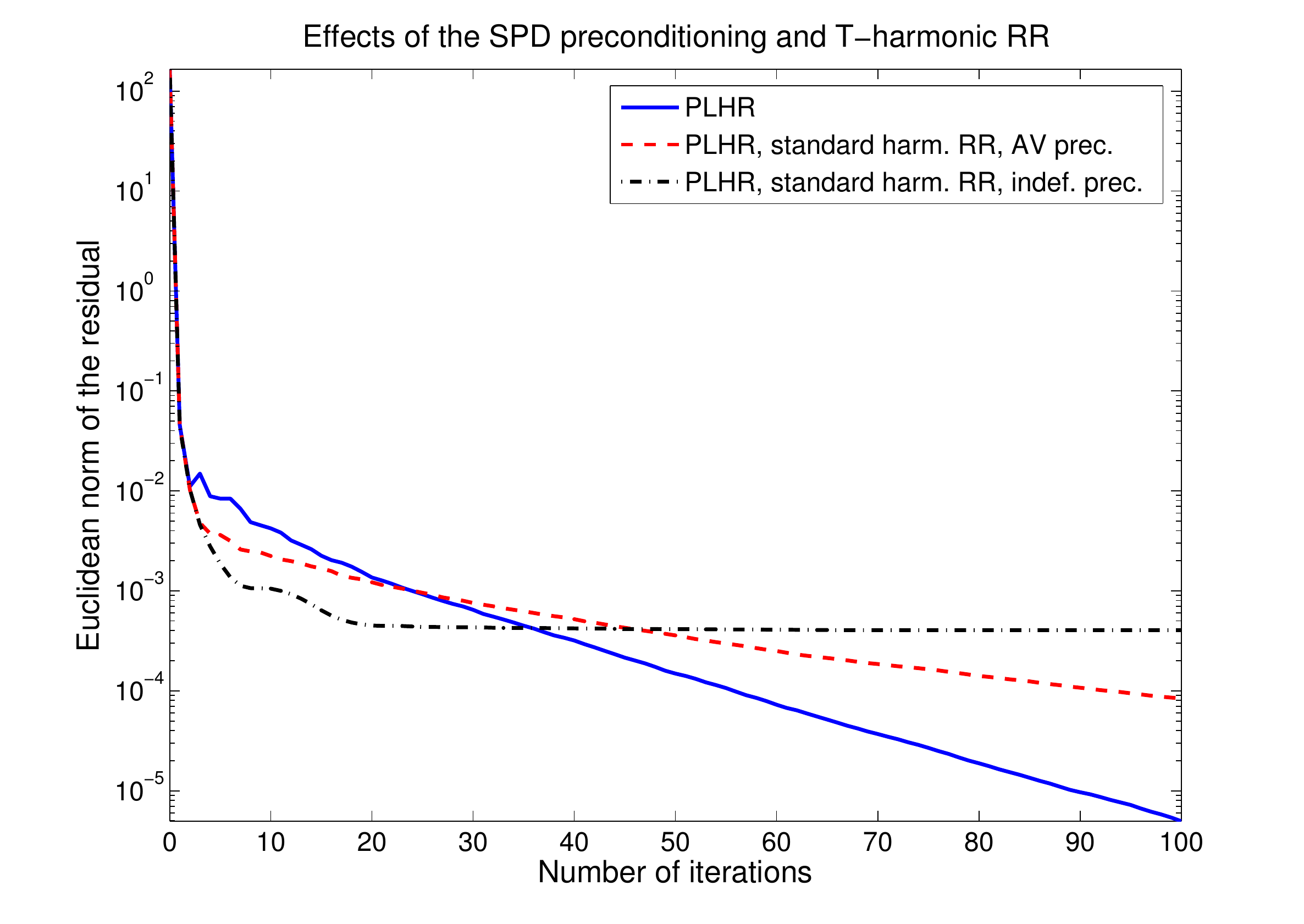}
\end{tabular}
\end{center}
\caption{Effects of the SPD AV preconditioning and $T$-harmonic RR; $\sigma = 980$. 
Preconditioners generated by a random SPD perturbation of 
$|A - \sigma B|^{-1}$ (for the AV preconditioner) 
and of $(A - \sigma B)^{-1}$ (for the indefinite preconditioner) with $\epsilon = 10^{-4}$
(left) and $\epsilon = 10^{-3}$ (right).}  
\label{fig:THarm}%
\end{figure}

%Both ``PLHR-harm-AV'' and ``PLHR-harm-INV''are compared to PLHR in 
%Figure~\ref{fig:THarm}. In this test, the AV preconditioner for PLHR
%and ``PLHR-harm-AV'' is again obtained by perturbing $|A - \sigma B|^{-1}$.
%The indefinite preconditioner for ``PLHR-harm-INV'', is similarly generated
%by a perturbation of $(A - \sigma B)^{-1}$, i.e., $T = (A - \sigma B)^{-1} + E$,
%where $E$ is random, such that $\|E\| \leq \epsilon \|(A - \sigma B)^{-1}\|$.
%In fact, the same SPD perturbation $E$ is used for both the AV and indefinite 
%preconditioner to maintain a similar preconditioning quality.
 
Figure~\ref{fig:THarm} illustrates the collective impact of the $T$-harmonic RR and 
AV preconditioning on the convergence of the new eigensolver.
Here, we set $\sigma = 980$ and apply all the three PLHR variants 
%in Table~\ref{tab:plmrs} 
to compute the corresponding eigenpair. 

In contrast to our previous tests, we now consider
preconditioners of different quality.    
In particular, in Figure~\ref{fig:THarm} (left) we set $\epsilon = 10^{-4}$,
%which gives a high-quality preconditioner, 
whereas in Figure~\ref{fig:THarm} (right) we choose $\epsilon = 10^{-3}$, which gives weaker preconditioners. 
%are obtained. 
As previously, the AV preconditioning 
%for PLHR and PLHR-HARM-AV'' 
is obtained by perturbing $|A - \sigma B|^{-1}$.
The indefinite preconditioner 
%in PLHR-HARM-INV 
is generated in a similar manner by
a random perturbation $E$ of $(A - \sigma B)^{-1}$, i.e., $T = (A - \sigma B)^{-1} + E$,
%where $E$ is random, 
where $\|E\| \leq \epsilon \|(A - \sigma B)^{-1}\|$.
In our tests, the same SPD perturbation $E$ is used for both the AV and indefinite preconditioners 
to maintain a similar preconditioning quality.

Figure~\ref{fig:THarm} (left) suggests that if a sufficiently strong preconditioner is at hand
then neither the $T$-harmonic projection nor the SPD preconditioning lead to any 
significant improvement. In this case, the fastest convergence is attained by the version of 
PLHR with the standard harmonic RR and indefinite preconditioner.
%; see Figure~\ref{fig:THarm} (left). 
%
However, if preconditioning is not as strong,
%, as in Figure~\ref{fig:THarm} (right), 
then the use of the SPD AV preconditioners along with the $T$-harmonic RR procedure becomes crucial. 
As can be seen in Figure~\ref{fig:THarm}~(right), the original PLHR version given by Algorithm~\ref{alg:plmr}
is the only scheme that is able to maintain convergence under the degraded preconditioning quality.   
%%%%%%%%%% REFEREE
Note that, in this example, the convergence of PLHR iterations could be preserved for values of $\epsilon$ up to $10^{-2}$.
For $\epsilon=10^{-2}$, the convergence behavior of all methods is similar to that in Figure~\ref{fig:THarm}~(right),
though with an increased iteration count due to preconditioning deterioration.

%the version ``PLHR-harm-INV'' 
%with indefinite preconditioning stagnates and the convergence of the 
%SPD-preconditioned ``PLHR-harm-AV'' slows down, so that the 
%PLHR algorithm results in the best run.       
 
%We are also interested in assessing the affects of the absolute value preconditioning 
%on the quality
%of trial subspaces, as opposed to using a standard preconditioning approach
%based on approximating the inverted operator $(A - \sigma B)^{-1}$. 

%As we  further demonstrate in our examples, 
The case where the preconditioner is only of a moderate quality 
is not uncommon in realistic applications, especially if the wanted eigenpairs are located deeper in the 
spectrum's interior. Thus, the decrease in the sensitivity to deterioration of the preconditioning
quality, demonstrated by PLHR in Figure~\ref{fig:THarm} (right),
%PLHR to maintain
%convergence in this situation, as demonstrated in Figure~\ref{fig:THarm} (right), 
is of a practical interest. In the remaining experiments, we reaffirm this finding 
on the example of several model problems for which AV preconditioners can be constructed in practice.      
%In many cases, 
%this is due to the fact that the construction of high-quality preconditioners 
%$T \approx (A - \sigma B)^{-1}$ is prohibitively
%costly in practice, especially for problems where $A - \sigma B$ 
%is highly indefinite. Thus, the approach given by PLHR may provide means for reliable computations of  
%interior eigenpairs with relatively inexpensive preconditioners. 

%Below, we show that this observation extends to the block 
%variant of PLHR (BPLHR) applied to more realistic examples for
%which practical SPD preconditioners are known and can be efficiently 
%constructed.   

\subsection{A model problem: the finite-difference (FD) Laplacian}
\label{subsec:laplace}
Let us consider the same continuous problem~\eqref{eqn:evp_cont}, but now apply a standard 5-point FD discretization
with step $h = 2^{-7}$. This gives an eigenvalue problem $L v = \lambda v$, where $A \equiv L$ 
is a discrete Laplacian of size $n = 16,129$.

We are interested in computing a subset of eigenvalues closest to the shift $\sigma$ using the block PLHR
iteration. Since $L$ is SPD, we employ the real arithmetic version of BPLHR in Algorithm~\ref{alg:bplhr_real}. 
In contrast to the previous example, where an 
artificial preconditioner has been constructed, we now utilize the practical SPD preconditioner 
$T \approx |L - \sigma I|^{-1}$ introduced in~\cite{Ve.Kn:13}.
% for solving shifted linear systems $(L - \sigma I) x = b$. 
For convenience, we state this preconditioning procedure in Algorithm~\ref{alg:avp-mg} 
of Appendix A.

The objective of the current experiment is two-fold. On the one hand, we would like to compare BPLHR to a well-established
solution scheme, such as the BGD method. On the other hand, similar to our previous test, 
we are interested in demonstrating effects of the $T$-harmonic extraction and AV preconditioning 
on the eigensolver convergence.

We consider two preconditioning options for the BGD algorithm. The first 
approach is exactly the same MG AV preconditioner~\cite{Ve.Kn:13} 
(Algorithm~\ref{alg:avp-mg} of Appendix A) as the one used in the BPLHR algorithm, i.e., $T \approx |L - \sigma I|^{-1}$
and it is SPD.
%Below, this scheme is referred to as BGD-AV.
%
The second preconditioner 
%approach represents an indefinite preconditioner 
is given by a standard MG solve~\cite{Briggs.Henson.McCormick:00, Trottenberg.Oosterlee.Schuller:01} 
for the shifted matrix $L - \sigma I$ (see Algorithm~\ref{alg:inv-mg} of Appendix A), 
which corresponds to an indefinite ``shift-and-invert'' type preconditioner $T \approx (L - \sigma I)^{-1}$.  
%We call this scheme BGD-INV.

The preconditioners in Algorithms~\ref{alg:avp-mg} and~\ref{alg:inv-mg}
%, used in BGD-AV and BGD-INV, 
%respectively, 
are of the same nature. They represent a standard MG V-cycle, with the difference that the former is applied to
solve the SPD system $|L- \sigma I|w  = r$~\cite{Ve.Kn:13}, whereas the latter seeks to approximate the solution of 
the indefinite $(L- \sigma I)w  = r$. 
%%%%%%%%%%%%%%%%%%%REFEREE
Note that preconditioners stronger than Algorithm~\ref{alg:inv-mg} are available for the indefinite matrix 
$L- \sigma I$, such as, e.g., in~\cite{Livshits.Brandt:06}. However, due to the algorithmic similarity 
to the employed MG AV preconditioner, for demonstration purposes, we use Algorithm~\ref{alg:inv-mg}
as a reference indefinite MG preconditioner.     

%Note that the preconditioner in Algorithm~\ref{alg:avp-mg} uses Laplacians 
%to approximate the absolute value operators on the intermediate grids, and the construction of $|L - \sigma I|^{-1}$ 
%appears only on the coarsest level.

%%%%%%%%%%%%%%%%%%%
%
%In essence, the AV preconditioner~\cite{Ve.Kn:13}, referred to as AV-MG, represents a standard MG $V$-cycle,
%where the $l$-th level operators $|L_l - \sigma I_l|$ are properly approximated using Laplacians, and the 
%exact $|L_l - \sigma I_l|^{-1}$ is computed only on the coarsest grid. In our tests, we set the size of the coarsest grid
%problem to $225$. 
%The whole preconditioning procedure is summarized in Algorithm~\ref{alg:avp-mg} listed in Appendix. 
%
%We compare the BPLHR algorithm with the block version of the GD (BGD) algorithm. 
%We consider two preconditioning options for BGD. The first option is to use the same AV-MG preconditioner
%as in BPLHR, so that $T \approx |L - \sigma I|^{-1}$ and is SPD. 
%Another option is to use a more traditional 
%MG preconditioner that aims at inverting $L - \sigma I$ rather than its absolute value. 
%In this case, $T \approx (L - \sigma I)^{-1}$ and is indefinite. The preconditioning scheme, which we call INV-MG, 
%is stated in Algorithm~\ref{alg:av-inv} in Appendix. 
 
In order to ensure a comparable 
(in terms of the approximate solves for the corresponding linear systems $|L- \sigma I|w  = r$
and $(L- \sigma I)w  = r$) 
preconditioning quality in the two variants of BGD, we require that the coarsest grid problems are of the same size 
(225 by 225) and that the same smoothing schemes (a single step of Richardson's iterations) 
are used on every level. Note that the computational costs of both preconditioners is essentially the same,
with the AV preconditioner performing slightly more arithmetic operations because of the polynomial 
approximations of the absolute value operators on intermediate 
levels. However, these additional expenses are negligible relative to the overall preconditioning cost.   
%Below, we refer to the BGD algorithm with the AV-MG preconditioner as BGD-AV. 
%If INV-MG is used, we call the corresponding algorithm BGD-INV.

Recall that at every iteration the BGD algorithm expands the search subspace with a set of 
preconditioned residuals. Thus, an increased amount of memory is required by the method 
at every new step. This is in contrast to the BPLHR algorithm, where the requested storage size 
is fixed at every iteration. In particular, in our implementation, BPLHR has to
store at most $12k$ vectors corresponding to 
$Z (= [V^{(i)}, W^{(i)}, S^{(i)}, P^{(i)}])$, $AZ$, and $T(A - \sigma B) Z$.

To maintain the same memory requirement in BGD, we restart the method once the size $m$ of
its search subspace becomes sufficiently large, reaching some prespecified $m_{\max}$. 
In this case, we collapse the search subspace, so that it only contains $k$ available eigenvector approximations.  
Since in standard BGD implementations each iteration of the method stores
the search subspace $Z_m$ together with the block $AZ_m$ (i.e., up to $2m$ vectors total), 
we want $2m$ not to exceed $12k$. Hence, to ensure the same memory requirement for BPLHR and BGD, 
we set $m_{\max} = 6k$. This maximum size of the BGD subspace is somewhat larger than the size of the 
trial subspaces in BPLHR which is $4k$. Nevertheless, as we demonstrate below, BPLHR can be more robust,
even though the extraction is performed with respect to the smaller subspaces.

While our main focus is on the comparison of the BPLHR and BGD methods, we are also interested in the question
of how much the AV preconditioning affects the eigensolver's convergence and if any benefit is received from
the $T$-harmonic RR. For this reason, along with the original BPLHR version in Algorithm~\ref{alg:bplhr_real}, 
%similar to the previous experiment,
we also consider its versions based on the standard harmonic RR with the AV and
indefinite preconditioning options, similar to the previous section.  
%-HARM-AV and PLHR-HARM-INV schemes (see Table~\ref{tab:plmrs}), 
Here, the AV and indefinite preconditioning strategies are based on the MG schemes in 
Algorithms~\ref{alg:avp-mg} and~\ref{alg:inv-mg}, respectively. 
%Here, both variants are obtained
%from the AV-MG preconditioned BPLHR by replacing the $T$-harmonic RR with the standard harmonic procedure. 
%In PLHR-harm-INV, additionally, the SPD AV-MG preconditioner is replaced by the indefinite AV-INV. 

%\begin{comment}
%{\small
%\begin{table}
%\centering
%\begin{tabular}{|c|c|c|c|c|c|c|c|}
%\cline{2-8} 
%\multicolumn{1}{c|}{} & \multicolumn{7}{c|}{Shifts ($\sigma$)}\tabularnewline
%\hline 
%Eigensolvers & 400 & 450 & 500 & 550 & 600 & 650 & 700 \tabularnewline
%\hline 
%\hline 
%BPLHR & 57 & 81 & 68 & 133 & 117 & 190 & 278 \tabularnewline
%\hline 
%BPLHR-HARM-AV & 563 & - & 493 & 635 & - & - & - \tabularnewline
%\hline 
%BPLHR-HARM-INV & 30 & 40 & 45 & - & 59 & 338 & 424 \tabularnewline
%\hline 
%BGD-AV & 209 & - & 533 & - & 493 & - & -  \tabularnewline
%\hline 
%BGD-INV & 36 & 46 & 57 & - & 376 & 763 & - \tabularnewline
%\hline 
%\end{tabular}
%\caption{Iteration numbers required by different eigensolvers to converge to 10
%eigenpairs of the 2D FD Laplacian closest to the shifts $\sigma$;
%``-'' corresponds to the cases where no
%convergence was reached within $1000$ iteration; $n = 16,129$.}\label{tab:laplace10}
%\end{table}
%}
%\end{comment}

{\small
\begin{table}
\centering
\begin{tabular}{|c|c|c|c|c|c|c|c|c|c|}
\cline{4-10} 
\multicolumn{3}{c|}{} & \multicolumn{7}{c|}{Shifts ($\sigma$)}\tabularnewline
\hline 
Iter. scheme & Prec. & RR & 400 & 450 & 500 & 550 & 600 & 650 & 700 \tabularnewline
\hline 
\hline 
BPLHR & AV & $T$-harm. & 57 & 81 & 68 & 133 & 117 & 190 & 278 \tabularnewline
\hline 
BPLHR & AV & harm. & 563 & - & 493 & 635 & - & - & - \tabularnewline
\hline 
BPLHR & Indef. & harm. & 30 & 40 & 45 & - & 59 & 338 & 424 \tabularnewline
\hline 
BGD & AV & harm. & 209 & - & 533 & - & 493 & - & -  \tabularnewline
\hline 
BGD & Indef. & harm. & 36 & 46 & 57 & - & 376 & 763 & - \tabularnewline
\hline 
\end{tabular}
\caption{Iteration numbers required by different eigensolvers to converge to 10
eigenpairs of the 2D FD Laplacian closest to the shifts $\sigma$;
``-'' corresponds to the cases where no
convergence was reached within $1000$ iteration; $n = 16,129$.
The AV and indefinite preconditioners are given by Algorithms~\ref{alg:avp-mg} and~\ref{alg:inv-mg},
respectively. 
}\label{tab:laplace10}
\end{table}
}

In Table~\ref{tab:laplace10} we report the numbers of iterations required 
by different eigensolvers to converge to 10 eigenpairs closest to the given shift. 
The schemes
were compared for a number of shifts in the range from $400$ to $700$. 
The convergence tolerance for the residual norms was set to $10^{-6}$ and the same 
initial guess was used for each run corresponding to the same $\sigma$.
Since the real arithmetic version of BPLHR is invoked, 
according to the discussion in Section~\ref{subsec:plmr_real}, we increase the block size by one,
i.e., apply Algorithm~\ref{alg:bplhr_real} with $k = 11$, but track the convergence only of the ten wanted eigenpairs.
The maximum size $m_{\max}$ of the BGD search subspace is $66$.  

Table~\ref{tab:laplace10} shows that the BPLHR algorithm (i.e., the original version with the AV 
preconditioner and $T$-harmonic RR) is robust with respect to the
choice of the shift. It is the only method among the compared schemes that was able to converge
all eigenpairs for each $\sigma$ in the prescribed range. Note that the $T$-harmonic extraction is 
crucial---its replacement by the standard harmonic approach (while preserving the same AV preconditioner) resulted in a
significant increase of the iteration count or a total loss of convergence. 
%for a number of shift values. 
The demonstrated results also suggest that the standard harmonic extraction leads 
to more satisfactory results if an indefinite preconditioner is employed. However,
this combination was still unable to maintain convergence for all shifts and required
a noticeably larger amount of iterations for $\sigma \geq 650$.  
    
Regardless of the choice of the preconditioner, both BGD based schemes fail to converge 
for a number of shift values. Note that for smaller values of $\sigma$ ($400$ to $500$), 
BGD with the indefinite preconditioning gives the lowest number of iterations. However, as $\sigma$ increases, either 
the iteration count grows dramatically or the convergence of the method is lost.   

%%%%%%%%REFEREE
As has been discussed in Section~\ref{sec:bplhr}, in the reported runs, the cost of each 
BPLHR iteration is dominated by 2 matrix-block multiplications and 4 block preconditioning operations.
This is clearly more expensive than, e.g., in the BGD method, where only one of each is needed.
%at every step. 
However, as seen in Table~\ref{tab:laplace10}, BGD fails to maintain convergence under the given memory constraint,
whereas BPLHR succeeds, i.e., the increased iteration cost results in an improved robustness of the overall computation.
%though at a price of a higher iteration cost.  
%trading a higer iteration cost for robustness.

It is well known (e.g.,~\cite{Elman.Ernst.OLeary:01, Ve.Kn:13}) that the the increase of $\sigma$ generally
leads to the deterioration of the MG solves for $|L - \sigma I|$ and $L - \sigma I$ in 
Algorithms~\ref{alg:avp-mg} and~\ref{alg:inv-mg}. Therefore, the observed shift 
robustness of BPLHR indicates that the method is more stable with respect to the loss of preconditioning
quality compared to the other methods tested. Remarkably, as $\sigma$ increases, the number of BPLHR iterations
does not grow too fast.
%, which is the case, e.g., for BGD-INV.    

%{\small
%\begin{table}
%\centering
%\begin{tabular}{|c|c|c|c|c|c|c|c|}
%\cline{2-8} 
%\multicolumn{1}{c|}{} & \multicolumn{7}{c|}{Shifts ($\sigma$)}\tabularnewline
%\hline 
%Eigensolvers & 800 & 900 & 1000 & 1100 & 1200 & 1300 & 1400 \tabularnewline
%\hline 
%\hline 
%BPLHR & 270 & 168 & 177 & 344 &  365 & 363 &  192\tabularnewline
%\hline 
%BPLHR-HARM-AV & 590 & 417 & 377 & 625 & 437 &  217 & 287 \tabularnewline
%\hline 
%BPLHR-HARM-INV & - & - & - & - & - & - &  -\tabularnewline
%\hline 
%BGD-AV & 331 & 305 & 356 & 666 & 509 & 481 &   443\tabularnewline
%\hline 
%BGD-INV & 230 & 818 & 837 & - &  - &  - &  -\tabularnewline
%\hline 
%\end{tabular}
%\caption{Iteration numbers required by different eigensolvers to converge to 20
%eigenpairs of the 2D FD Laplacian closest to the shifts $\sigma$;
%``-'' corresponds to the cases where no
%convergence was reached within $1000$ iteration; $n = 16,129$.}\label{tab:laplace20}
%\end{table}
%}

{\small
\begin{table}
\centering
%\begin{tabular}{|c|c|c|c|c|c|c|c|}
%\cline{2-8} 
%\multicolumn{1}{c|}{} & \multicolumn{7}{c|}{Shifts ($\sigma$)}\tabularnewline
\begin{tabular}{|c|c|c|c|c|c|c|c|c|c|}
\cline{4-10} 
\multicolumn{3}{c|}{} & \multicolumn{7}{c|}{Shifts ($\sigma$)}\tabularnewline
\hline 
Iter. scheme & Prec.& RR & 800 & 900 & 1000 & 1100 & 1200 & 1300 & 1400 \tabularnewline
\hline 
\hline 
BPLHR & AV & $T$-harm. & 270 & 168 & 177 & 344 &  365 & 363 &  192\tabularnewline
\hline 
BPLHR & AV & harm. & 590 & 417 & 377 & 625 & 437 &  217 & 287 \tabularnewline
\hline 
BPLHR & Indef. & harm. & - & - & - & - & - & - &  -\tabularnewline
\hline 
BGD & AV & harm. & 331 & 305 & 356 & 666 & 509 & 481 &   443\tabularnewline
\hline 
BGD & Indef. & harm. & 230 & 818 & 837 & - &  - &  - &  -\tabularnewline
\hline 
\end{tabular}
\caption{Iteration numbers required by different eigensolvers to converge to 20
eigenpairs of the 2D FD Laplacian closest to the shifts $\sigma$;
``-'' corresponds to the cases where no
convergence was reached within $1000$ iteration; $n = 16,129$.
The AV and indefinite preconditioners are given by Algorithms~\ref{alg:avp-mg} and~\ref{alg:inv-mg},
respectively. 
}\label{tab:laplace20}
\end{table}
}

In Table~\ref{tab:laplace20} we report a similar experiment, where the number $k$ of targeted eigenpairs has been 
increased to $20$; the maximum size of the BGD search subspace is set to $126$. The range of shifts is between
800 and 1400.
Again, we can see that BPLHR was able to converge for all values of $\sigma$ and in most 
cases exhibited the lowest iteration count. The schemes based on BPLHR and BGD with the indefinite MG preconditioning
%,BPLHR-HARM-INV and BGD-INV, 
failed to provide satisfactory convergence. We relate it to a 
deteriorated quality of the INV-MG preconditioner in Algorithm~\ref{alg:inv-mg} for larger shift values. In particular, this 
shows that the AV preconditioning is more robust if the targeted eigenpairs are 
deeper in the spectrum's interior. Also, note that the increased block size, 
compared to the case in Table~\ref{tab:laplace10}, allows PLHR to handle larger 
values of $\sigma$, due to the corresponding increase of the size of the trial subspaces.    
%{\color{blue} finish this table, add a description paragraph here} 
% !!! MENTION 1 buffer vector ==> m_max = (10+1)*6 = 66

{\small
\begin{table}
\centering
\begin{tabular}{|c||c|c|c|c|}
\hline 
$\omega$        & 6  & 7  & 8  & 9  \tabularnewline
\hline 
$\#$ iterations & 41 & 42 & 43 & 42   \tabularnewline 
\hline 
\end{tabular}
\caption{Mesh independence of the BPLHR algorithm with the MG AV preconditioner~\cite{Ve.Kn:13} for the 2D FD Laplacian. 
The mesh parameter is given by $h = 1/(2^{\omega} +1 )$. Four eigenpairs corresponding to the 
eigenvalues closest to $\sigma = 400$ are computed.  
}\label{tab:meshindep}
\end{table}
}

It was demonstrated in~\cite{Ve.Kn:13} that, in the context of solving linear systems $(L - \sigma I)w = r$, 
the MG AV preconditioner in Algorithm~\ref{alg:avp-mg} leads to a mesh-independent convergence of an iterative solver.
In Table~\ref{tab:meshindep}, we show that this property also holds if the same AV preconditioning is combined with the 
BPLHR scheme for computing interior eigenpairs. In particular, we decrease the mesh size $h = 1/(2^{\omega}+1)$ by varying the 
parameter $\omega$ between $6$ and $9$, and observe that the number of steps required to obtain the residual norm of $10^{-4}$ 
is about the same for each run. Here, in order to mitigate the effects of deflation, we compute fewer ($k = 4$) eigenpairs. 
The shift $\sigma$ is set to $400$.

\subsection{Interior eigenpairs of the Kohn--Sham Hamiltonians}
\label{subsec:kssolv}

In this concluding set of experiments we apply BPLHR to several Hermitian matrices that arise
in the context of electronic structure calculations. These matrices correspond
to plane wave discretizations of the Hamiltonian operators in the framework
of the Kohn--Sham (KS) density functional theory~\cite{Kohn.Sham:65}.
%
%All matrices in these examples are generated using the KSSOLV {\sc matlab} 
%toolbox~\cite{kssolv:09} for solving Kohn--Sham equations.
All tests are performed within the KSSOLV package~\cite{kssolv:09}---a {\sc matlab} toolbox
for solving the KS equations.   

\begin{figure}[!ht]
\centering
\begin{subfigure}{.33\textwidth}\includegraphics[width=\textwidth]{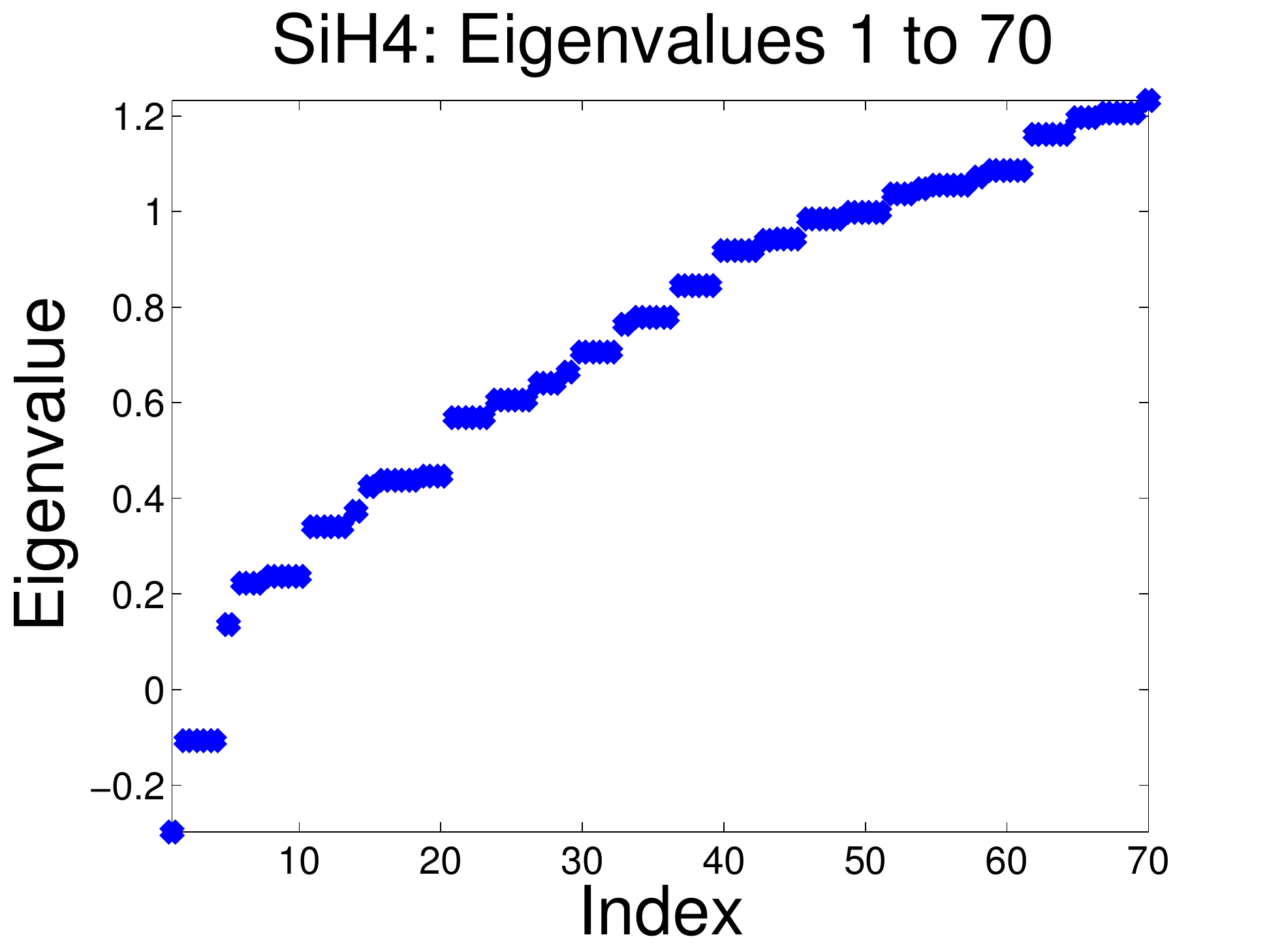}
\caption{SiH4.}\end{subfigure}\hfill
\begin{subfigure}{.33\textwidth}\includegraphics[width=\textwidth]{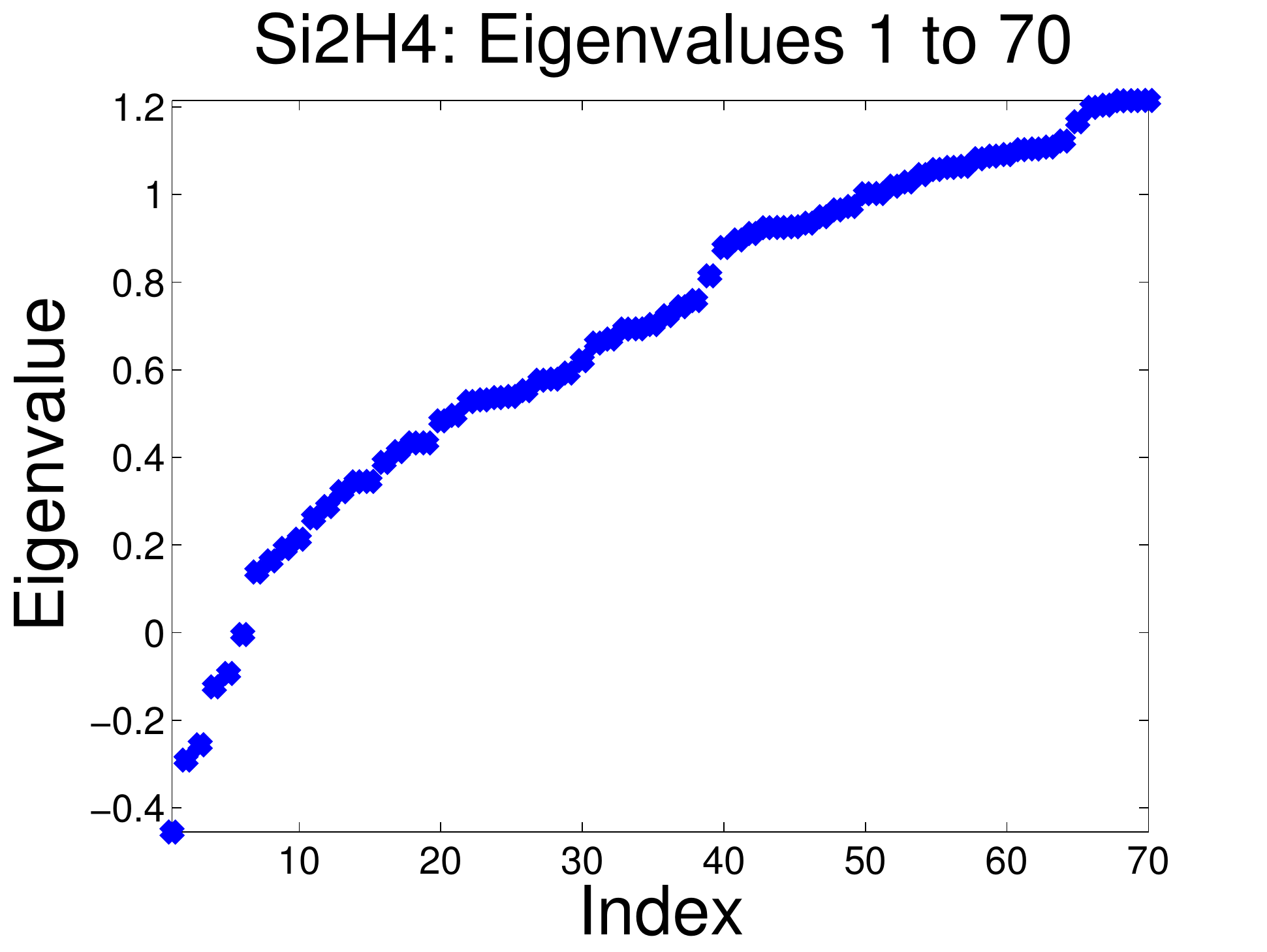}
\caption{Si2H4.}\end{subfigure}\hfill
\begin{subfigure}{.33\textwidth}\includegraphics[width=\textwidth]{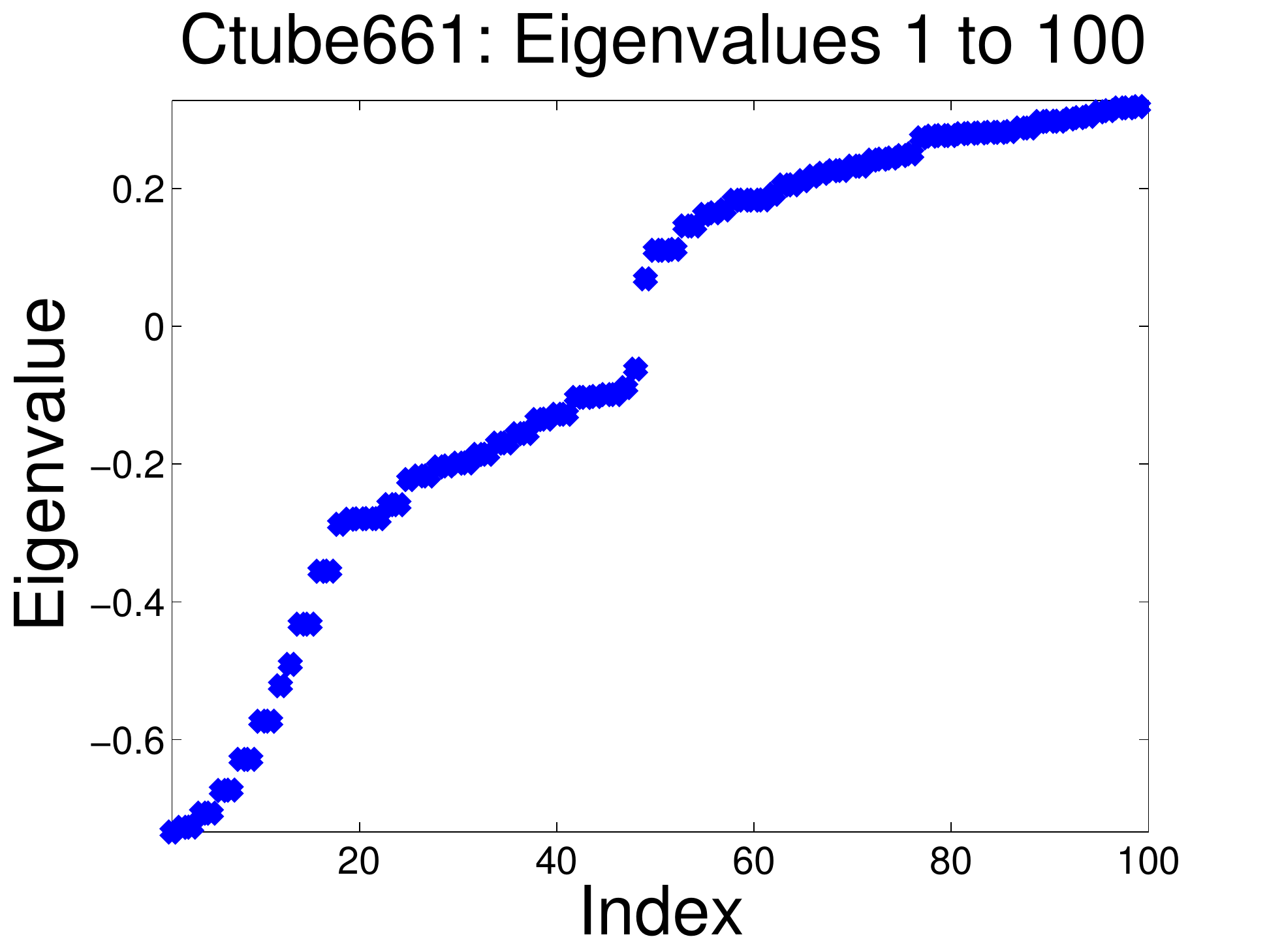}
\caption{Ctube661.}\end{subfigure}
\caption{Parts of spectra of the Hamiltonian matrices.}
\label{fig:spectr}
\end{figure}

We consider three model systems of particles: the silane (SiH2) and planar singlet silysilylene 
(Si2H4) molecules, and a carbon nanotube (Ctube661). For each system, we set up the corresponding KS equations,
which give a nonlinear eigenvalue problem. 
The problem is discretized 
%with the plane wave basis 
and then solved using a self-consistent field iteration. 
As a result of this iteration, the nonlinear operator (Hamiltonian) of the discrete KS problem converges 
to a matrix which represents the Hamiltonian for the converged electron density.
The spectrum of this converged Hamiltonian describes the electronic structure of the underlying system. Its eigenvalues 
represent different energy levels and the eigenvectors define the associated wavefunctions, or orbitals.   
%In our tests, we compute eigenpairs of the converged Hamiltonian 
%matrices that are closest to some given shifts.  

%%In our examples, we are interested in computing $k = 10$ eigenpairs closest to a given shift.
%In our experiments, we are interested in computing $10$ eigenpairs, closest to given shifts,
%of three model Hamiltonian matrices. These Hamiltonians correspond to the silane (SiH2) molecule, 
%planar singlet silysilylene (Si2H4) molecule, and a carbon nanotube (Ctube661). The three matrices are generated
%by KSSOLV, where a self-consistent field iteration has been initially 
%applied to obtain the converged electron density, i.e., the converged Hamiltonians are used in our examples.

In our tests, we are interested in computing interior eigenpairs
of the converged Hamiltonians corresponding to the three model systems.
In particular, for each problem, we would like to compare the convergence of the BPLHR algorithm to that of BGD,
where the schemes are applied to find $k$ eigenpairs around a given shift $\sigma$. 
Throughout, we let $k$ be equal to $10$, i.e., ten eigenpairs closest to $\sigma$ are sought in every test.  
Since the Hamiltonians are complex, all our computations are performed in the complex arithmetic, i.e.,
Algorithm~\ref{alg:bplmr} is employed. 

We use the Teter--Payne--Allan preconditioner developed in~\cite{Teter.Payne.Allan:89}. 
This preconditioning approach is known to be effective for eigenvalue computations in the context of the
plane wave electronic structure analysis, and is readily available in KSSOLV.
Although the preconditioner is more commonly used for computing a number of lowest eigenpairs, 
it can also be applied to the interior eigenvalue computations provided that the targeted eigenvalues
are not too deep inside the spectrum. For example, combining the preconditioning of~\cite{Teter.Payne.Allan:89}
with a (block) Davidson algorithm based on the harmonic extraction  
was suggested for computing interior eigenpairs in~\cite{Jordan.Marsman.Kim.Kresse:12}. 
%Note that 
The BGD scheme used in this section 
%compared with BPLHR in this paper, 
is similar to this approach.
%the Davidson's approach in~\cite{Jordan.Marsman.Kim.Kresse:12}.   

%However, the same preconditioning can be used in the 
%context of the interior eigenvalue computations provided that $\sigma$ targets the eigenvalues that
%are not too deep inside the spectrum; see, e.g.,~\cite{Jordan.Marsman.Kim.Kresse:12} where
%preconditioner~\cite{Teter.Payne.Allan:89} was used within a (block) Davidson algorithm.  
%Below, we compare BPLHR to a similar scheme that is based on the BGD iteration in Algorithm~?.
%In both cases, we use the preconditioner~\cite{Teter.Payne.Allan:89}.
%%, which is generated by KSSOLV.
%
%The preconditioner represents a diagonal matrix with positive entries.

The preconditioner in~\cite{Teter.Payne.Allan:89} represents a diagonal matrix with positive entries. 
Hence, the preconditioning is SPD and its application is extremely fast.
This is especially beneficial for the BPLHR algorithm, where additional preconditioning operations are needed 
to accomplish the $T$-harmonic extraction.
Since the cost of the diagonal preconditioning is negligible, the total cost of each BPLHR iteration is dominated by two 
matrix-block multiplications. The similar consideration applies to the BGD algorithm, whose iteration cost 
is dominated by a single matrix-block product.

Following the discussion in the preceding subsection, we choose the restart parameter in BGD
to be at least $6k = 60$. In this case BGD and BPLHR have the same memory requirement.
In some of our tests, however, we will allow BGD to construct search subspaces that
are larger than $60$. For this reason, in order to distinguish between different subspace sizes, 
let us denote each BGD run by ``BGD($m_{\max}$)'', where $m_{\max}$ specifies the corresponding 
restart parameter. 

As previously, in order to address the effects of the $T$-harmonic extraction, we also consider the 
BPLHR variant with a standard harmonic RR procedure. 
%We refer to this scheme as BPLHR-HARM. 
It is combined with exactly the same SPD diagonal preconditioner~\cite{Teter.Payne.Allan:89} as the one
used in BPLHR and BGD.

%Namely, in the context of plane wave electronic structure calculations, 
%the preconditioning overhead
%related to the $T$-harmonic extraction is negligibly low, 
%and the cost of each BPLHR iteration is dominated by two matrix-block multiplications. 

%
To discretize the three model Hamiltonians, we use the energy cut-off of 75 Ry for the SiH4 and Si2H4 systems, 
and 25 Ry for the carbon tube. This leads to the eigenvalue problems of size $n = 11,019$ 
(for SiH4 and Si2H4) and $n = 12,599$ (for Ctube661). The parts of spectrum that we are interested in
for each of the three Hamiltonian matrices are plotted in Figure~\ref{fig:spectr}.

\begin{figure}[ptbh]
\begin{center}%
\begin{tabular}
[c]{cc}%
	\includegraphics[width=6.5cm]{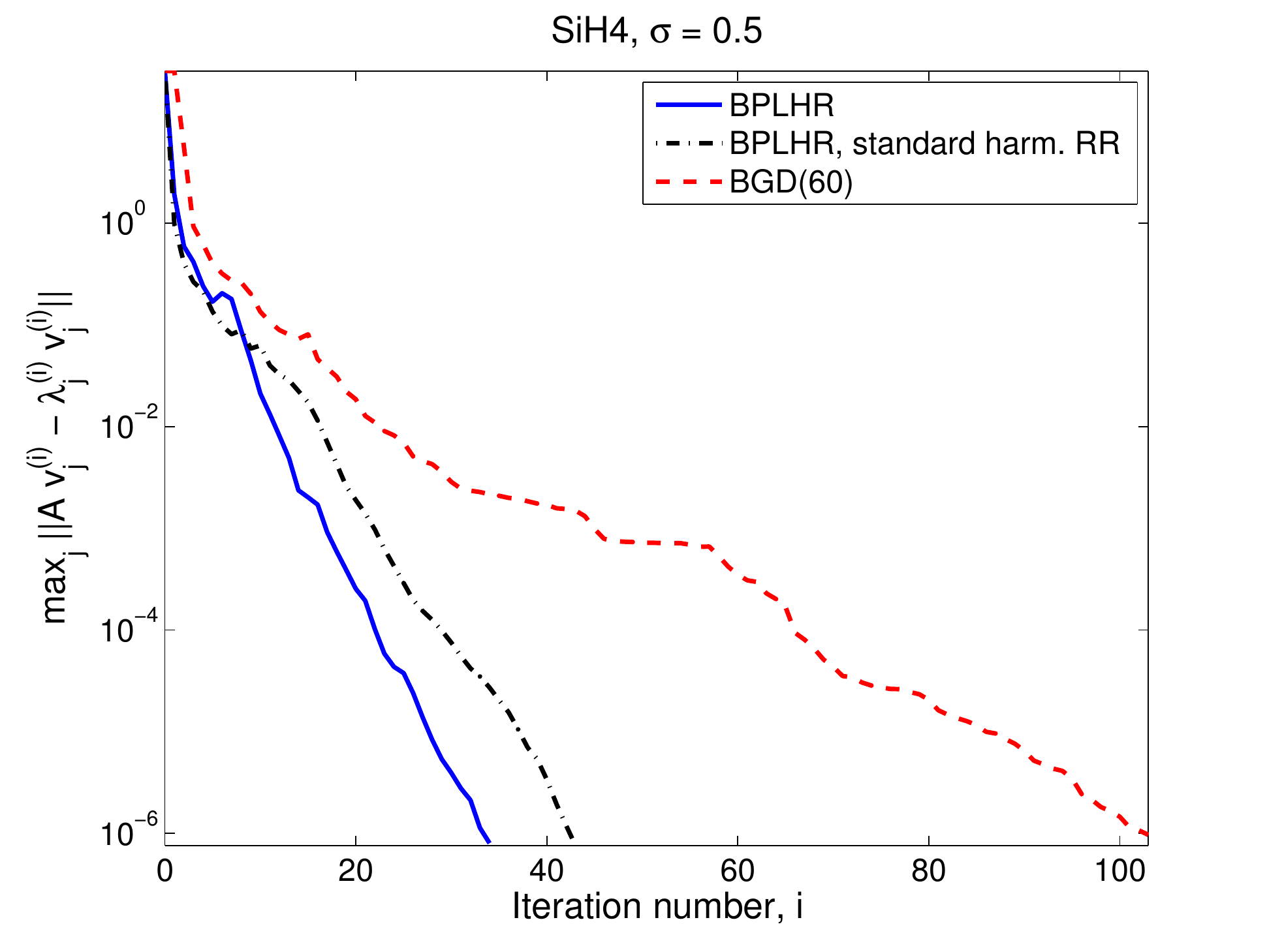}
	\includegraphics[width=6.5cm]{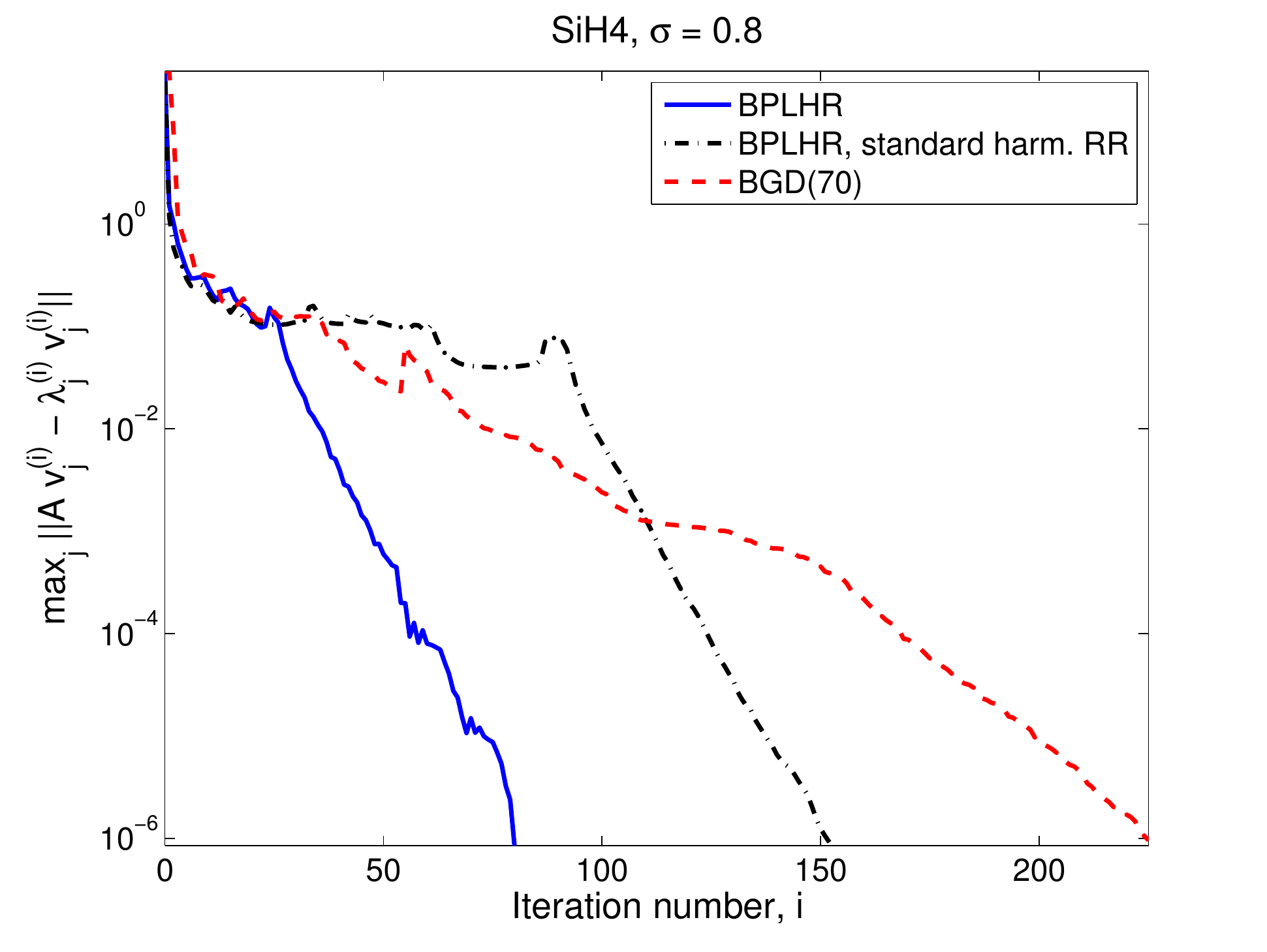}
\end{tabular}
\end{center}
\caption{Computing $10$ eigenpairs closest to $\sigma = 0.5$ (left) and $\sigma = 0.8$ (right) 
 of the converged Hamiltonian of the SiH4 system; $n = 11,019$.}%
\label{fig:sih4}%
\end{figure}

Figure~\ref{fig:sih4} shows the convergence of the three schemes for the SiH4 example.
Here and below, we assess the convergence by monitoring the largest 
residual norm in the block. Note that in our experiments all the $10$ targeted pairs 
converge essentially at the same rate. Therefore, tracking only the largest norm
is indicative of the convergence behavior of each eigenpair in the block. 
 
It can be seen from Figure~\ref{fig:sih4} that for both shifts, $\sigma = 0.5$ (left) and $\sigma = 0.8$ (right),
the BPLHR algorithm results in almost a three times reduction in the number of iterations compared to BGD.    
Even though each BPLHR iteration is (roughly) twice as expensive as the BGD step, 
the overall decrease of the computational work is evident. Additionally, note that the BGD run in
Figure~\ref{fig:sih4} (right) requires more memory than BPLHR, i.e., BPLHR gives a faster
convergence while consuming less storage. Our experiments below will make this observation
yet more pronounced. 
 
We can see from Figure~\ref{fig:sih4} (left) that the introduction of the $T$-harmonic projection
results in a minor improvement of the convergence for $\sigma = 0.5$. In this case, the number of BPLHR iterations is
only slightly decreased compared to its version with the standard harmonic RR. 
However, for the larger shift ($\sigma = 0.8$) in
Figure~\ref{fig:sih4} (right) the situation is substantially different. The $T$-harmonic projection in BPLHR
allows reducing the iteration count by more than a factor of two. Thus, the $T$-harmonic RR procedure
makes the scheme less sensitive to the choice of the shift and more stable with respect 
to deterioration of the preconditioning quality. Note that this is consistent with our observations in the
previous subsections.

\begin{figure}[ptbh]
\begin{center}%
\begin{tabular}
[c]{cc}%
	\includegraphics[width=6.5cm]{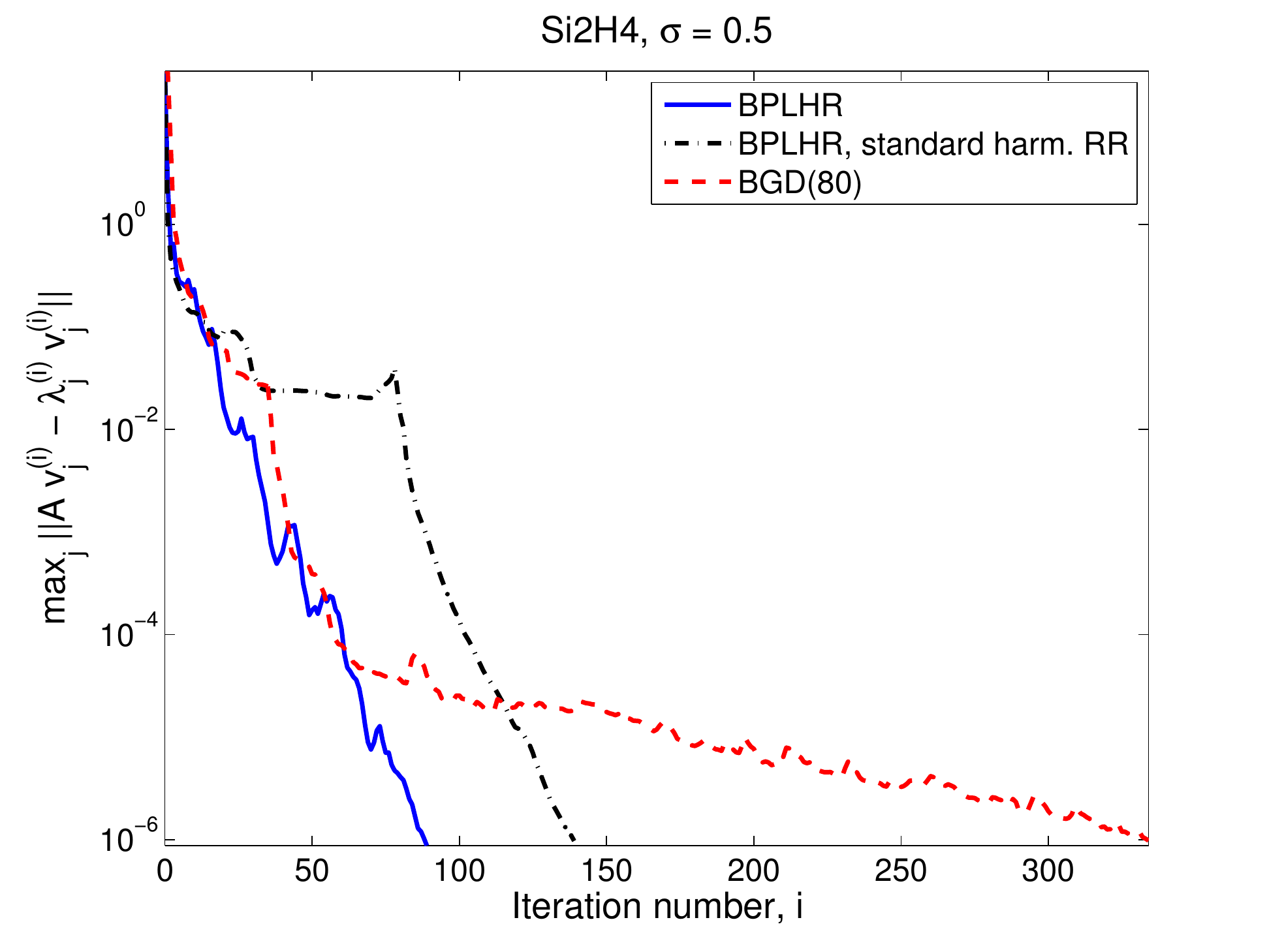}
	\includegraphics[width=6.5cm]{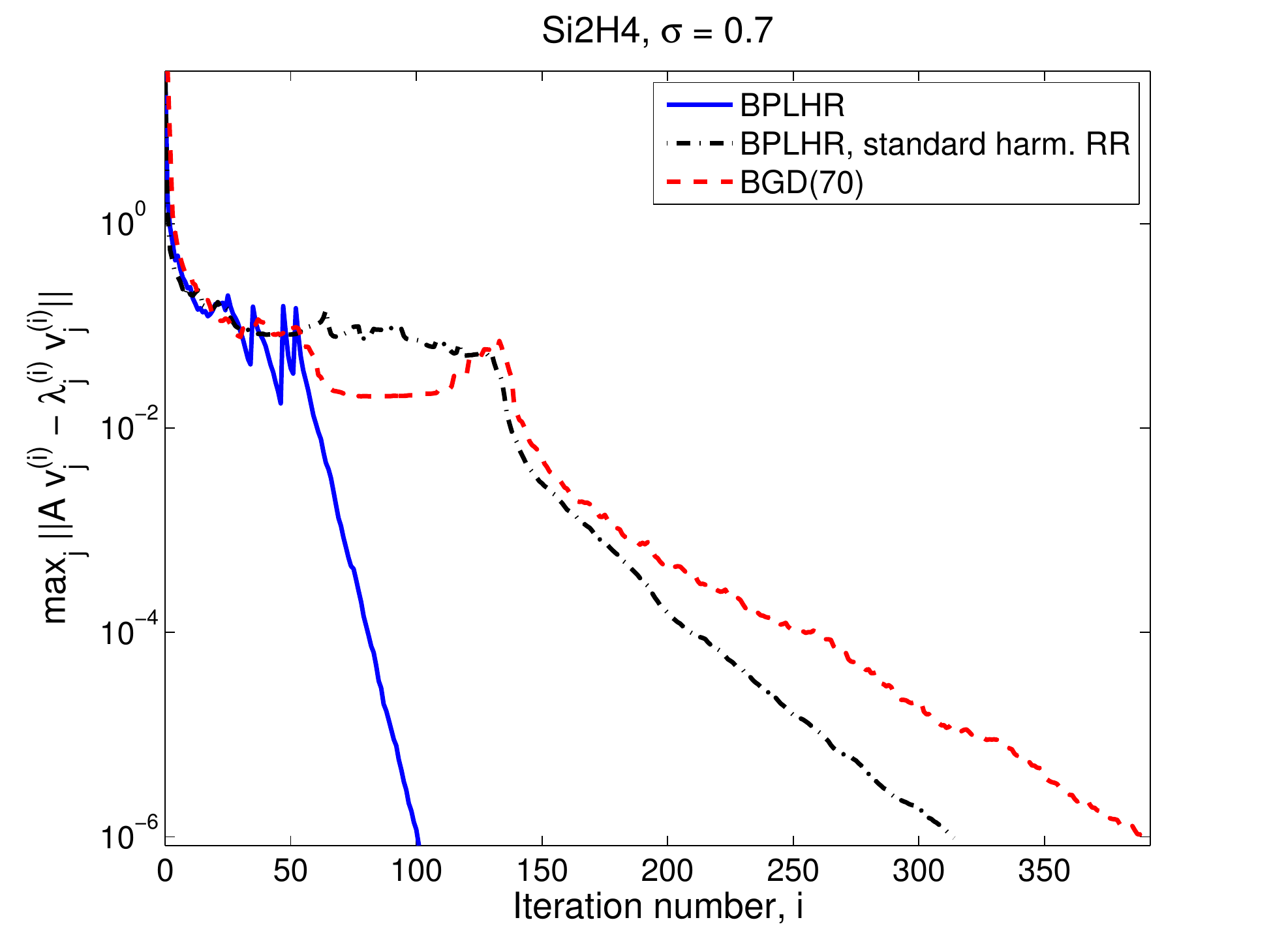}
\end{tabular}
\end{center}
\caption{Computing $10$ eigenpairs closest to $\sigma = 0.5$ (left) and $\sigma = 0.7$ (right) 
 of the converged Hamiltonian of the Si2H4 system; $n = 11,019$.}%
\label{fig:si2h4}%
\end{figure}

In Figure~\ref{fig:si2h4}, we apply the eigensolvers to the Hamiltonian of the Si2H4 system.
The targeted energy shifts are $\sigma = 0.5$ (left) and $\sigma = 0.7$ (right). In both cases, the BPLHR
algorithm gives the smallest iteration count and a significant decrease of the overall computational
work. For example, we can see around 4 time reduction in the number of iterations compared to BGD. 
The impact of the $T$-harmonic projection on the BPLHR convergence can be observed by 
comparing the method with its variant based on the standard harmonic RR. Similar to the previous test, note that
BGD requires more memory than the BPLHR schemes.

\begin{figure}[ptbh]
\begin{center}%
\begin{tabular}
[c]{cc}%
	\includegraphics[width=6.5cm]{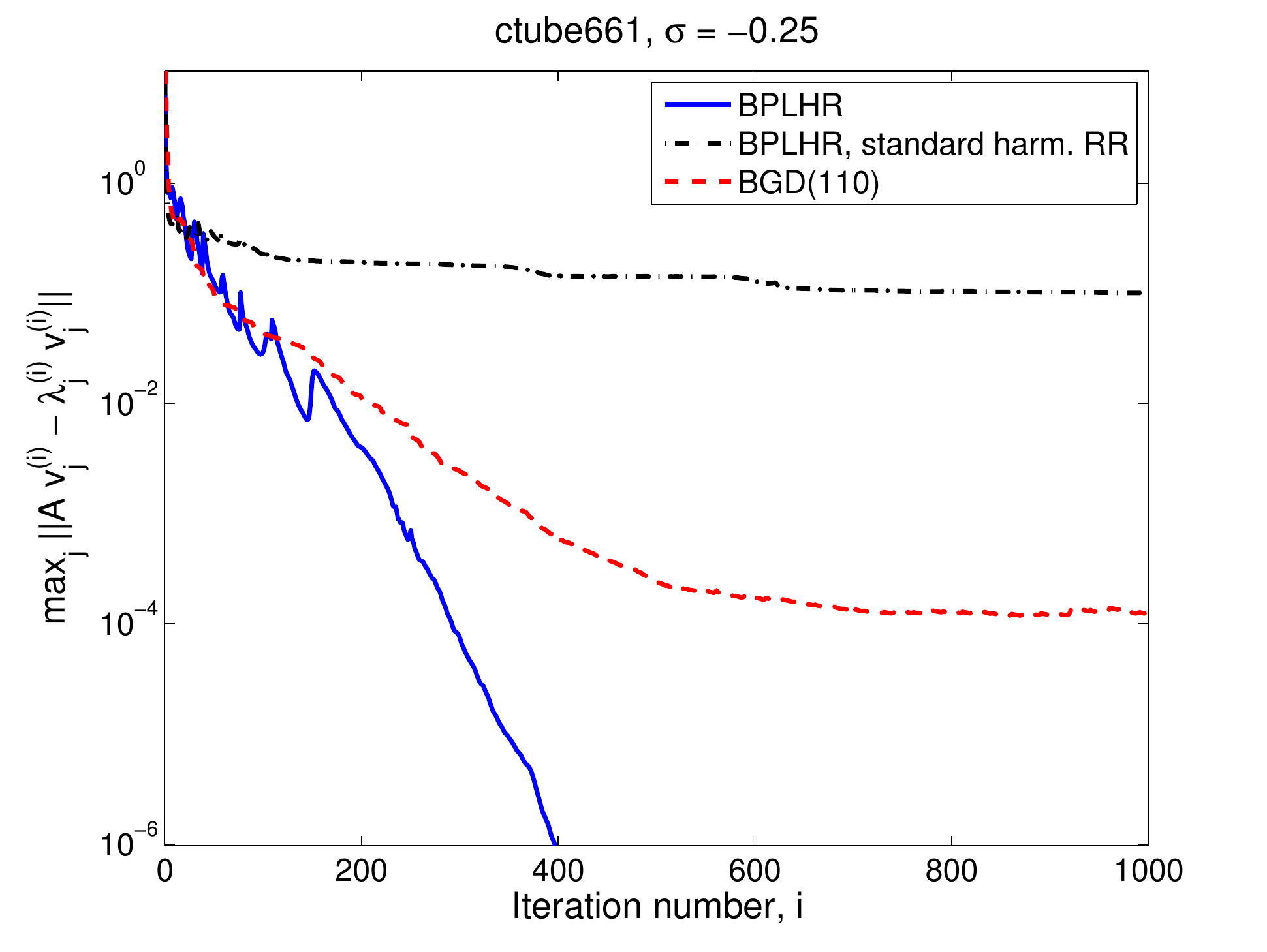}
	\includegraphics[width=6.5cm]{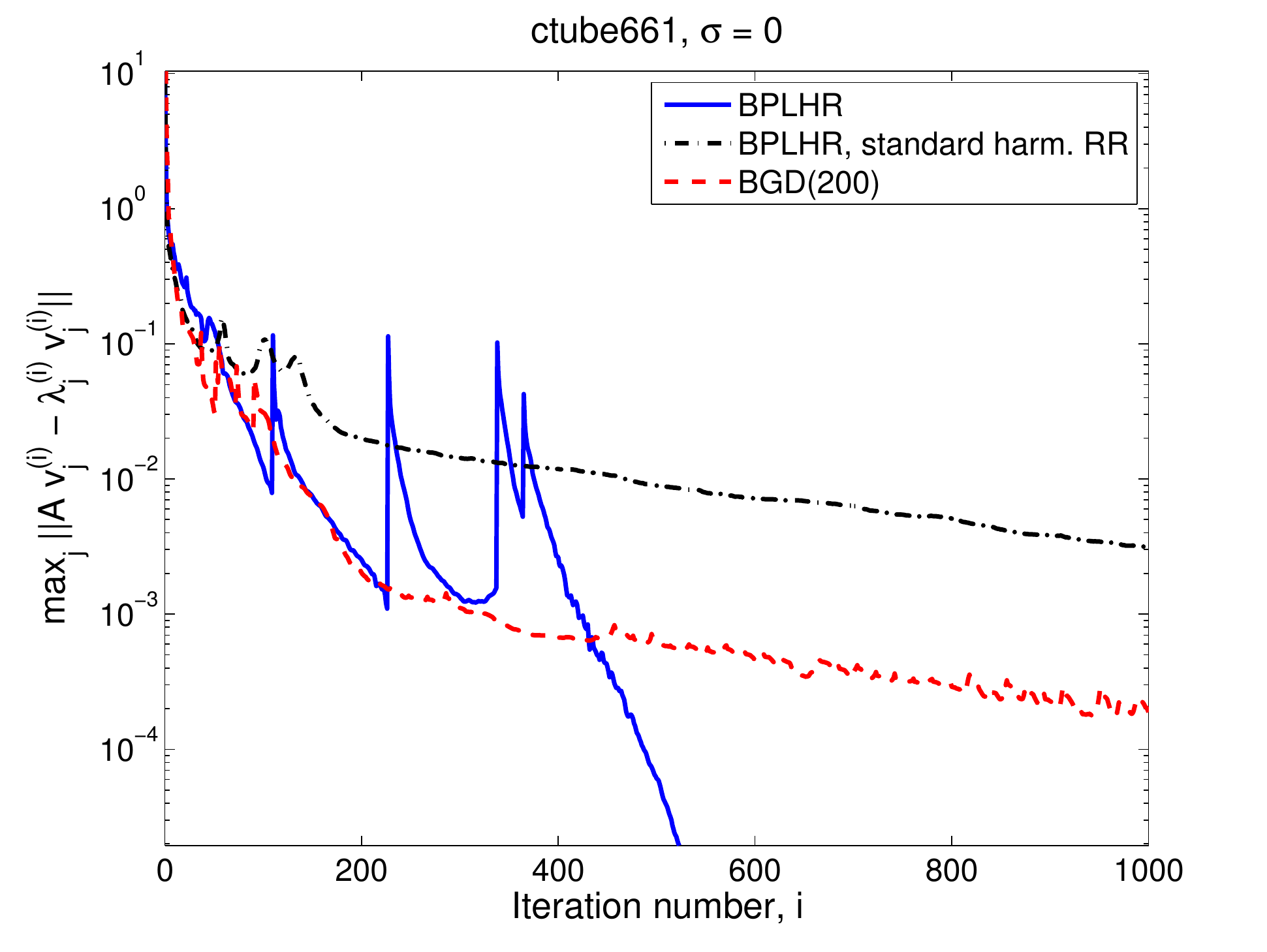}
\end{tabular}
\end{center}
\caption{Commuting $10$ eigenpairs closest to $\sigma = -0.25$ (left) and $\sigma = 0$ (right) 
 of the converged Hamiltonian of the carbon tube; $n = 12,599$.}%
\label{fig:ctube661}%
\end{figure}

Our last experiment for the carbon tube is presented in Figure~\ref{fig:ctube661}. 
In this example, the BPLHR algorithm is the only scheme that is able to reach 
convergence of the eigenpairs near the shift $\sigma = -0.25$ (left) and 
$\sigma = 0$ (right). Remarkably, the BGD search subspaces are allowed to be 3 to 5 times
larger than the BPLHR trial subspaces. Nevertheless, the increased memory consumption
does not allow BGD to maintain the convergence, whereas the BPLHR algorithm computes the
solution with a much tighter storage. The example also clearly demonstrates the importance of 
the $T$-harmonic RR. Substituting the procedure by the standard harmonic RR leads to the
loss of convergence.
%, as demonstrated by the run of BPLHR-HARM.  

\section{Conclusions}\label{sec:concl_interior}

We have presented the Preconditioned Locally Harmonic Residual (PLHR) algorithm for
computing interior eigenpairs closest to the shift~$\sigma$. 
The method represents a preconditioned four-term recurrence and can be easily extended
to the block case (BPLHR). It is equally applicable to the standard and generalized eigenvalue problems.
%and both problems are treated in the same manner. 
The algorithm is based on the 
$T$-harmonic RR procedure, and does not require shift-and-invert or folded spectrum transformations. 

The proposed approach has 
been tested for a number of model problems, including the Laplacian and Hamiltonian 
(arising from the density functional theory based electronic structure calculations) 
matrices. (B)PLHR has been shown to exhibit a lower sensitivity to the preconditioning
quality and has been able to maintain convergence under stringent memory requirements.    

A possible limitation of the method is given by the need to provide an HPD AV preconditioner.
However, as demonstrated in the paper, such preconditioners are available for a number of important  
applications, such as the plane wave electronic structure calculations. In the case where an AV
preconditioner is unavailable, one can use the algorithm version with the standard harmonic RR instead
of the proposed $T$-harmonic scheme. However, we anticipate that the (B)PLHR approach 
would greatly benefit from further progress in developing efficient AV preconditioning~techniques.

%In this chapter we have proposed a novel approach, which we call the PLHR method, 
%for computing an approximation
%of the smallest, in the absolute value, eigenvalue and the corresponding eigenvector
%of a symmetric matrix pencil. The method represents a four-term recurrent iterative
%scheme, with iteration parameters determined by solving small auxiliary eigenvalue problems.
%The method is preconditioned. It requires an SPD preconditioner, which can be constructed
%according to the idea of the absolute value preconditioning described 
%in the context of symmetric indefinite linear systems in the previous chapter.
%In fact, this allows to use the same SPD preconditioners for both symmetric indefinite linear systems and
%the corresponding interior eigenvalue problems.

%We have applied the PLHR method
%to approximate an eigenpair of the two-dimensional discrete negative Laplace operator,
%which corresponds to the eigenvalue, closest to a given shift value. As a preconditioner 
%we have reused the (geometric) MG absolute value preconditioner, constructed for the corresponding 
%linear system (the model problem) in the previous chapter. For a significant number of the initial steps,
%the PLHR method has demonstrated convergence behavior, comparable to that of an \textit{idealized optimal}
%preconditioned eigenvalue~solver.   

\appendix
\section{The AV-MG and INV-MG preconditioners}
The idea behind the AV-MG preconditioner is to apply the formal MG V-cycle
to the system $|L - \sigma I| w = r$, where $L$ is the Laplacian operator.
A hierarchy of grids is introduced, and at each level $l$ the corresponding absolute value operator
$|L_l - \sigma I_l|$ is approximated by some $B_l$. At finer grids, $B_l$ is chosen to be 
simply the Laplacian, i.e., $B_l = L_l$, whereas at coarser
levels polynomial approximations are employed, so that $B_l = p_{m} (L_l - \sigma I_l)$,
where $m$ is a given degree of the polynomial. The (Richardson's) smoothing is performed with respect
to $B_l$. The restriction and prolongation are carried out in a standard way. 
The actual construction of $|L - \sigma I|$ appears only on the coarsest level, where
the coarse grid solve is performed. The whole scheme is summarized in Algorithm~\ref{alg:avp-mg}.
For more detail we refer the reader to~\cite{Ve.Kn:13}.      
%---
\begin{algorithm}[htbp]
\begin{small}
\begin{center}
  \begin{minipage}{5in}
\begin{tabular}{p{0.5in}p{4.5in}}
{\bf Input}:  &  \begin{minipage}[t]{4.0in}
The residual $r_l$; parameters $\delta$, $\nu$, and $m$;
                  \end{minipage} \\
{\bf Output}:  &  \begin{minipage}[t]{4.0in}
             $w_l \approx |L_l - I_l|^{-1} r_l$; 
                  \end{minipage}
\end{tabular}
\begin{algorithmic}[1]

\STATE Set $B_l = L_l$ if $\sqrt{\sigma} h_l < \delta$. Otherwise define $B_l$ as a polynomial approximation
$p_{m} (L_l - \sigma I_l)$ of $| L_l - \sigma I_l |$, where $m$ is the degree of the polynomial. 
\STATE \emph{Presmoothing}. Apply $\nu$ smoothing steps, $\nu \geq 1$:
%\begin{equation}\label{eqn:mg-pre}
\[
w_l^{(i+1)} = w_l^{(i)} + M_l^{-1} (r_l - B_l w_l^{(i)}),  \ i = 0,\ldots,\nu - 1, \ w_l^{(0)} = 0,
\]
%\end{equation}
where $M_l$ defines a smoother on level~$l$. 
Set $w_l^{pre} = w_l^{(\nu)}$. 
	\STATE \emph{Coarse grid correction}. Restrict ($R_{l-1}$) $r_l - B_l w_l^{pre}$ to the grid $l-1$,
        recursively apply AV-MG, and prolongate ($P_l$) back to the fine grid.
         This delivers the coarse grid correction added to $w_l^{pre}$:
%\begin{equation}\label{eqn:mg-cgc-1}
\[
w_{l-1} = \left\{ 
\begin{array}{ll}
\displaystyle  \left|L_0 - \sigma I_0 \right|^{-1} R_0 \left(r_1 - B_1 w_1^{pre}\right), & l = 1, \\
\displaystyle  \mbox{AV-MG}\left(R_{l-1} \left(r_l - B_l w_l^{pre}\right) \right), & l > 1; 
\end{array}
\right.
%\end{equation}
\]
%\begin{equation}\label{eqn:mg-cgc-2}
\[
  w_l^{cgc}  =  w_l^{pre} + P_l w_{l-1}.
\]
%\end{equation}
	\STATE \emph{Postsmoothing}. Apply $\nu$ smoothing steps:
\begin{equation}\label{eqn:mg-post}
w_l^{(i+1)} = w_l^{(i)} + M_l^{-*} (r_l - B_l w_l^{(i)}),  \ i = 0,\ldots,\nu - 1, \ w_l^{(0)} = w_l^{cgc},
\end{equation}
where $M_l$ and $\nu_l$ are the same as in step 2. 
Return $w_l = w_l^{post} =  w_l^{(\nu)}$.
%Set $w_l^{post} = w_l^{(\nu)}$.   
% \end{enumerate}

\end{algorithmic}
\end{minipage}
\end{center}
\end{small}
  \caption{AV-MG($r_l$): the MG AV preconditioner}
  \label{alg:avp-mg}
\end{algorithm}

The INV-MG preconditioner represents a standard $V$-cycle for system $(L - \sigma I)w = r$
and results in an indefinite preconditioner. The preconditioning scheme is stated in Algorithm~\ref{alg:inv-mg}.
Note that, in this paper, the choice of the main MG components, such as smoothers, restriction and prolongation operators,
for Algorithm~\ref{alg:inv-mg} is identical to AV-MG in Algorithm~\ref{alg:avp-mg}.
\begin{algorithm}[htbp]
\begin{small}
\begin{center}
  \begin{minipage}{5in}
\begin{tabular}{p{0.5in}p{4.5in}}
{\bf Input}:  &  \begin{minipage}[t]{4.0in}
The residual $r_l$; $\nu$;
                  \end{minipage} \\
{\bf Output}:  &  \begin{minipage}[t]{4.0in}
             $w_l \approx (L_l - I_l)^{-1} r_l$; 
                  \end{minipage}
\end{tabular}
\begin{algorithmic}[1]

\STATE \emph{Presmoothing}. Apply $\nu$ smoothing steps, $\nu \geq 1$:
%\begin{equation}\label{eqn:mg-pre}
\[
w_l^{(i+1)} = w_l^{(i)} + M_l^{-1} (r_l - (L_l - \sigma I_l) w_l^{(i)}),  \ i = 0,\ldots,\nu - 1, \ w_l^{(0)} = 0,
\]
%\end{equation}
where $M_l$ defines a smoother on level~$l$. 
Set $w_l^{pre} = w_l^{(\nu)}$. 
	\STATE \emph{Coarse grid correction}. Restrict ($R_{l-1}$) $r_l - (L_l - \sigma I_l) w_l^{pre}$ to the grid $l-1$,
        recursively apply INV-MG, and prolongate ($P_l$) back to the fine grid.
         This delivers the coarse grid correction added to $w_l^{pre}$:
%\begin{equation}\label{eqn:mg-cgc-1}
\[
w_{l-1} = \left\{ 
\begin{array}{ll}
\displaystyle  \left(L_0 - \sigma I_0 \right)^{-1} R_0 \left(r_1 - B_1 w_1^{pre}\right), & l = 1, \\
\displaystyle  \mbox{AV-MG}\left(R_{l-1} \left(r_l - (L_l - \sigma I_l) w_l^{pre}\right) \right), & l > 1; 
\end{array}
\right.
%\end{equation}
\]
%\begin{equation}\label{eqn:mg-cgc-2}
\[
  w_l^{cgc}  =  w_l^{pre} + P_l w_{l-1}.
\]
%\end{equation}
	\STATE \emph{Postsmoothing}. Apply $\nu$ smoothing steps:
%\begin{equation}\label{eqn:mg-post}
\[
w_l^{(i+1)} = w_l^{(i)} + M_l^{-*} (r_l - (L_l - \sigma I_l) w_l^{(i)}),  \ i = 0,\ldots,\nu - 1, \ w_l^{(0)} = w_l^{cgc},
\]
%\end{equation}
where $M_l$ and $\nu_l$ are the same as in step 2. 
Return $w_l = w_l^{post} =  w_l^{(\nu)}$.
%Set $w_l^{post} = w_l^{(\nu)}$.   
% \end{enumerate}

\end{algorithmic}
\end{minipage}
\end{center}
\end{small}
  \caption{INV-MG($r_l$): the MG preconditioner}
  \label{alg:inv-mg}
\end{algorithm}

\bibliographystyle{siam} 
 
\bibliography{plmr}

%\end{thebibliography}
\end{document}